\documentclass[12pt,a4paper,reqno]{amsart}
\usepackage{amssymb,amscd}
\usepackage{pb-diagram,lamsarrow,pb-lams}
\makeatletter
\@addtoreset{equation}{section}
\renewcommand{\@cite}[2]{{\rm[{{\textbf{#1}}\if@tempswa , #2\fi}]}}
\renewcommand{\section}{\@startsection{section}{1}{\parindent}{3.5ex plus 1ex minus .2ex}{2.3ex plus.2ex}{\sc}}
\renewcommand{\subsection}{\@startsection{subsection}{1}{0pt}{-3.25ex plus -1ex minus-.2ex}{1.5ex plus.2ex}{\normalfont\it}}
\renewcommand{\phi}{\varphi}
\renewcommand{\epsilon}{\varepsilon}

\renewcommand{\kappa}{\varkappa}
\DeclareMathOperator{\Hom}{Hom}
\DeclareMathOperator{\Lex}{Lex}
\DeclareMathOperator{\zg}{Zg}
\DeclareMathOperator{\End}{End}
\DeclareMathOperator{\fg}{fg}
\DeclareMathOperator{\fp}{fp}
\DeclareMathOperator{\coh}{coh}
\DeclareMathOperator{\kr}{Ker}
\DeclareMathOperator{\Ext}{Ext}
\DeclareMathOperator{\im}{Im}
\DeclareMathOperator{\coker}{Coker}
\DeclareMathOperator{\md}{mod}
\DeclareMathOperator{\Mod}{Mod}
\DeclareMathOperator{\modd}{mod}
\DeclareMathOperator{\Mor}{Mor}
\newcommand{\ten}[2]{{#1\otimes_A #2}}
\newcommand{\tena}[2]{{#1\otimes_{\cc A} #2}}
\newcommand{\lat}[1]{\text{\sf L}(\cc #1)}
\newcommand{\ol}{\overline}
\newcommand{\tn}[2]{{#1\otimes #2}}

\newcommand{\wt}{\widetilde}
\newcommand{\ifff}{if and only if }
\newcommand{\ccc}[1]{\cc #1^{\cc A}}
\newcommand{\ses}{a short exact sequence }
\newcommand{\rfp}{\modd A}
\newcommand{\lfp}{\modd A^{\text{\rm op}}}
\newcommand{\Lfp}{\Mod A^{\text{\rm op}}}
\newcommand{\Rfp}{\Mod A}
\newcommand{\lgn}[1]{{}_{#1}\cc C}

\newcommand{\mc}[1]{\Mod{\mathcal #1}}
\newcommand{\mcl}[1]{\Mod{\mathcal #1^{\text{\rm op}}}}
\newcommand{\fpl}[1]{\md{\mathcal #1^{\text{\rm op}}}}
\newcommand{\fpr}[1]{\md{\mathcal #1}}

\newcommand{\rcc}{{\cc C}_{\cc A}}
\newcommand{\lcc}{{}_{\cc A}{\cc C}}
\newcommand{\lc}{{}_A{\cc C}}
\newcommand{\rc}{{\cc C}_A}
\newcommand{\rcohh}[1]{#1\in\coh\rcc}
\newcommand{\lcohh}[1]{#1\in\coh\lcc}
\newcommand{\rcoh}[1]{#1\in\coh\rc}

\newcommand{\fc}[1]{{({\cc #1}^{\text{\rm op}},\text{\rm Ab})}}
\newcommand{\mmc}[1]{\Mod{\cc #1}^{\text{\rm op}}}
\newcommand{\ffc}[1]{{(\cc #1,\text{\rm Ab})}}
\newcommand{\fcat}[1]{{({\cc #1},\text{\rm Ab})}}
\newcommand{\ii}[3]{0\longrightarrow#1\longrightarrow#2\longrightarrow#3}
\newcommand{\pp}[3]{#1\longrightarrow#2\longrightarrow#3\longrightarrow 0}
\newcommand{\p}[3]{#1\to#2\to#3\to 0}
\newcommand{\lp}{\varinjlim}
\newcommand{\lo}{\varprojlim}
\newcommand{\ps}{\oplus}

\newcommand{\cc}{\mathcal}
\newcommand{\f}{\mathfrak}
\newcommand{\bb}{\mathbb}
\newcommand{\negl}[1]{$\cc #1$-neg\-li\-gible}
\newcommand{\dec}[3]{{#1{\Pi}_{#2}#3}}
\newcommand{\homa}[2]{\Hom_{\cc A}(#1,#2)}
\newcommand{\hm}[2]{\Hom_A(#1,#2)}
\newcommand{\iso}{\thickapprox}
\newcommand{\rqq}[2]{{}_{\rcc}(#1,#2)}
\newcommand{\lqq}[2]{{}_{\lcc}(#1,#2)}
\newcommand{\qc}[2]{{}_{\cc C}(#1,#2)}
\newcommand{\q}[2]{{}_{\cc C/\cc S}(#1,#2)}
\newcommand{\qq}[2]{#1_{\cc #2}}
\newcommand{\lra}[1]{\bl{#1}\longrightarrow\relax}
\newcommand{\les}[3]{0\longrightarrow#1\longrightarrow#2\longrightarrow#3\longrightarrow 0}
\newcommand{\ls}[4]{0\longrightarrow#1\longrightarrow#2\longrightarrow#3\longrightarrow#4\longrightarrow 0}
\newcommand{\es}[3]{0\to #1\to #2\to #3\to 0}
\newcommand{\bl}[1]{\buildrel #1\over}
\newcommand{\cl}[1]{$\cc #1$-closed}
\newcommand{\tr}[1]{$\cc #1$-torsionfree}
\newcommand{\env}[1]{$\cc #1$-envelope}
\newcommand{\eqv}{if and only if }

\newtheorem{thm}{Theorem}[section]
\newtheorem{prop}[thm]{Proposition}
\newtheorem{cor}[thm]{Corollary}
\newtheorem{lem}[thm]{Lemma}
\newtheorem*{rem}{Remark}
\newtheorem*{thmm}{Theorem}
\newtheorem*{exs}{Examples}
\newtheorem*{example}{Example}
\newtheorem*{defs}{Definitions}
\newtheorem*{question}{Question}
\makeatother

\begin{document}

\title{Grothendieck categories}
\author{Grigory Garkusha}
\address{Saint-Petersburg State University\\ Higher Algebra and Number Theory
Department\\ Faculty of Mathematics and Mechanics\\ Bibliotechnaya Sq. 2\\
198904\\ Russia}
\email{ggarkusha@hotmail.com}
\maketitle

\setcounter{section}{-1}
\tableofcontents

\thispagestyle{empty}
\section{Introduction}

The theory has its origin in the work of Grothendieck~\cite{Gro} who
introduced the following notation of properties of abelian
categories:

\begin{itemize}
\item[Ab3.] An abelian category with coproducts or equivalently, a
	    cocomplete abelian category.
\item[Ab5.] Ab3-category, in which for any directed family
	    $\{A_i\}_{i\in I}$ of subobjects of an arbitrary object
	    $X$ and for any subobject $B$ of $X$ the following relation holds:
	    $$(\sum_{i\in I}A_i)\cap B=\sum_{i\in I}(A_i\cap B)$$
\end{itemize}

Ab5-categories possessing a family of generators are called {\em
Grothendieck categories}. They constitute a natural extension of the
class of module categories, with which they share a great number of
important properties.

Although our interest lies exclusively with the category $\Rfp$ with $A$
a ring (associative, with identity), here by a {\em module category\/}
$\mc A$ we mean a Grothendieck category possessing a family of finitely
generated projective generators $\cc A=\{P_i\}_{i\in I}$. We also refer to the
family $\cc A$ as a {\em ring.} It is equivalent to the category $\fc A$
of additive functors from the category $\cc A^{\text{op}}$ to the category
of abelian groups Ab. Conversely, every functor category $\ffc B$ with
$\cc B$ a preadditive category is a module category $\mc A$ with a ring
of representable functors $\cc A=\{h^B=(B,-)\}_{B\in\cc B}$. When
$\cc A=\{A\}$ with $A$ a ring we write $\mc A=\Rfp$.

In the fundamental work of Gabriel~\cite{Ga} it was presented perhaps the
basic tool for studying Grothendieck categories: localization theory.
Such concepts as a {\em localizing subcategory,} a {\em quotient
category\/} (and arising in this case the respective concepts of
a {\em torsion functor\/} and a {\em localizing functor\/}) play
a key role in our analysis.

Recall that the full subcategory $\cc S$ of a Grothendieck category $\cc C$
is a {\em Serre subcategory\/} if for every short exact sequence
   $$\les{C'}C{C''}$$
in $\cc C$ the object $C\in\cc S$ \ifff $C'$, $C''\in\cc S$. A Serre
subcategory $\cc S$ of $\cc C$ is {\em localizing\/} if it is closed
under coproducts.

Gabriel observed that there is a bijection on the class of all idempotent
topoligizing sets of right ideals of the ring $A$ (``The Gabriel topologies''),
and the class of all localizing subcategories of $\Rfp$. Here we define
a Gabriel topology on an arbitrary ring $\cc A=\{P_i\}_{i\in I}$ as follows.
Let $\f F^i$ be some family of subobjects of $P_i\in\cc A$. Then the family
$\f F=\{\f F^i\}_{i\in I}$ is a {\em Gabriel topology on $\cc A$\/}
if it satisfies the following axioms:

\begin{itemize}
\item[$T1.$] $P_i\in\f F$ for each $i\in I$.\par
\item[$T2.$] If $\f a\in\f F^i$ and $\mu\in\homa{P_j}{P_i}$, $P_j\in\cc A$, then
	$\{\f a:\mu\}=\mu^{-1}(\f a)$ belongs~to~$\f F^j$.\par
\item[$T3.$] If $\f a$ and $\f b$ are subobjects of
	$P_i$ such that $\f a\in\f F^i$ and $\{\f b:\mu\}\in\f F^j$ for any
	$\mu\in\homa{P_j}{P_i}$ with $\im\mu\subset\f a$, $P_j\in\cc A$,
	then $\f b\in\f F^i$.
\end{itemize}

\begin{thmm}[Gabriel]\label{Gabriel}
The map
   $$\cc S\longmapsto\f F(\cc S)=\{\f a\subseteq P_i\mid i\in I, P_i/\f a\in\cc S\}$$
establishes a bijection between Gabriel to\-po\-lo\-gies on $\cc A$ and the
class of localizing subcategories of $\mc A$.
\end{thmm}

Similar to the category of modules $\Rfp$ the respective localizing functor
   $$\qq{(-)}S:\mc A\lra{}\mc A/\cc S$$
is constructed in the following way:
   $$\qq MS(P_i)=\lp_{\f a\in\f F^i}\homa{\f a}{M/\qq tS(M)}.$$
Here $\qq tS$ denotes an $\cc S$-torsion functor that takes each right
$\cc A$-module $M\in\mc A$ to the largest subobject $\qq tS(M)$ of $M$
belonging to $\cc S$.

The famous Popescu-Gabriel Theorem~\cite{PG} reduces (at least in principle)
the Grothendieck categories theory to a study of quotient categories
of the module category $\Rfp$. It asserts that every Grothendieck category
$\cc C$ with a family of generators $\cc U=\{U_i\}_{i\in I}$ is equivalent
to the quotient category $\Rfp/\cc S$ with $A$ the endomorphism ring $\End U$
of a generator $U=\ps_{i\in I}U_i$ of $\cc C$ and $\cc S$ some localizing
subcategory of $\Rfp$.

However the main inconvenience of working with Popescu-Gabriel's Theorem
is that, as a rule, the endomorphism ring $\End U$ has enough complicated
form. It often arises a necessity to study Grothendieck categories locally,
with the help of ``visible'' subcategories of $\cc C$. In this vein,
in~\cite{GG} it has been shown that if $P$ is a finitely generated projective
object of $\cc C$ (if such an object exists), then the full subcategory
generated by the object $P$
   $$\cc S=\{C\in\cc C\mid\Hom_{\cc C}(P,C)=0\}$$
is localizing in $\cc C$ and the quotient category $\cc C/\cc S$ of $\cc C$
with respect to $\cc S$ is equivalent to the category of modules $\Rfp$
with $A=\End P$. In particular, every module category $\mc A$ with
$\cc A=\{P_i\}_{i\in I}$ can be partially covered with the categories
$\Mod A_i$ where $A_i=\End P_i$, i.e.\ the following relation holds:
   $$\bigcup_{i\in I}\Mod A_i\subseteq\mc A.$$

We extend this result to arbitrary projective families of objects. We say
that a family $\cc U$ of the Grothendieck category $\cc C$ {\em generates\/}
a full subcategory $\cc B$ of $\cc C$ if for every object $C$ of $\cc B$
there is an exact sequence
   $$\pp{\ps_{i\in I}U_i}{\ps_{j\in J}U_j}C$$
with $U_i\in\cc U$.

\begin{thmm}
Let $\cc C$ be a Grothendieck category and $\cc U=\{P_i\}_{i\in I}$
some family of projective objects of $\cc C$. Then the
subcategory $\cc S=\{C\in\cc C\mid\qc PC=0\text{ for all
$P\in\cc U$}\}$ is localizing in $\cc C$ and $\cc C/\cc S$ is equivalent
to the quotient category $\mc A/\cc P$, where $\cc A=\{h_P=(-,P)\}_{P\in\cc U}$,
$\cc P$ is some localizing subcategory of $\mc A$.
Moreover, $\cc C/\cc S$ is equivalent to a subcategory of $\cc C$ generated
by $\cc U$. If, in addition, each $P\in\cc U$ is finitely generated, then
$\cc C/\cc S$ is equivalent to the module category $\mc A$.
\end{thmm}

The Popescu-Gabriel Theorem is generalized as follows.

\begin{thmm}[Popescu and Gabriel]
Let $\cc C$ be a Grothendieck category with a family of generators
$\cc U=\{U_i\}_{i\in I}$ and $T=(-,?):\cc C\to\mc A$ be the representation
functor that takes each $X\in\cc C$ to $(-,X)$, where
$\cc A=\{h_{U_i}=(-,U_i)\}_{i\in I}$. Then:

$(1)$ $T$ is full and faithful.

$(2)$ $T$ induces an equivalence between $\cc C$ and the quotient category
$\mc A/\cc S$, where $\cc S$ denotes the largest localizing subcategory
in $\mc A$ for which all modules $TX=(-,X)$ are $\cc S$-closed.
\end{thmm}

The advantage of the Theorem is that we can freely choose a family of
generators $\cc U$ of $\cc C$. To be precise, if $\cc M$ is an arbitrary
family of objects of $\cc C$, then the family $\ol{\cc U}=\cc U\cup\cc M$
is also a family of generators.

We say that an object $C$ of $\cc C$ is {\em $\cc U$-finitely generated\/}
({\em $\cc U$-finitely presented\/}) if there is an epimorphism
$\eta:\ps_{i=1}^nU_i\to C$ (if there is an exact sequence
$\ps_{i=1}^nU_i\to\ps_{j=1}^mU_i\to C$) where $U_i\in\cc U$. The full
subcategory of $\cc U$-finitely generated ($\cc U$-finitely presented)
objects of $\cc C$ is denoted by $\fg_{\cc U}\cc C$ ($\fp_{\cc U}\cc C$).
When every $U_i\in\cc U$ is finitely generated (finitely presented), that is
the functor $(U_i,-)$ preserves direct unions (limits), we write
$\fg_{\cc U}\cc C=\fg\cc C$ ($\fp_{\cc U}\cc C=\fp\cc C$). Then every
Grothendieck category is locally $\cc U$-finitely generated
(locally $\cc U$-finitely presented) that means every object $C$ of $\cc C$
is a direct union (limit) $C=\sum_{i\in I}C_i$ ($C=\lp_{i\in I}C_i$)
of $\cc U$-finitely generated ($\cc U$-finitely presented) objects $C_i$.

Recall also that a localizing subcategory $\cc S$ of $\cc C$ is of prefinite
(finite) type provided that the inclusion functor $\cc C/\cc S\to\cc C$
commutes with direct unions (limits). So the following assertion holds.

\begin{thmm}[Breitsprecher]
Let $\cc C$ be a Grothendieck category with a family of generators
$\cc U=\{U_i\}_{i\in I}$. Then the representation functor
$T=(-,?):\cc C\to((\fp_{\cc U}\cc C)^{\text{\rm op}},\text{\rm Ab})$
defines an equivalence between $\cc C$ and
$((\fp_{\cc U}\cc C)^{\text{\rm op}},\text{\rm Ab})/\cc S$,
where $\cc S$ is some localizing subcategory of $((\fp_{\cc U}\cc C)^{\text{\rm op}},\text{\rm Ab})$.
Moreover, $\cc S$
is of finite type \ifff $\fp_{\cc U}\cc C=\fp\cc C$. In this case, $\cc C$
is equivalent to the category
$\Lex((\fp\cc C)^{\text{\rm op}},\text{\rm Ab})$ of contravariant
left exact functors from $\fp\cc C$ to $\text{\rm Ab}$.
\end{thmm}

The Gabriel topology stated above is taken up again in expanded and more
general form, to be used for the characterization of finiteness conditions
on localizing subcategories. Let $\cc U=\{U_i\}_{i\in I}$ be a family of
generators of the Grothendieck category $\cc C$ and $\cc S$ a localizing
subcategory. By Gabriel topology
$\f G$ on $\cc U$ we mean a collection of the following sets:
   $$\f G^i=\{\f a\subseteq U_i\mid U_i/\f a\in\cc S\}$$
where $U_i\in\cc U$. Then $\f G$ satisfies axioms $T1-T3$ (see above)
and similar to the module category the localizing functor
$\qq{(-)}S:\cc C\lra{}\cc C/\cc S$
is defined by the rule
   $$\qq XS(U_i)=\lp_{\f a\in\f G^i}(\f a,X/\qq tS(X))$$
where $U_i\in\cc U$ and $\qq tS$ is an $\cc S$-torsion functor.

At present in modern theory of rings and modules and also in theory
of abelian categories there is a number of fundamental concepts having
come from model theory of modules. The model theory has brought essentially
new principles and statements of questions that touch upon purely algebraic
objects. For this reason, it has stipulated a number of investigations,
which are, on the one hand, a translation of model-theoretic idioms into
algebraic language and, on the other hand, methods obtained so prove,
in a certain context, a very convenient tool for studying the category of
modules.

It turned out that the most important model-theoretic conceptions are
realized in the category
   $$\rc=(\lfp,\text{Ab})$$
whose objects are additive covariant functors from the category $\lfp$ of
finitely presented left $A$-modules into the category of abelian
groups Ab~\cite{JL}. We refer to the category $\rc$ as the category of
{\em generalized right $A$-modules\/} on account of the right exact, fully
faithful functor $M_A\mapsto\ten M-$ from the category $\Rfp$ of right
$A$-modules to $\rc$. This functor identifies pure-injective right
$A$-modules with injective objects of $\rc$~\cite{GJ}. Auslander~\cite{Au}
observed that the full subcategory $\fp\rc$ of finitely presented functors
of $\rc$ is abelian. Equivalently, $\rc$ is a locally coherent Grothendieck
category.

Note separately that the category $\rc$ plays an important role in modern
representation theory of finite dimensional algebras (see e.g.~\cite{Au1,Kr2}).

One of the fundamental model-theoretic concepts is a {\em Ziegler spectrum\/}
of a ring introduced by Ziegler~\cite{Zi} in model-theoretic terms. Recently
Herzog~\cite{He} and Krause~\cite{Kr1} have defined (algebraically) the
Ziegler spectrum of an arbitrary locally coherent Grothendieck category.
An extension of the Ziegler spectrum to arbitrary locally coherent categories
is especially  constructive in the study of purity in arbitrary locally
finitely presented (not only Grothendieck) categories. Generally, a
locally coherent Grothendieck category plus the couple $(\qq tS,\qq{(-)}S)$
of torsion/localizing functors turns out a useful tool every time in model
theory it arises something new about the category of modules! This is,
if you will, a good heuristic method in the indicated context.

So then the Ziegler spectrum $\zg\cc C$ of a locally coherent Grothendieck
category $\cc C$ is a topological space whose points are the isomorphism
types of the indecomposable injective objects of $\cc C$ and an open basis
of $\zg\cc C$ is given by the collection of subsets
   $$\cc O(C)=\{E\in\zg\cc C\mid\Hom_{\cc C}(C,E)\ne 0\}$$
as $C\in\coh\cc C$ ranges over the coherent objects of $\cc C$. If $\cc C=\rc$,
this space is called a Ziegler spectrum of a ring.

The second half of the paper, for the most part, is devoted to a
generalization of results of~\cite{GG,GG2}. At first, given a ring
$\cc A=\{P_i\}_{i\in I}$ of finitely generated projective generators of
the module category $\mc A=\fc A$, we define, similar to the category $\rc$,
the category of generalized $\cc A$-modules
   $$\rcc=(\fpl A,\text{Ab}).$$
Here $\fpl A=\fp\ffc A$ denotes a full subcategory of finitely presented
left $\cc A$-modules. The most important notions of the category $\rc$ are
easily carried over to the same notions for the category $\rcc$.

The following Theorem is the main result for this part.

\begin{thmm}
Let $\cc C$ be a Grothendieck category with a family of generators
$\cc U=\{U_i\}_{i\in I}$ and $\cc A=\{h_{U_i}\}_{i\in I}$ a ring generated
by $\cc U$. Then $\cc C$ is equivalent to the quotient category of $\rcc$
with respect to some localizing subcategory $\cc S$ of $\rcc$. In particular,
every module category $\mc A$ with $\cc A=\{P_i\}_{i\in I}$ a ring is
equivalent to the quotient category of $\rcc$ with
respect to the localizing subcategory
$\ccc P=\{F\in\rcc\mid F(P)=0\text{ for all $P\in\cc A$}\}$.
\end{thmm}

The preceding Theorem leads to a description of different classes
of rings and modules in terms of torsion/localizing functors in the
category $\rcc$. The aim of such a description is to give a criterion of
a duality for categories of finitely presented left and right $\cc A$-modules.

The paper is organized as follows. The first section is preliminary,
collecting the necessary category-theoretic background. In the second section
we discuss Gabriel topologies and localization in module categories. The
following section is of purely technical interest. The results of this
section will be needed in proving Popescu-Gabriel's Theorem. The principle
section of the paper is fourth, in which we prove the Popescu-Gabriel
Theorem and also discuss projective generating sets. In the fifth section
we show how various finiteness conditions of the localizing subcategory
$\cc S$ of $\cc C$ are reflected by properties of the family of generators
$\cc U$ of $\cc C$. We also prove here Breitsprecher's Theorem. The basic
facts about categories of generalized $\cc A$-modules, including
Auslander-Gruson-Jensen's Duality and Theorem of Herzog is presented in the
sixth section. In the remaining sections we present Grothendieck categories
as quotient categories of $\rcc$ and illustrate how localizing subcategories
of $\rcc$ are used to study rings and modules.

When there is no doubt about the ring $A$ or the category $\cc B$,
we usually abbreviate $\Hom_A(M,N)$ or $\Hom_{\cc B}(M,N)$ as $(M,N)$.
We shall freely invoke the fact that every object $X\in\cc C$ of
a Grothendieck category $\cc C$ has an injective envelope $E(X)\in\cc C$
(for details see~\cite[Chapter~VI]{Fr}). This fact is also
discussed in section~\ref{rgc}.
If $\cc B$ is a category, then by a subcategory $\cc A$ of $\cc B$
we shall always mean a {\em full\/} subcategory of $\cc B$. For
concepts such as subobject, epimorphism, injectivity, etc.\ we shall
use the prefix $\cc A$-subobject or $\cc B$-subobject to indicate the
context. This prefix can be omitted if the concept in question is
absolute with respect to the inclusion $\cc A\subseteq\cc B$. To
indicate the context of an operation, for example $\kr\mu$, $E(X)$
or $\lp X_i$, we shall use a subscript, for example, $\kr_{\cc A}\mu$,
$E_{\cc A}(X)$, etc.\ which may also be omitted in case of absoluteness.

\subsubsection*{Acknowledgements.}
I should like to thank A.~I.~Generalov and M.~Prest for helpful discussions.
Many thanks also to Vladimir Kalinin for  his help in preparing the
manuscript and constant interest in the work.

\section{Preliminaries}\label{prel}
In this preliminary section, we collect the basic facts about
Gro\-then\-dieck categories. Some constructions can be defined for more
general (than Grothendieck) categories. For details and proofs we
refer the reader to~\cite{Ga,Fa}.

\subsection{Ab-conditions}\label{Ab}

Recall that an abelian category $\cc C$ is {\em cocomplete\/} or
{\em Ab3-category\/} if it has arbitrary direct sums. The cocomplete
abelian category $\cc C$ is said to be {\em Ab5-category\/} if for any
directed family $\{A_i\}_{i\in I}$ of subobjects of $A$ and for any
subobject $B$ of $A$, the relation
$$(\sum_{i\in I}A_i)\cap B=\sum_{i\in I}(A_i\cap B)$$ holds.

The condition Ab3 is equivalent to the existence of arbitrary direct
limits. Also, Ab5 is equivalent to the fact that there exist inductive
limits and the inductive limits over directed families of indices
are exact, i.e.\ if $I$ is a directed set and
$$\les{A_i}{B_i}{C_i}$$ is an exact sequence for any $i\in I$, then
$$\les{\lp A_i}{\lp B_i}{\lp C_i}$$ is an exact sequence.

Let $\cc C$ be a category and $\cc U=\{U_i\}_{i\in I}$ a family of
objects of $\cc C$. The family $\cc U$ is said to be a {\em family of generators\/}
of the category $\cc C$ if for any object $A$ of $\cc C$
and any subobject $B$ of $A$ distinct from $A$ there exists at least
an index $i\in I$ and a morphism $u:U_i\to A$ which cannot be factorized
through the canonical injection $i:B\to A$ of $B$ into $A$. An object
$U$ of $\cc C$ is said to be a {\em generator\/} of the category $\cc C$
provided that the family $\{U\}$ is a family of generators of the
category $\cc C$.

Let $\cc C$ be a cocomplete abelian category; then $\cc U=
\{U_i\}_{i\in I}$ is a family of generators for $\cc C$ if and only
if the object $\ps_{i\in I} U_i$ is a generator of
$\cc C$~\cite[Proposition~5.33]{BD}. According to~\cite[Proposition~5.35]{BD}
the cocomplete abelian category $\cc C$ which possesses a family
of generators $\cc U$ is locally small and similar to~\cite[Proposition~5.34]{BD}
one can be proved that any object of $\cc C$ is isomorphic to a
quotient of an object $\ps_{j\in J}U_j$, where $J$ is some set of
indices, $U_j\in\cc U$ for any $j\in J$.

An abelian category which satisfies the condition Ab5 and which
possesses a family of generators is called a {\em Grothendieck category.}
\begin{exs}{\rm
(1) The category of left (or right) $A$-modules, where $A$ is a ring,
and the category of (pre-)sheafs of $A$-modules on an arbitrary
topological space are Grothendieck categories.

(2) Let $\cc B$ be a preadditive small category. We denote by
${(\cc B,\text{Ab})}$ the category whose objects are the additive functors
$F:\cc B\to\text{Ab}$ from $\cc B$ to the category of abelian groups Ab and
whose morphisms are the natural transformations between functors. That it is
Grothendieck follows from~\cite[Example~V.2.2]{St}. Besides, the family of
representable functors $\{h^B=(B,-)\}_{B\in\cc B}$ is a family of projective
generators for ${(\cc B,\text{Ab})}$~\cite[Corollary~IV.7.5]{St}.
}\end{exs}

Let $\cc C$ be an abelian category and $\cc U=\{U_i\}_{i\in I}$ some set
of objects of $\cc C$. Consider $\cc U$ as a small preadditive category
and let $\fc U$ be the category of contravariant functors from $\cc U$ to Ab.
By $T:\cc C\to\fc U$ we denote a functor defined as follows.
$$TX=\qc-X,\ Tf=(-,f)$$ where $X\in\cc C$ and $f$ is a morphism of $\cc C$.

\begin{prop}\label{11}
The functor $T:\cc C\to\fc U$ defined above is faithful if and only
if $\cc U$ is a family of generators of $\cc C$.
\end{prop}

\begin{proof}
Assume $T$ is faithful; if $i:X'\to X$ is a monomorphism of $\cc C$ which
is not an isomorphism, and let $j:X\to X/X'$ is the cokernel of $i$, then,
if $Ti$ is an isomorphism, it follows that $Tj=0$ since $T$ is left
exact. Therefore $j=0$, which is a contradiction.

Conversely, assume $\cc U$ is a family of generators of $\cc C$ and let
$f:X\to Y$ be a morphism in $\cc C$ such that $Tf=0$. If $f=pj$ is the
canonical decomposition of $f$ with $p$ monomorphism and $j$ epimorphism,
then $Tj=0$ and if $i:\kr f\to X$ is the kernel of $j$, we get that
$Ti$ is an isomorphism. Since $\cc U$ is a family of generators, $i$
is an isomorphism and therefore $j=0$, hence $f=0$.
\end{proof}

\subsection{Localization in Grothendieck categories}

A subcategory $\cc S$ of the Grothendieck category $\cc C$
is {\em closed under extensions\/} if for any short exact
sequence
   $$\les {X'}X{X''}$$
in $\cc C$ in which $X',X''$ belong to
$\cc S$ the object $X$ belongs to $\cc S$. The subcategory $\cc S$
is a {\em Serre subcategory\/} provided that it is closed under extensions,
subobjects, and quotient objects. The corresponding {\em quotient
category\/} $\cc C/\cc S$ is constructed as follows. The objects
of $\cc C/\cc S$ are those of $\cc C$ and
   $$\q XY=\lp\qc{X'}{Y/Y'}$$
with $X'\subseteq X$, $Y'\subseteq Y$ and $X/X',Y'\in\cc S$.
The set of such pairs is a partially ordered directed set with
respect to the relation $(X',Y/Y')\le (X'',Y/Y'')$, which holds if
and only if $X''\subseteq X'$ and $Y'\subseteq Y''$. The direct
limit is indexed by this partial order.
Again $\cc C/\cc S$ is abelian and there is canonicaly defined the
{\em quotient functor\/} $q:\cc C\to\cc C/\cc S$ such that $q(X)=X$;
it is exact with $\kr q=\cc S$ (see~\cite{Ga,Fa}). Here the {\em kernel\/} $\kr f$
of a functor $f:\cc C\to\cc D$ is, by definition, the subcategory
of all objects $X$ such that $f(X)=0$.

A Serre subcategory $\cc S$ of $\cc C$ is called {\em localizing\/}
provided that the corresponding quotient functor admits a right adjoint
$s:\cc C/\cc S\to\cc C$. Note that $\cc S$ is localizing if and only if
$\cc S$ is closed under taking coproducts~\cite[Theorem~15.11]{Fa}.
In this case, $\cc S$ and $\cc C/\cc S$ are again Grothendieck
categories~\cite[Proposition~III.4.9]{Ga}. Besides, the inclusion
functor $\cc S\to\cc C$ admits a right adjoint $t=\qq tS:\cc C\to\cc S$
which assigns to $X\in\cc C$ the largest subobject $t(X)$ of $X$
belonging to $\cc S$~\cite[Lemma~2.1]{Kr1}. The functor $\qq tS$ is
called a {\em torsion functor} and an object $X$ is called
{\em $\cc S$-periodic\/} or simply {\em periodic\/} provided that
$\qq tS(X)=X$. Furthermore, for any object $X\in\cc C$ there is
a natural morphism $\lambda_X:X\to sq(X)$ such that $\kr\lambda_X$,
$\coker\lambda_X\in\cc S$ and $\kr\lambda_X=\qq tS(X)$ (see~\cite[Chapter~XV]{Fa}).

An object $X\in\cc C$ is said to be {\em $\cc S$-closed\/}
(respectively {\em $\cc S$-torsionfree\/}) provided that $\lambda_X$ is an
isomorphism (respectively a monomorphism). Thus the {\em section functor $s$\/}
induces an equivalence between $\cc C/\cc S$ and the subcategory
of \cl S objects in $\cc C$~\cite[Proposition~15.19B]{Fa}. Moreover, the quotient category $\cc C/\cc S$
is equivalent to the {\em perpendicular category\/} ${\cc S}^\perp$ consisting
of all objects $X\in\cc C$ satisfying $\qc SX=0$ and $\Ext_{\cc C}^1(S,X)=0$
for any $S\in\cc S$~\cite[Lemma~2.2]{Kr1}.
Henceforth, the object $sq(X)$ is denoted by $\qq XS$ and the morphism
$sq(\alpha)$ is denoted by $\qq\alpha S$ for every $X\in\cc C$ and $\alpha
\in\Mor\cc C$; the morphism $\lambda$ we shall call an {\em \env S\/} of
the object $X$. Thus for any object $X$ of $\cc C$ there is an
exact sequence
   \begin{equation}\label{15.19A}
      \ls {A'}{X\bl{\lambda_X}}{\qq XS}{A''}
   \end{equation}
with $A'$, $A''\in\cc S$ and $\lambda_X$ the \env S of $X$.
Note that any two $\cc S$-envelopes
$\lambda^i_X:X\to\qq XS$, $i=1,2$, of $X\in\cc C$ are isomorphic and
$\qq XS\iso\qq {(\qq XS)}S$. Also, note that
$\qq XS=0$ \eqv the object $X$ belongs to $\cc S$.

\begin{prop}\cite[Proposition~15.19C]{Fa}\label{15C}
Let $X$ be an \tr S object; then a monomorphism $\mu:X\to Y$ is an
$\cc S$-envelope \ifff $Y$ is \cl S and $X/Y\in\cc S$. In this case, the
following properties hold:

$(1)$ $\mu$ is an essential monomorphism.

$(2)$ If $E$ is an essential extention of $Y$, then both $E$ and $E/Y$
are \tr S.

Conversely, if $(1)$ and $(2)$ hold and $Y/X\in\cc S$,
then $\mu$ is an $\cc S$-envelope. Moreover, if $E(X)$ is an injective
envelope of $X$ and $X$ is \tr S, then its $\cc S$-envelope is the
largest subobject $D$ of $E(X)$ containing $X$ such that $D/X\in\cc S$.
Thus an \tr S object $X$ is \cl S \ifff $E(X)/X$ is \tr S.
\end{prop}

Let ${\cc S}^\perp$ be the subcategory of \cl S objects. Consider
the {\em localizing functor\/} $\qq{(-)}S:\cc C\to{\cc S}^\perp$, $\qq{(-)}S=sq$;
then the inclusion functor $i:{\cc S}^\perp\to\cc C$, by definition, is fully
faithful and the localizing functor $\qq{(-)}S$ is exact since $q$ is
exact and the section functor $s$, as we have already noticed above,
induces an equivalence of ${\cc S}^\perp$ and $\cc C/\cc S$. Let
$X$, $Y$ be objects of $\cc C$ and $\alpha\in\qc XY$,
$\qq{(-)}S(\alpha)=\qq\alpha S=\qq{(\lambda_Y\alpha)}S$ with $\lambda_Y$
the \env S of $Y$. Clearly that $\qq\alpha S=0$ \ifff $\im\alpha
\subseteq\qq tS(Y)$ whence it easily follows that for $X\in\cc C$,
$Y\in{\cc S}^\perp$ there is an isomorphism
$\qc XY\iso{}_{_{{\cc S}^\perp}}(\qq XS,Y)$, that is $i$ is right adjoint
to the localizing functor $\qq{(-)}S$. On the other hand, if
$\cc C$ and $\cc D$ are Grothendieck categories, $q':\cc C\to\cc D$
is an exact functor and the functor $s':\cc D\to\cc C$ is fully
faithful and right adjoint to $q'$, then $\kr q'$ is a localizing
subcategory and there exists an equivalence
$\cc C/\kr q'\bl H\iso\cc D$ such that $Hq'=q$ with $q$ the canonical
functor~\cite[Proposition~15.18]{Fa}.

Later on, the quotient category $\cc C/\cc S$ always means the
subcategory of \cl S objects ${\cc S}^\perp$ with the pair
of functors $(i,\qq{(-)}S)$, where $i:\cc C/\cc S\to\cc C$ is an
inclusion functor, $\qq{(-)}S:\cc C\to\cc C/\cc S$ is a localizing
functor.

The following Lemma charecterizes \cl S injective objects.

\begin{lem}\label{clinj}
$(1)$ An object $E\in\cc C/\cc S$ is $\cc C/\cc S$-injective
      \ifff it is $\cc C$-injective.

$(2)$ An \tr S and $\cc C$-injective object $E$ is \cl S.
\end{lem}

\begin{proof}
(1). The inclusion functor $i:\cc C/\cc S\to\cc C$ preserves
injectivity since it is right adjoint to the exact functor
$\qq{(-)}S$. If $E\in\cc C/\cc S$ is $\cc C$-injective, then
any $\cc C/\cc S$-monomorphism $\mu:E\to X$ is also a $\cc C$-monomorphism,
and so splits.

(2). It follows from Proposition~\ref{15C}.
\end{proof}

As for $\cc C/\cc S$-morphisms, it holds the following.

\begin{lem}\label{yyy}
Let $\alpha:X\to Y$ be a morphism in $\cc C/\cc S$. Then:

$(1)$ The $\cc C$-kernel of $\alpha$ is \cl S.

$(2)$ $\alpha$ is a $\cc C/\cc S$-epimorphism \ifff $Y/\im_{\cc C}\alpha\in\cc S$.
\end{lem}

\begin{proof}
(1). It suffices to notice that the inclusion functor $i:\cc C/\cc S\to\cc C$,
being right adjoint to the localizing functor $\qq{(-)}S$, is left exact.

(2). Localizing the exact sequence
   $$\pp {X\bl{\alpha}}{Y\bl\beta}{Y/\im_{\cc C}\alpha}$$
with $\beta=\coker\alpha$, we get that $\qq\beta S=0$ that implies
$\qq{(Y/\im_{\cc C}\alpha)}S=0$.
\end{proof}

In particular, a $\cc C/\cc S$-morphism is a mo\-no\-mor\-phism \ifff it is
a $\cc C$-mo\-no\-mor\-phism. We shall refer to this as the
{\em absoluteness\/} of monomorphism. So for $A$, $B\in\cc C/\cc S$
the relation $A\le B$ holds in $\cc C/\cc S$ \ifff it holds in $\cc C$.
Also, it is easily shown that for a $\cc C$-morphism $\alpha:X\to Y$
the $\cc C/\cc S$-morphism $\qq\alpha S$ is a $\cc C/\cc S$-monomorphism
\ifff $\kr\alpha\in\cc S$ and $\qq\alpha S$ is a $\cc C/\cc S$-epimorphism
\ifff $Y/\im\alpha\in\cc S$. Finally $\qq\alpha S$ is a
$\cc C/\cc S$-isomorphism \ifff $\kr\alpha\in\cc S$ and
$Y/\im\alpha\in\cc S$.

\subsection{Lattices of localizing subcategories}

The results of this paragraph are of purely technical interest,
but they will be needed later. Let $\cc C$ be a Grothendieck
category with a family of generators $\cc U=\{U_i\}_{i\in I}$.
Denote by $\lat C$ a lattice consisting of localizing subcategories
of $\cc C$ ordered by inclusion.

Recall that $X\in\cc C$ is $\cc U$-finitely generated provided
that there is an epimorphism $\ps_{i=1}^nU_i\to X$ with $U_i\in\cc U$.
A subcategory consisting of $\cc U$-finitely generated objects
is denoted by $\fg_\cc U\cc C$. The fact that $\lat C$ is a set
follows from that any localizing subcategory $\cc S$ is generated
by its intersection $\fg_\cc U\cc S=\cc S\cap\fg_\cc U\cc C$ with
$\fg_\cc U\cc C$. This means that every object $X\in\cc S$
can be written as a direct union $\sum X_i$ of objects from
$\fg_\cc U\cc S$. Because the category $\fg_\cc U\cc C$
is skeletally small, $\lat C$ is indeed a set.

\begin{prop}\cite[Proposition~2.7]{GG}\label{gri}
Let $\cc P$ and $\cc S$ be localizing subcategories of $\cc C$;
then $\cc P\subseteq\cc S$ \ifff $\cc C/\cc S$ is a quotient
category of $\cc C/\cc P$ with respect to the
localizing in $\cc C/\cc P$ subcategory
$\cc S/\cc P=\{X\in\cc C/\cc P\mid\qq XS=0\}$.
\end{prop}

\begin{prop}\label{dalee}
Let $\cc P$ be a localizing subcategory of $\cc C$ and
$\cc A$ a localizing subcategory of $\cc C/\cc P$; then
there is a localizing subcategory $\cc S$ of $\cc C$
containing $\cc P$ such that $\cc S/\cc P=\cc A$.
\end{prop}

\begin{proof}
Suppose the pair
   $$\qq iP:\cc C/\cc P\lra{}\cc C,\ \qq{(-)}P:\cc C\lra{}\cc C/\cc P$$
defines $\cc C/\cc P$ as a quotient category of $\cc C$.
Also, suppose
   $$\qq iA:(\cc C/\cc P)/\cc A\lra{}\cc C/\cc P,\
     \qq{(-)}A:\cc C/\cc P\lra{}(\cc C/\cc P)/\cc A$$
defines $(\cc C/\cc P)/\cc A$ as a quotient category of
$\cc C/\cc P$.

Denote by
   \begin{gather*}
    Q=\qq{(-)}A\circ\qq{(-)}P:\cc C\lra{}(\cc C/\cc P)/\cc A\\
    I=\qq iP\circ\qq iA:(\cc C/\cc P)/\cc A\lra{}\cc C.
   \end{gather*}
Then $Q$, being a composition of exact fuctors, is an exact
functor. Similarly, $I$, being a composition of fully faithful
fuctors, is a fully faithful functor. Furthermore, given
$X\in\cc C$ and $Y\in(\cc C/\cc P)/\cc A$, we have
   $$\qc XY\iso{}_{\cc C/\cc P}(\qq XP,Y)
     \iso{}_{(\cc C/\cc P)/\cc A}(Q(X),Y).$$
Hence $Q$ is a left adjoint functor to $I$. Thus the pair
$(I,Q)$ defines $(\cc C/\cc P)/\cc A$ as a quotient category
of $\cc C$ with respect to the localizing subcategory
$\cc S=\kr Q$. By construction of $\cc S$ it is easily
seen that $\cc P\subseteq\cc S$ and $\cc S/\cc P=\cc A$.
\end{proof}

Given a localizing subcategory $\cc P$ of $\cc C$, consider
the following sublattice of $\lat C$:
   $$\text{\sf L}_{\cc P}(\cc C)=
     \{\cc S\in\lat C\mid\cc S\supseteq\cc P\}.$$

\begin{cor}
If $\cc P$ is a localizing subcategory of $\cc C$, then
the map
   $$L:\text{\sf L}_{\cc P}(\cc C)\lra{}
    \text{\sf L}(\cc C/\cc P),\ \cc S\longmapsto\cc S/\cc P$$
is a lattice isomorphism.
\end{cor}

Note also that for any $\cc S\in\text{\sf L}_{\cc P}(\cc C)$
   $$\cc S/\cc P=\qq{\cc S}P=\{\qq SP\mid S\in\cc S\}.$$
Indeed, clearly that $\cc S/\cc P\subset\qq{\cc S}P$. In turn,
for $S\in\cc S$ consider exact sequence~\eqref{15.19A}
   $$\ls{A'}S{\qq SP}{A''}$$
with $A'$, $A''\in\cc P$. Because $\cc P\subseteq\cc S$, it follows
that $A'$, $A''\in\cc S$. Hence $\qq SP\in\cc S$ and since
$\qq SP$ is \cl P, one gets $\qq SP\in\cc S/\cc P$.

\subsection{Locally finitely presented Grothendieck categories}

Throughout this paragraph we fix a Grothendieck category $\cc C$.
We define here the most impotant subcategories of $\cc C$, essentially
used further. Namely we describe the subcategories consisting of
finitely generated, finitely presented and coherent objects respectively.
These categories are orded by inclusion as follows:
   $$\cc C\supseteq\fg\cc C\supseteq\fp\cc C\supseteq\coh\cc C.$$

Recall an object $A\in\cc C$ is {\em finitely generated\/} if whenever
there are subobjects $A_i\subseteq A$ for $i\in I$ satisfying
$A=\sum_{i\in I}A_i$, then there is already a finite subset
$J\subset I$ such that $A=\sum_{i\in J}A_i$. The category
of finitely generated subobjects of $\cc C$ is denoted by $\fg\cc C$.
The category is {\em locally finitely generated\/} provided that
every object $X\in\cc C$ is a directed sum
$X=\sum_{i\in I}X_i$ of finitely generated subobjects $X_i$ or
equivalently,
$\cc C$ possesses a family of finitely generated generators.

\begin{thm}\label{fingen}\cite[Proposition~V.3.2]{St}
An object $C\in\cc C$ is finitely generated if and only if the canonical
homomorphism $\varPhi:\lp\qc C{D_i}\to\qc C{\sum D_i}$ is an isomorphism
for every object $D\in\cc C$ and directed family $\{D_i\}_I$ of subobjects
of $D$.
\end{thm}

 A finitely generated object $B\in\cc C$ is {\em finitely presented\/}
provided that every epimorphism $\eta:A\to B$ with $A$
finitely generated has a finitely generated kernel
$\kr\eta$. The subcategory of finitely presented objects of $\cc C$
is denoted by $\fp\cc C$. The respective categories of finitely presented
left and right $A$-modules over the ring $A$ are denoted by $\lfp=\fp (\Lfp)$
and $\rfp=\fp (\Rfp)$. Notice that the subcategory $\fp\cc C$ of $\cc C$
is closed under extensions. Besides, if
   $$\les ABC$$
is a short exact sequence in $\cc C$ with $B$ finitely presented,
then $C$ is finitely presented \ifff $A$ is finitely generated.

The most obvious example of a finitely presented object of $\cc C$ is a
finitely generated projective object $P$. We say that $\cc C$ has
{\em enough\/} finitely generated projectives provided that
every finitely generated object $A\in\cc C$ admits an epimorphism
$\eta:P\to A$ with $P$ a finitely generated projective object.
If $\cc C$ has enough finitely generated projectives, then by the
remarks above, every finitely presented object $B\in\cc C$ is isomorphic
to the cokernel of a morphism between finitely generated projective objects.
This is expressed by an exact sequence
   $$\pp {P_1}{P_0}B$$
called a {\em projective presentation\/} of $B$.

\begin{exs}{\rm
The category $\Lfp$ of left $A$-modules has enough finitely generated
projectives.

Another example of a category having enough finitely generated
projectives is the category of functors $\fcat B$ from the small
preadditive category $\cc B$ to the category of abelian groups Ab. In this
category every finitely generated projective object is a coproduct
factor of a finite coproduct of representable objects $\ps_{i=1}^n(B_i,-)$
(see~\cite[\S1.2]{He}). In addition, if $\cc B$ is an additive
category, that is $\cc B$ is preadditive, has finite products/coproducts
and idempotents split in $\cc B$, then every finitely generated
projective object in ${(\cc B,\text{Ab})}$
is representable~\cite[Proposition~2.1]{He}.
}\end{exs}

The category $\cc C$ is {\em locally finitely presented\/}
provided that every object $B\in\cc C$ is a direct
limit $B=\lp B_i$ of finitely presented objects $B_i$ or equivalently,
$\cc C$ possesses a family of finitely presented generators.
As an example, any locally finitely generated
Grothendieck category having enough finitely generated
projectives $\{P_i\}_{i\in I}$ is locally finitely
presented~\cite[Appendice]{La}. In this case, $\{P_i\}_{i\in I}$
are generators for $\cc C$. For instance, the set of representable
functors $\{h^B\}_{B\in\cc B}$ of the functor category $\fcat B$
with $\cc B$ a small preadditive category form a family of finitely
generated projective generators for $\fcat B$. Therefore $\fcat B$
is a locally finitely presented Grothendieck category (see~\cite[Proposition~1.3]{He}).

\begin{thm}\label{finpr}\cite[Proposition~V.3.4]{St}
Let $\cc C$ be locally finitely generated. An object $B\in\cc C$ is
finitely presented \ifff the functor $\qc B-:{\cc C}\to\text{\rm Ab}$
commutes with direct limits.
\end{thm}

A finitely presented object $C\in\cc C$ is {\em coherent\/} provided
that every finitely generated subobject $B\subseteq C$ is finitely
presented. Evidently, a finitely generated subobject of a coherent
object is also coherent. The subcategory of coherent objects of
$\cc C$ is denoted by $\coh\cc C$. The category $\cc C$ is
{\em locally coherent\/} provided that every object of $\cc C$ is a
direct limit of coherent objects. Equivalently, $\fp\cc C$ is
abelian~\cite[\S2]{Ro} or $\cc C$ possesses a family of coherent
generators. For example, a category
of left $A$-modules is locally coherent \ifff the ring $A$
is left coherent.

In order to characterize the fact that ${(\cc B,\text{Ab})}$ with $\cc B$
an additive category is locally
coherent, that is $\fp\cc C=\coh\cc C$~\cite[\S2]{Ro}, recall that a
morphism $\psi:Y\to Z$ is a {\em pseudo-cokernel\/} for $\phi:X\to Y$
in $\cc B$ if the sequence $h^Z\bl{(\psi,-)}\to h^Y\bl{(\phi,-)}\to h^X$
is exact, i.e.\ every morphism $\delta:Y\to Z'$ with $\delta\phi=0$
factors trough $\psi$.

\begin{lem}\label{pseudo}\cite[Lemma~C.3]{Kr2}
The following are equivalent for $\cc C$:

$(1)$ $\fp\cc C$ is abelian.

$(2)$ Every morphism in $\cc B$ has a pseudo-cokernel.
\end{lem}

The classical example of a locally coherent Grothendieck category
is the category of right (left) generalized $A$-modules
$\rc=(\lfp,\text{\rm Ab})$ ($\lc=(\rfp,\text{\rm Ab})$) consisting
of covariant additive functors from the category $\lfp$ ($\rfp$) of
finitely presented left (right) $A$-modules to Ab.
By the preceding Lemma the category $\fp\rc$, henceforth the category
of coherent functors, is abelian. As we have already said, the
finitely generated projective objects of $\rc$ are the
representable functors $(M,-)=\Hom_A({}_AM,-)$ for some
$M\in\lfp$ and they are generators for $\rc$.

There is a natural right exact and fully faithful functor
   \begin{equation}\label{onn}
    \ten?-:\Rfp\lra{}\rc
   \end{equation}
which takes each module $M_A$ to the tensor functor $\ten M-$.
Recall that \ses
   $$\les XYZ$$
of right $A$-modules is {\em pure\/} if for any $M\in\rfp$
the sequence of abelian groups
   $$\les{\hm MX}{\hm MY}{\hm MZ}$$
is exact. Equivalently, the $\rc$-sequence
   $$\les{\ten X-}{\ten Y-}{\ten Z-}$$
is exact. The module $Q\in\Rfp$ is {\em pure-injective\/} if the functor
$\hm-Q$ takes pure-monomorphisms to epimorphisms.

Functor~\eqref{onn} identifies pure-injective $A$-modules with
injective objects of $\rc$~\cite[Proposition~1.2]{GJ}
(see also~\cite[Proposition~4.1]{He}). Furthermore, the
functor $\rcoh{\ten M-}$ if and only if $M\in\rfp$~\cite{Au}
(see also~\cite{He}).

The category $\coh\rc$ has enough injectives and they
are precisely objects of the form $\ten M-$ with
$M\in\rfp$~\cite[Proposition~5.2]{He}. Thus every coherent
object $\rcoh C$ has both a projective presentation in $\rc$
   $$\pp{(K,-)}{(L,-)}C$$
and an injective presentation in $\coh\rc$
   $$\ii C{\ten M-}{\ten N-}.$$
Here $K$, $L\in\lfp$ and $M$, $N\in\rfp$.

It should be remarked that most important for the applications
in representation theory of finite dimensional algebras
is the concept of purity because the pure-injective modules
play a prominent role among non-finitely generated modules.
It is therefore that many concepts and problems of the
theory are naturally formulated and solved in the category
$\rc$. For this subject we recommend the reader
Auslander's work~\cite{Au1} and Krause's thesis~\cite{Kr2}.

Another important application came from model theory of
modules since the main its conceptions are realized in
$\rc$ (see~\cite{JL,He}). One of such concepts (``The
Ziegler spectrum'') will be discussed in section~\ref{walker}.

\section{Grothendieck categories possessing finitely generated projective
generators}\label{poss}

The following terminology is inspired from the classical theory for categories
of modules $\Rfp$ where $A$ is a ring.
Similar to $\Rfp$, Grothendieck categories $\cc C$ possessing finitely
generated projective generators $\cc A=\{P_i\}_{i\in I}$ are denoted by
$\mc A$ and $\cc A$ we call a {\em ring of projective generators\/}
$\{P_i\}_{i\in I}$
or simply a {\em ring}. The category $\mc A$ is called a {\it category of
right $\cc A$-modules}.
Finally any submodule $\f a$ of $P_i\in\cc A$ is called an {\em ideal\/}
of the ring $\cc A$ corresponding to the object $P_i$.

\subsection{The Gabriel topology}

Let us consider a family of ideals $\f F=\{\f F^i\}_{i\in I}$,
where $\f F^i$ is some family of ideals of $\cc A$ corresponding to
the object $P_i$. Then $\f F$ is a {\em Gabriel topology on $\cc A$\/}
if it satisfies the following axioms:

\begin{itemize}
\item[$T1.$] $P_i\in\f F$ for each $i\in I$.\par
\item[$T2.$] If $\f a\in\f F^i$ and $\mu\in\homa{P_j}{P_i}$, $P_j\in\cc A$, then
	$\{\f a:\mu\}=\mu^{-1}(\f a)$ belongs~to~$\f F^j$.\par
\item[$T3.$] If $\f a$ and $\f b$ are ideals of $\cc A$ corresponding to
	$P_i$ such that $\f a\in\f F^i$ and $\{\f b:\mu\}\in\f F^j$ for any
	$\mu\in\homa{P_j}{P_i}$ with $\im\mu\subset\f a$, $P_j\in\cc A$,
	then $\f b\in\f F^i$.
\end{itemize}

If $\cc A=\{A\}$ is a ring and $\f a$ a right ideal of $A$,
for an arbitrary endomorphism $\mu:A\to A$ of the module $A_A$
   $$\mu^{-1}(\f a)=\{\f a:\mu(1)\}=\{a\in A\mid \mu(1)a\in\f a\}.$$
On the other hand, given an element $x\in A$, one has $\{\f a:x\}=\mu^{-1}(\f a)$
with $\mu\in\End A$ such that $\mu(1)=x$.

\begin{rem}{\rm
We shall need the following properties of Gabriel topologies
$\f F=\{{\f F}^i\}_{i\in I}$ on the ring $\cc A=\{P_i\}_{i\in I}$:

(1). If $\f a\in\f F^i$ and $\f b$ is an ideal of $\cc A$ corresponding to
$P_i\in\cc A$ containing $\f a$, then $\f b\in\f F^i$. Indeed, if
$\mu\in (P_j,P_i)$ such that $\im\mu\subseteq\f a$,
then $\{\f b:\mu\}=P_j\in\f F^j$.

(2). If $\f a$, $\f b\in\f F^i$, then $\f a\cap\f b\in\f F$. Indeed,
since $\{\f a\cap\f b:\mu\}=\{\f a:\mu\}\cap\{\f b:\mu\}$ for any
$\mu\in (P_j, P_i)$, it follows that $\{\f a\cap\f b:\mu\}=\{\f a:\mu\}\in\f F^j$
for $\mu\in (P_j,P_i)$ such that $\im\mu\subset\f b$.

Thus every $\f F^i$, $i\in I$, is a downwards directed system of
ideals.\label{dw}
}\end{rem}

\begin{thm}[Gabriel]\label{Gabriel}
There is a bijective correspondence between Gabriel
to\-po\-lo\-gies on $\cc A$ and localizing subcategories in $\mc A$.
\end{thm}

\begin{proof} Suppose $\f F=\{\f F^i\}_{i\in I}$ is a Gabriel topology on $\cc A$;
then by $\cc S$ denote the following subcategory in $\mc A$:
$A\in\cc S$ \ifff for any $\delta:P_i\to A$ the kernel $\kr\delta\in~\f F^i$.
We claim that $\cc S$ is a localizing subcategory. Indeed, let $A\in\cc S$ and
$i:A'\to A$ be a monomorphism, and $\delta:P_i\to A'$; then $\kr\delta=
\kr (i\delta)\in\f F^i$. Suppose now $p:A\to A''$ is an epimorphism,
$\delta:P_i\to A''$. Since $P_i$ is projective, there is $\gamma:P_i\to A$
such that $p\gamma=\delta$ whence $\kr\gamma\subset\kr\delta$, and hence
$\kr\delta\in\f F^i$.

Let us show now that $\cc S$ is closed under extensions. To see this,
consider a short exact sequence
   $$0\lra{}A'\lra i A\lra p A''\lra{}0$$
with $A',A''\in\cc S$.
Let $\delta:P_i\to A$; then $\f a=\kr (p\delta)\in\f F^i$. Consider
$\gamma:P_j\to P_i$ such that $\im\gamma\subset\f a$. As $p\delta\gamma=0$,
there exists $\alpha:P_j\to A'$ such that $\delta\gamma=i\alpha$.
Therefore $\{\kr\delta:\gamma\}
=\gamma^{-1}(\kr\delta)=\kr(\delta\gamma)=\kr\alpha\in\f F^j$.
>From $T3$ it follows that $\kr\delta\in\f F^i$. Thus $\cc S$ is a Serre
subcategory. If
$\delta_i:P_i\to\oplus A_k,\,A_k\in\cc S$, then there is $k_1,\ldots,k_s$
such that $\im\delta\subset\oplus_{j=1}^sA_{k_j}$ (since $P_i$
is finitely generated). Because $\cc S$ is a Serre subcategory, it
follows that any finite direct sum of objects
from $\cc S$ belongs to $\cc S$, and so $\kr\delta\in\f F^i$,
hence $\cc S$ is a localizing subcategory.

Conversely, assume that $\cc S$ is a localizing subcategory in $\mc A$. Let
$\f F^i=\{\f a\subset P_i\mid P_i/\f a\in\cc S\}$. Obviously
that $P_i\in\f F^i$. If $\f a\in\f F^i$ and $\delta:P_j\to P_i$ is a morphism,
then $\{\f a:\delta\}=\kr (p\delta)$ with
$p:P_i\to P_i/\f a$ the canonical epimorphism. Since $P_i/\f a\in\cc S$,
it follows that $\{\f a:\delta\}\in\f F^j$.
It remains to check $T3$. Assume that $\f a\in\f F^i$ and let
$\f b\subset P_i$ be
such that for any $\mu:P_j\to P_i$ with $\im\mu\subset\f a$ the ideal
$\mu^{-1}(\f b)\in\f F^j$. Let us
consider an exact sequence
   $$0\lra{}{\f a+\f b/\f b}\lra{}P_i/\f b\lra{}P_i/\f a+\f b\lra{}0.$$
Since $\f a\subset\f a+\f b$, one has $\f a+\f b\in\f F^i$.
Let $p:P_i\to P_i/\f b$ be the canonical
epimorphism, $\gamma_{\mu}=p\mu$ for $\mu\in (P_j,P_i)$.
Because $\mu^{-1}(\f b)=\kr\gamma_{\mu}$, it follows that
$P_j/\mu^{-1}(\f b)=\im\gamma_{\mu}=p(\mu(P_j)+\f b)$. In particular,
if $\im\mu\subset\f a$, then $p(\mu(P_j)+\f b)=P_j/\mu^{-1}(\f b)\in\cc S$,
and hence we obtain then that
  $$\f a+\f b/\f b=\sum\limits_{\mu\in (P_j,P_i):\im\mu\subset\f a,\atop{P_j\in\cc A}} p(\mu(P_j)+\f b)$$
belongs to $\cc S$. Since $\cc S$ is closed under extensions,
we conclude that $P_i/\f b\in\cc S$.
\end{proof}

The following result was obtained by Freyd~\cite[p.~120]{Fr} for
cocomplete abelian categories possessing a family of finitely generated
projective generators (see also~\cite[Corollary~III.6.4]{Pop}).

\begin{prop}\label{Freyd}
Let $\mc A$ be the category of right $\cc A$-modules with
$\cc A=\{P_i\}_{i\in I}$ a ring of
finitely generated projective generators,
$T:\mc A\to\fc A$ a functor which assigns to the right $\cc A$-module
$M$ the functor $\homa-M$. Then $T$ is an equivalence.
\end{prop}

\begin{proof}
By definition the funcor $T$ is exact, faithful by Proposition~\ref{11},
and by Theorem~\ref{finpr} it preserves direct limits. Clearly that
$T$ also preserves finitely presented modules.
>From Yoneda Lemma it easily follows that
   $$\homa{\ps_{i=1}^nP_i}{\ps_{j=1}^mP_j}\iso (T(\ps_{i=1}^nP_i),T(\ps_{j=1}^mP_j)).$$
Let us show that for any finitely presented $\cc A$-modules
$M,N\in\fpr A$ there is an isomorphism
   $$\homa MN\iso (T(M),T(N)).$$
To see this, we consider the following commutative diagram with exact rows:
   $$\begin{CD}
       \ps_{i=1}^nh_{P_i}@>q>>\ps_{j=1}^mh_{P_j}@>p>>TM@>>>0\\
       @.@.@VVfV\\
       \ps_{l=1}^sh_{P_l}@>>q'>\ps_{k=1}^th_{P_k}@>>p'>TN@>>>0.\\
     \end{CD}$$
Since $\ps h_{P_j}$ is projective, there exists $g$ such that
$p'g=fp$. In turn, because $p'gq=fpq=0$, it follows that
$\im(gq)\subset\kr p'=\im q'$. Since $\ps h_{P_i}$ is projective,
there exists $h$ such that $q'h=gq$. Thus there results a commutative
diagram
   $$\begin{CD}
       \ps P_i@>>>\ps P_j@>>>M@>>>0\\
       @Vh'VV@VVg'V@.\\
       \ps P_l@>>>\ps P_k@>>>N@>>>0.\\
     \end{CD}$$
where $T(g')=g$, $T(h')=h$. We get then that there exists $f':M\to N$
such that $T(f')=f$.

Consider now right $\cc A$-modules $M$ and $N$. Write them
as direct limits $M=\lp_IM_i$ and $N=\lp_JN_j$ of finitely presented
$\cc A$-modules $M_i$, $N_j$. One has
   \begin{multline*}
     \homa MN\iso\lo_I\lp_J\homa{M_i}{N_j}\iso\lo_I\lp_J (TM_i,TN_i)\\
     \iso (\lp TM_i,\lp TN_j)\iso (T(\lp M_i),T(\lp N_j))\iso (TM,TN),
   \end{multline*}
that is $T$ is a fully faithful functor. It thus remains to check that
any functor $F\in\fc A$ is isomorphic to $TM$ for some $M\in\mc A$.
To see this, choose for $F$ a projective presentation
   $$\pp{\ps_Ih_{P_i}=T(\ps_IP_i)\bl\alpha}{\ps_Jh_{P_j}=T(\ps_JP_j)}F$$
with $I$, $J$ some sets of indices. Because $T$ is fully faithful, there
exists $\beta:\ps_IP_i\to\ps_JP_j$ such that $T(\beta)=\alpha$.
Define $M$ by the exact sequence
   $$\pp{\ps_IP_i\bl\beta}{\ps_JP_j}M.$$
We result in the following commutative diagram
   $$\begin{CD}
     \ps_Ih_{P_i}@>\alpha>>\ps_Jh_{P_j}@>>>F@>>>0\\
     @|@|@.\\
     T(\ps_IP_i)@>>T(\beta)>T(\ps_JP_j)@>>>TM@>>>0.\\
     \end{CD}$$
Therefore we obtain that $F\iso TM$, i.e.\ $T$ is an equivalence.
\end{proof}

In fact, the technique of the proof of the preceding Proposition
allows to establish more strong result used later. Before
formulating it, we make some notation. Let us consider an
arbitrary Grothendieck category $\cc C$ and suppose
$\cc A=\{P_i\}_{i\in I}$ is some family of finitely generated
projective objects of $\cc C$ (if such a family exists).
For clearness we also denote by $\mc A$
the following subcategory of $\cc C$: an object $M\in\mc A$
\ifff there exists an exact sequence
   $$\pp{\ps_{k\in K}P_k}{\ps_{j\in J}P_j}M$$
with $J$, $K$ some sets of indices, $P_i$, $P_k\in\cc A$.
As usual, let us consider $\cc A$ as a preadditive category
and let $T:\mc A\to\fc A$ be the functor taking $M$ to $(-,M)$ and
$f$ to $(-,f)$.

\begin{prop}\label{da}
The functor $T$ establishes an equivalence of $\mc A$ and
$\fc A$. In particular, if $\cc A=\{P\}$, where $P$ is some
finitely generated projective object of $\cc C$, then $\mc A$
is equivalent to the category of right $A$-modules with
$A=\End P$ the endomorphism ring of $P$.
\end{prop}

\begin{proof}
By the slight modification the proof of the first part repeats the
proof of Proposition~\ref{Freyd}. Furthermore, if $\cc A=\{P\}$,
then $h_P$ is a finitely generated projective generator
for $\fc A$ whence it easily follows that $\fc A$, in view of
the Mitchell Theorem, is equivalent to the category of
right $A$-modules $\Rfp$ with $A=\End P$.
\end{proof}

It will be shown in section~\ref{rgc} that the category $\mc A$
is in fact equivalent to a quotient category $\cc C/\cc S$, where
$\cc S=\{C\in\cc C\mid\qc{P_i}C=0\text{ for all $P_i\in\cc A$}\}$.

>From Proposition~\ref{Freyd} it also follows that in order to
define a right $\cc A$-module
$M$ (respectively an $\cc A$-homomorphism), it suffices to define $M$
as a functor from ${\cc A}^{\text{\rm op}}$ to Ab
(respectively as a natural transformation between functors). And
conversely, any functor $F:{\cc A}^{\text{\rm op}}\to\text{Ab}$
(respectively a natural transformation between functors from
$\fc A$) can be considered as a right
$\cc A$-module (respectively an $\cc A$-homomorphism). Later on, we
shall not distinguish $\cc A$-modules and functors of $\fc A$
and shall freely make use of this fact without additional reserves.

\subsection{Localization in module categories}

Fix a localizing subcategory $\cc S$ of $\mc A$, $\cc A=\{P_i\}_{i\in I}$.
Let $\f F=\{\f F^i\}_{i\in I}$ be the respective Gabriel topology on $\cc A$.
As we have already noticed on p.~\pageref{dw}, $\f F^i$ is a downwards
directed system. Let $X$ be an arbitrary right $\cc A$-module and $t=\qq tS$
the $\cc S$-torsion functor.
For every pair $\f a$, $\f b\in\f F^i$ such that
$\f b\subset\f a$ there is a homomorphism
   $$\homa{\f a}{X/t(X)}\lra{}\homa{\f b}{X/t(X)}$$
induced by inclusion of $\f b$ into $\f a$. Clearly that the abelian
groups $\homa{\f a}{X/t(X)}$ with these homomorphisms form an inductive
system over $\f F^i$.

Consider a functor $H:\mc A\to\mc A$ defined as follows. For every
$M\in\mc A$ and every $P_i\in\cc A$ we put
   \begin{equation}\label{123}
     H(M)(P_i)=\lp_{\f a\in\f F^i}\homa{\f a}{M/t(M)}
   \end{equation}
Let us show that abelian groups~\eqref{123} define $H(M)$ as
a functor from $\cc A^{\text{op}}$ to Ab. To see this, we consider
a morphism $\mu:P_j\to P_i$ and an element $m$ of $H(M)(P_i)$. Let
$u:\f a\to M/t(M)$ be a morphism representing the element $m$ of
the direct limit. We then define $H(M)(\mu)(m)\in H(M)(P_j)$ to be
represented by the composed map
   \begin{equation}\label{234}
     \mu^{-1}(\f a)\lra{\bar\mu}\f a\lra u M/t(M)
   \end{equation}
It is easy to see that $H(M)(\mu)$ is well defined, i.e.\ is independent
of the choice of the representing morphism $u$. Thus $H(M)$ becomes
a right $\cc A$-module.

Let now $f:M\to N$ be a morphism in $\mc A$. It is obvious that $f(t(M))$
is contained in $t(N)$. Thus $f$ induces a unique morphism
$f':M/t(M)\to N/t(N)$. In turn, $f'$ induces a unique morphism
$H(f):H(M)\to H(N)$ which is a homomorphism of right $\cc A$-modules.
This concludes the construction of the functor $H$.

The Gabriel topology $\f F$ ordered by inclusion is a directed set
(see remark on p.~\pageref{dw}). From the fact that Ab satisfies
Ab5-condition and from the construction of $H$ we deduce that
$H$ is a left exact functor. Moreover, if $M\in\cc S$, then $H(M)=0$
since $t(M)=M$.

Let $\zeta_i$ be the canonical morphism from $\homa{P_i}{M/t(M)}$
into the direct limit
$H(M)(P_i)$, $p:M\to M/t(M)$ be the canonical epimorphism. Consider
a map of $\cc A$-modules $\varPhi_M:M\to H(M)$
defined as follows. For $P_i\in\cc A$, $\alpha\in\homa{P_i}M$ we put
$\varPhi_M(P_i)(\alpha)=\zeta_i p\alpha$. Let $\mu:P_j\to P_i$ be a
morphism. From construction of the map $H(M)(\mu)$ it easily follows
that the diagram
    $$\begin{CD}
      \homa{P_i}M@>\varPhi_M(P_i)(\mu)>>H(M)(P_i)\\
      @V(\mu,M)VV@VVH(M)(\mu)V\\
      \homa{P_j}M@>>\varPhi_M(P_j)(\mu)>H(M)(P_j)\\
      \end{CD}$$
is commutative, and so $\varPhi_M$ is an $\cc A$-homomorphism. It is
directly verified that $\varPhi_M$ is functorial in $M$. Thus
we obtain a functorial morphism $\varPhi:1_{\mc A}\to H$.

Concerning the functoral morphism $\varPhi$ we prove the following.

\begin{prop}\label{Phi}
$\kr\varPhi_M$ and $\coker\varPhi_M$ belong to $\cc S$ for any right
$\cc A$-module $M$.
\end{prop}

\begin{proof}
As above, we can construct a morphism $\varPsi_M:M/t(M)\to H(M)$,
$\varPsi_M(P_i)(\mu)=\zeta_i\mu$ with $\mu\in\homa{P_i}{M/t(M)}$. One
analogously varifies that $\varPsi_M$ is an $\cc A$-ho\-mo\-mor\-phism.
Then from definitions of the morphisms $\varPhi_M$ and $\varPsi_M$ it
easily follows that the diagram with exact rows
    $$\begin{CD}
      0@>>>t(M)@>>>M@>>>M/t(M)\\
      @.@.@|@VV\varPsi_MV\\
      0@>>>\kr\varPhi_M@>>>M@>>\varPhi_M>H(M)\\
      \end{CD}$$
is commutative whence one gets $t(M)\subset\kr\varPhi_M$. Let us show
that $\varPsi_M$ is a monomorphism or equivalently, $\varPhi_M$ is
a monomorphism for any \tr S object $M$. Indeed, assume that $M$ is
\tr S and let $\mu:P_i\to M$ be such that $\varPhi_M(P_i)(\mu)=0$.
Then there exists an element $\f a\in\f F^i$ such that the restriction
$\mu|_{\f a}=0$. But this implies that $\f a$ is contained in $\kr\mu$
whence $\kr\mu\in\f F^i$ (see remark on p.~\pageref{dw}), and so
$\im\mu\in\cc S$. Because $M$ is \tr S, it follows that $\im\mu=0$,
that is $\mu=0$. Thus $\kr\varPhi_M=t(M)$. Therefore,
if $H(M)=0$, then $M\in\cc S$. Indeed, we obtain then that
$\kr\varPhi_M=M$, hence $M\in\cc S$. Therefore $H(M)=0$ \ifff $M\in\cc S$.

It remains to check that $\coker\varPhi_M\in\cc S$.
Let $\mu:P_i\to H(M)$. We suffice to show that the ideal
$\mu^{-1}(\im\varPhi_M)$ is an element of $\f F^i$. Indeed,
if $p:P_i\to\coker\varPhi_M$, then there exists $\mu:P_i\to H(M)$
such that $p=\coker\varPhi_M\circ\mu$ since $P_i$ is projective.
Further, because the sequence
   $$\ii{\mu^{-1}(\im\varPhi_M)}{P_i\bl p}{\coker\varPhi_M}$$
is exact and $\mu^{-1}(\im\varPhi_M)$, by assumtion, belongs to
$\f F^i$, it will follow then that $\coker\varPhi_M$ belongs to
$\cc S$. Without loss a generality we can assume that
$\varPhi_M$ is a monomorphism and identify $M$ with
$\im\varPhi_M$. We put $\f b=\mu^{-1}(M)$.

Let $u:\f a\to M$ be an $\cc A$-homomorphism representing $\mu$ in the
direct limit $H(M)(P_i)$ and $\xi:P_j\to P_i$ be such that
$\im\xi\subseteq\f a$. Let us consider the following commutative
diagram
   $$\begin{CD}
     \xi^{-1}(\f b)@>>>{\f b}@>>>M\\
     @V{\iota}VV@VVV@VV{\varPhi_M}V\\
     P_j@>>\xi>P_i@>>\mu>H(M).
     \end{CD}$$
Recall that the element $\mu\xi$ is represented by the composed morphism
$u\bar\xi:\{\f a:\xi\}\to M$ (see sequence~\eqref{234}). In view of
that $\im\xi\subseteq\f a$, we have $\{\f a:\xi\}=P_j$, and therefore
$\varPhi_M(u\bar\xi)=\mu\xi$. Since both squares of the diagram are
pullback, it follows that the outer square is pullback, and so there
exists a morphism $\kappa:P_j\to\xi^{-1}(\f b)$ such that $\iota\kappa=1_{P_j}$
whence $\xi^{-1}(\f b)=P_j\in\f F^j$. By $T3$ we deduce that $\f b\in\f F^i$.
\end{proof}

\begin{thm}
For any right $\cc A$-module $M$ the module $H(M)$ is \cl S. Moreover,
the homomorphism $\varPhi_M$ is an \env S of $M$.
\end{thm}

\begin{proof}
To begin with, we shall show that $H(M)$ is \tr S. Let $S$ be a subobject
of $H(M)$ belonging to $\cc S$, $\mu:P_i\to S$ and let $\mu$ be represented by
$u:\f a\to M/t(M)$ in the direct limit $H(M)(P_i)$. Suppose also
$\xi:P_j\to P_i$ is such that $\im\xi\subseteq\kr\mu$; then the
equality $\mu\xi=0$ implies that the image of the composed map
$u\bar\xi:\{\f a:\xi\}\to M/t(M)$ in $H(M)(P_i)$ equals to zero, that is
(using the properties of direct limits) there is an ideal $\f b\in\f F^j$
such that the restriction $u\bar\xi$ representing $\mu\xi$ to $\f b$
equals to zero, so that $\f b\subseteq\kr(u\bar\xi)$. Then
$\kr(u\bar\xi)\in\f F^j$, i.e.\ $\im(u\bar\xi)\in\cc S$. But $M/t(M)$
is \tr S whence $\im(u\bar\xi)=0$, hence $u\bar\xi=0$. Because this
holds for every $\xi:P_j\to P_i$ such that $\im\xi\subseteq\kr\mu$,
we infer that $\kr\mu\subseteq\kr u$ and since $\kr\mu\in\f F^i$,
it follows that $\kr u$ is also an element of $\f F^i$. But in that
case, $\mu$, being an image of zero homomorphism $u\circ\kr u$ from
$\homa{\kr u}{M/t(M)}$, equals to zero. Since it holds for any
$\mu\in\homa{P_i}S$, we deduce that $S=0$.

Let us prove that a module $M$ is \cl S \ifff $\varPhi_M$ is an isomorphism.
Indeed, if $M$ is \cl S, then, in view of the preceding Proposition,
$\kr\varPhi_M=t(M)=0$, that is $\varPhi_M$ is a monomorphism. Since
$\Ext^1(\coker\varPhi_M,M)=0$, there exists a morphism $\alpha:H(M)\to M$
such that $\varPhi_M\alpha=1_M$, i.e.\ $\coker\varPhi_M$ is a direct
summand of $H(M)$ and since $H(M)$ is \tr S, we conclude that
$\coker\varPhi_M=0$ that imply $\varPhi_M$ is an isomorphism.

On the other hand, if for $M$ the morphism $\varPhi_M$ is an isomorphism,
then $t(M)=\kr\varPhi_M=0$, that is $M$ is \tr S. Now if we showed
that every short exact sequence
    $$\les {M\bl i}NS$$
with $S\in\cc S$ splits, it would follow then that $\Ext^1(S,M)=0$ for
any $S\in\cc S$ that implies $M$ would be \cl S. To see this, consider
a commutative diagram
    $$\begin{CD}
	0@>>>M@>i>>N@>>>S@>>>0\\
	@.@V{\varPhi_M}VV@VV{\varPhi_N}V@.@.\\
	0@>>>H(M)@>>H(i)>H(N)@>>>H(S)@.
      \end{CD}$$
where the bottom row is exact and $H(S)=0$. We deduce that $H(i)$
is an isomorphism. Hence $\varPhi_M^{-1}H(i)^{-1}\varPhi_N i=1_M$ that
implies $i$ is a split monomorphism.

Thus to see that $H(M)$ is \cl S for any module $M$, it
suffices to show that $\varPhi_{H(M)}:H(M)\to H(H(M))$ is an isomorphism.
To begin, let us prove that $H(\varPhi_M)$ is an isomorphism.
>From construction of $H(M)$ it follows that $H(M)=H(M/t(M))$.
Let $p:M\to M/t(M)$ be the canonical epimorphism. If we apply the
functor $H$ to the commutative diagram
    $$\begin{CD}
	M@>p>>M/t(M)\\
	@V{\varPhi_M}VV@VV\varPhi_{M/t(M)}V\\
	H(M)@=H(M/t(M)),
      \end{CD}$$
we obtain $H(\varPhi_M)=H(\varPhi_{M/t(M)})$. Since
$H(\coker\varPhi_{M/t(M)})=0$, it follows that $H(\varPhi_{M/t(M)})$ is an
isomorphism, and so $H(\varPhi_M)$ is an isomorphism.

Further, since $\varPhi$ is a functorial morphism, the following
relations
   \begin{equation}\label{345}
     \varPhi_{H(M)}\varPhi_M=H(\varPhi_M)\varPhi_M
   \end{equation}
and
   \begin{gather}
     \varPhi_{H^2(M)}\varPhi_{H(M)}=H(\varPhi_{H(M)})\varPhi_{H(M)}\notag\\
     \label{456}\varPhi_{H^2(M)}H(\varPhi_M)=H^2(\varPhi_M)\varPhi_{H(M)}
   \end{gather}
hold. Applying the functor $H$ to~\eqref{345}, one gets
$H(\varPhi_{H(M)})H(\varPhi_M)=H^2(\varPhi_M)H(\varPhi_M)$ and since
$H(\varPhi_M)$ is an isomorphism, one has $H(\varPhi_{H(M)})=H^2(\varPhi_M)$.
Then from equalities~\eqref{456} it follows that
    \begin{equation}\label{567}
      \varPhi_{H^2(M)}\varPhi_{H(M)}=\varPhi_{H^2(M)}H(\varPhi_M).
    \end{equation}
Because $H^2(M)$ is \tr S, according to the first part of the proof
of the preceding Proposition $\varPhi_{H^2(M)}$ is a monomorphism, and
so from~\eqref{567} it follows that $\varPhi_{H(M)}=H(\varPhi_M)$, that is
$\varPhi_{H(M)}$ is an isomorphism, hence $H(M)$ is \cl S.

In particular, if we consider an exact sequence
   $$\ls{\kr\varPhi_M=t(M)}{M\bl{\varPhi_M}}{H(M)}{\coker\varPhi_M}$$
and apply the exact localizing functor $\qq{(-)}S$, one obtains
   $$\qq MS\iso\qq{(H(M))}S\iso H(M)$$
whence it immediately follows that $\varPhi_M$ is an \env S of $M$.
\end{proof}

Suppose now $i:\mc A/\cc S\to\mc A$ is an inclusion functor. The functor
$H$ is left adjoint to $i$. An adjunction morphism of $i$ with $H$
is furnished by the functorial morphism $\varPhi:1_{\mc A}\to H=i\circ H$
constructed above and a quasi-inverse morphism to $\varPhi$ is a
functorial morphism $\varPsi:H\circ i\to 1_{\mc A/\cc S}$ defined by
the equality $\varPsi_M=(\varPhi_M)^{-1}$ for all $M\in\mc A/\cc S$.
The fact that $\varPsi_{H(M)}H(\varPhi_M)=1_{H(M)}$ is deduced from
the equality $\varPhi_{H(M)}=H(\varPhi_M)$ proved in the preceding
Theorem. The fact that $i(\varPsi_M)\varPhi_{i(M)}=1_{i(M)}$ is
trivial, so that $H$ is indeed left adjoint to $i$.
But the localizing functor $\qq{(-)}S$ is also left adjoint to
$i$, and therefore the functors $\qq{(-)}S$ and $H$ are equivalent.

\section{The structure of localizing subcategories}

The results of this section are of purely technical interest,
but they will be needed in proving the Popescu-Gabriel Theorem.

\begin{prop}\label{negl}
Suppose that $\cc C$ is an abelian category, $\cc M$ is some class of objects
of $\cc C$, $\cc S=\{C\in\cc C\mid\qc CM=0,
\Ext^1_{\cc C}(C,M)=0 \text{ for all }M\in\cc M\}$. Then the following
assertions hold:
\begin{itemize}
\item[(1)] $\cc S$ is closed under extensions and $S\in\cc S$ \ifff for any
$M\in\cc M$, $X\in\cc C$, and epimorphism $f:X\to S$ the canonical
homomorphism $(X,M)\to (\kr f,M)$ is an isomorphism.
\end{itemize}
In addition, if $\cc C$ is cocomplete, then $\cc S$ is closed
under taking coproducts.
\begin{itemize}
\item[(2)] For $\cc S$ the following conditions are equivalent:
\begin{itemize}
\item[(a)] $\cc S$ is a Serre subcategory.
\item[(b)] $\cc S$ is closed under subobjects.\par
\item[(c)] $S\in\cc S$ \ifff for any $M\in\cc M$, $X\in\cc C$,
and morphism $f:X\to S$ the canonical homomorphism $(X,M)\to
(\kr f,M)$ is an isomorphism.
\end{itemize}\end{itemize}
In addition, if $\cc C$ is a Grothendieck category, then
$\cc S$ is a localizing subcategory.
\end{prop}

\begin{proof}
(1). Let us consider a short exact sequence in $\cc C$
\begin{equation}
\les {S'\bl i}{S\bl j}{S''},\label{ses}
\end{equation}
which induces an exact sequence
\begin{multline}
0\to (S'',M)\to (S,M)\to (S',M)\to\\\Ext^1(S'',M)
\to\Ext^1(S,M)\to\Ext^1(S',M)\label{*}
\end{multline}
If $M\in\cc M$ and $S',S''\in\cc S$, then evidently $S\in\cc S$,
so $\cc S$ is closed under extensions.\par
Let $S\in\cc S$, $M\in\cc M$, $X\in\cc C$, $f:X\to S$ be an
epimorphism. If we consider an exact sequence
\begin{equation}
0\to (S,M)=0\to (X,M)\to (\kr f,M)\to\Ext^1(S,M)=0\label{**}
\end{equation}
one obtains $(X,M)\iso (\kr f,M)$.\par
Conversely, consider the identity morphism $1_S:S\to S$ of
the object $S$. Then $0=(\kr 1_S,M)\iso (S,M)$. It
remains to check that $\Ext^1(S,M)=0$. It suffices to show
that any short exact sequence
   $$\les {M\bl h}{P\bl f}S$$
splits. By assumtion for $1_M:M\to M$
there is $g:P\to M$ such that $gh=1_M$. So $h$ splits.

In turn, if $\cc C$ is cocomplete, then the fact that $\cc S$
is closed under under taking coproducts follows from the functor
$\Ext^1_{\cc C}(-,M)$ commutes with direct sums.

(2). $\text{(a)}\Rightarrow\text{(b)}$ is trivial.

$\text{(b)}\Rightarrow\text{(c)}$: Suppose that $S\in\cc S$, $f:X\to S$
is an arbitrary morphism; then
$\im f$, by hypothesis, belongs to $\cc S$. Substituting
$S$ for $\im f$ in~\eqref{**}, we get that the canonical morphism
$(X,M)\to (\kr f,M)$ is an isomorphism for all $M\in\cc M$.
The converse is proved similar to~(1).

$\text{(c)}\Rightarrow\text{(a)}$: From the first assertion it
follows that $\cc S$ is closed under extensions. Suppose now
that in exact sequence~\eqref{ses} the object $S\in\cc S$
and $f:X\to S'$, $X\in\cc C$. Then $\kr f=\kr(if)$ whence it
easily follows that $S'\in\cc S$. Further, if we consider
exact sequence~\eqref{*} with $S',S\in\cc S$, we obtain that
also $S''\in\cc S$.

If $\cc C$ is a Grothendieck category, then $\cc S$ is closed under
taking coproducts, and so it is a localizing subcategory~\cite[Theorem~15.11]{Fa}.
\end{proof}

\subsection{Negligible objects and covering morphisms}

Throughout this paragraph $\cc C$ is assumed to
be a Grothendieck category. Let $\cc M$ be some class of objects of the
category $\cc C$; then an object $S\in\cc C$ is {\it \negl M\/} provided that
for any $M\in\cc M$, $X\in\cc C$, and $f:X\to S$ the canonical homomorphism
$(X,M)\to (\kr f,M)$ is an isomorphism.

\begin{example}{\rm
Any localizing subcategory $\cc S$ of $\cc C$ is a subcategory consisting
of $\cc C/\cc S$-negligible objects, where $\cc C/\cc S$ is the quotient
category of $\cc C$ relative to $\cc S$, since $\cc S={(\cc C/\cc S)}^\perp$.
}\end{example}

\begin{lem}\label{largest}
The subcategory $\cc S$ of $\cc C$
consisting of \negl M objects is the largest localizing
subcategory such that all $M\in\cc M$ are \cl S.
\end{lem}

\begin{proof}
Indeed, if $\cc P$ is a localizing subcategory such that any
object $M\in\cc M$ is \cl P, then for any $f:X\to S$,
$S\in\cc P$, the object $\im f\in\cc P$, and hence, if we consider exact
sequence~\eqref{**} with $S=\im f$, we get that the homomorphism
$(X,M)\to (\kr f,M)$ is an isomorphism. Therefore, in view of Proposition~\ref{negl},
$S\in\cc S$, that is $\cc P\subset\cc S$.
\end{proof}

It is useful to have available the following characterization of \negl M
objects in terms of generators $\cc U=\{U_i\}_{i\in I}$ of the category $\cc C$.

\begin{prop}\label{prf}
An object $S$ is \negl M if and only if for an arbitrary finite set $J$ and morphism
$f:U_J=\ps_{j\in J}U_j\to S$, $U_j\in\cc U$, the canonical homomorphism
$(U_J,M)\to (\kr f, M)$ is an isomorphism for every $M\in\cc M$.
\end{prop}

\begin{proof}
The necessary condition is straightforward. Assume the converse.
Let us consider the commutative diagram
$$\dgARROWLENGTH=0.25\dgARROWLENGTH
\begin{diagram}
 \node[2]{K_i=\kr (fp\psi_i)}\arrow[2]{e,t}{\gamma_i}\arrow[4]{s,l}{\delta_i}
 \node{}
 \node{U_i}\arrow[2]{e,t}{\psi_i}
 \node{}
 \node{U_I=\ps_{i\in I}U_i}\arrow[4]{s,r}{p}\\[4]
 \node{0}\arrow{e}
 \node{\kr f}\arrow[4]{e,b}{\phi}
 \node{}\node{}\node{}
 \node{X}\arrow{e,b}{f}
 \node{S}
\end{diagram}$$
where the bottom row is exact, the couple $(K_i,\gamma_i)$
is the kernel of $fp\psi_i$, $p$ is an epimorphism, $U_i\in\cc U$
for any $i\in I$ (here $I=\bigcup_{U\in\cc U}(U,X)$), $\psi_i$ are the canonical monomorphisms,
and $\delta_i$ is a unique morphism that makes the diagram
commute. Let $M\in\cc M$ and $u:X\to M$ a morphism
such that $u\phi=0$. Then $u\phi\delta_i=0$ for any $i\in I$.
By assumption it follows that $up\psi_i=0$ for any $i$, and
therefore $up=0$ that implies $u=0$ since $p$ is an epimorphism.
Thus the canonical homomorphism induced by $\phi$
$$(X,M)\lra{} (\kr f,M)$$ is a monomorphism.

Consider now a morphism $u:\kr f\to M$. Then $u\delta_i:K_i\to M$.
Let $\lambda_i:U_i\to M$ be such that $\lambda_i\gamma_i=
u\delta_i$. The morphisms $\lambda_i$ induce the morphism
$\lambda:U_I\to M$ such that $\lambda\psi_i=\lambda_i$.

Suppose now that $J$ is an arbitrary finite subset of $I$,
$\psi_J$ is the canonical monomorphism $U_J\to U_I$, the couple
$(K_J,\gamma_J)$ is the kernel of $fp\psi_J$, $\delta_J:K_J
\to\kr f$ is a unique morphism such that $p\psi_J\gamma_J=
\phi\delta_J$, $\lambda_J:U_J\to M$ is a unique morphism such that
$\lambda_J\gamma_J=u\delta_J$, and $\psi_i'$ is the canonical
map $U_i\to U_J$, so that $\psi_i=\psi_J\psi_i'$. As the diagram
   $$\begin{CD}
      K_J@>\gamma_J>>U_J\\
      @V\delta_JVV@VVp\psi_JV@.\\
      \kr f@>>\phi>M\\
     \end{CD}$$
is pullback and $p\psi_J\psi_i'\gamma_i=p\psi_i\gamma_i=\phi\delta_i$,
there exists a morphism $\delta_i':K_i\to K_J$ such that
$\gamma_J\delta_i'=\psi_i'\gamma_i$ and $\delta_J\delta_i'=\delta_i$.
Then we have $u\delta_i=u\delta_J\delta_i'=\lambda_J\gamma_J\delta_i'
=\lambda_J\psi_i'\gamma_i$. In turn, $u\delta_i=\lambda_i\gamma_i$
whence, using the hypothesis, we conclude that $\lambda_J\psi_i'=\lambda_i$
for all $i\in J$. But also $\lambda\psi_J\psi_i'=\lambda\psi_i
=\lambda_i$ for all $i\in J$. Consequently, $\lambda_J=\lambda\psi_J$.
On the other hand, since $p\psi_J\gamma_J=\phi\delta_J$ and the
square from the diagram
\begin{equation}\label{kerfp}
\dgARROWLENGTH=0.25\dgARROWLENGTH
\begin{diagram}
\node[4]{\kr p}\arrow[2]{ssw,l}{\kappa'}\arrow[4]{s,r}{\kappa}\\[4]
\node[2]{K=\kr (fp)}\arrow[2]{e,t}{\gamma}\arrow[4]{s,l}{\delta}
 \node{}
 \node{U_I}\arrow[4]{s,r}{p}\\[4]
 \node{0}\arrow{e}
 \node{\kr f}\arrow[2]{e,b}{\phi}
 \node{}
 \node{X}\arrow{e,b}{f}
 \node{S}
\end{diagram}
\end{equation}
is pullback, there exists $\delta_J':K_J\to K$ such that
$\gamma\delta_J'=\psi_J\gamma_J$ and $\delta\delta_J'=\delta_J$.
Then $\lambda\gamma\delta_J'=\lambda\psi_J\gamma_J=
\lambda_J\gamma_J=u\delta_J=u\delta\delta_J'$, i.e.\ (if we
identify $U_J$ with his image in $U_I$) $\lambda\gamma$
coincides with $u\delta$ on $K\cap U_J$ for all finite
$J\subset I$. Because $U_I=\sum U_J$, we conclude that
$\lambda\gamma=u\delta$. From the diagram~\eqref{kerfp} we
have then $\lambda\kappa=\lambda\gamma\kappa'=u\delta\kappa'$.
But since $\phi\delta\kappa'=p\kappa=0$ and $\phi$ is a monomorphism,
it follows that $\delta\kappa'=0$. Hence $\lambda\kappa=0$,
and therefore there is a morphism $v:X\to M$ such that $vp=
\lambda$. But then $v\phi\delta=vp\gamma=\lambda\gamma=u\delta$
and since $\delta$ is an epimorphism, it follows that $u=v\phi$
as was to be proved.
\end{proof}

Let $\cc S$ be a localizing subcategory of the category
$\cc C$. We say that the object $M\in\cc C$ {\em cogenerates\/}
$\cc S$ if $\cc S=\{C\in\cc C\mid\qc CM=0\}$.

\begin{lem}\label{cog}
The localizing subcategory $\cc S$ consisting of \negl M objects is cogenerated
by the objects $E(M)\ps E(E(M)/M)$ where $M\in\cc M$, $E(M)$ (respectively
$E(E(M)/M)$) is an injective envelope of $M$ (respectively $E(M)/M$).
\end{lem}

\begin{proof}
Denote by $\cc P$ the localizing subcategory cogenerated by $E(M)\ps E(E(M)/M)$.
We need to show that $\cc S=\cc P$.

Let $M\in\cc M$. Consider a short exact sequence
$$\les M{E(M)}{E(M)/M}.$$ It induces an exact sequence
$$0\to (S,M)\to (S,E(M))\to (S,E(M)/M)\to\Ext^1(S,M)\to 0$$
for any $S\in\cc C$. If $S\in\cc S$, then $\Ext^1(S,M)=0$ and also
$(S,E(M))=0$ since $M$ is an \tr S object. Therefore $(S,E(M)/M)=0$,
and hence $(S,E(E(M)/M))=0$ that implies $\cc S\subset\cc P$.

On the other hand, suppose $S\in\cc P$; then $(S,E(M)/M)=0$ since
$E(E(M)/M)$ is a \tr P obj\-ect.
Hence we obtain that $(S,M)=0$ and $\Ext^1(S,M)=0$
that means $M$ is \cl P. But $\cc S$, by hypothesis, is the largest
localizing subcategory such that any $M$ is \cl S, hence $\cc P\subset\cc S$.
\end{proof}

In the proof we made use of the fact that an object $M$ of $\cc C$ is \tr S
\ifff the injective envelope $E(M)$ of $M$ is \tr S.

We shall say that an object $C$ of the category $\cc C$ is {\it cyclic\/}
if there is an epimorphism $f:U_i\to C$ fore some $U_i\in\cc U$.

\begin{lem}\label{easy}
Suppose $C,D$ are objects of $\cc C$; then $\qc C{E(D)}=0$ if and only if
$\qc {C'}D=0$ for any cyclic subobject $C'$ of $C$.
\end{lem}

\begin{proof}
Easy.
\end{proof}

\begin{prop}\label{negligible}
Let $\mc A$ be the category of right $\cc A$-modules with $\cc A$
the ring of finitely generated projective generators $\{P_i\}_{i\in I}$ and let
$\cc M$ be some class of $\cc A$-modules. Then an object $S$ is \negl M if
and only if for an arbitrary morphism $f:P_i\to S$, $P_i\in\cc A$, the
canonical homomorphism $(P_i,M)\to (\kr f, M)$ with $M\in\cc M$ is an isomorphism.
\end{prop}

\begin{proof}
The necessary condition is straightforward. Let $S'$ be an arbitrary
cyclic subobject of $S$; then $(S',M)=0$.
Indeed, there is an epimorphism $f:P_i\to S'$ for some $P_i\in\cc A$ and if
we consider a short exact sequence $$\les{\kr f\bl i}{P_i\bl f}{S'}$$
it will follow then $(S',M)=0$ since the homomorphism $(i,M)$, by hypothesis,
is an isomorphism. Therefore, in view of Lemma~\ref{easy}, $(S,E(M))=0$.

Consider now a short exact sequence $$\les M{E(M)}{E(M)/M}.$$
It induces an exact sequence
    $$(S',E(M))\to (S',E(M)/M)\to\Ext^1(S',M)\to 0.$$
But $(S',E(M))=0$ (since $(S,E(M))=0$), and hence $(S',E(M)/M)\iso\Ext^1(S',M)$.
Our proof will be finished if we show that any short exact sequence
    $$\les {M\bl g}{N\bl h}{S'}$$
splits and then $\Ext^1(S',M)=0$.
To see this, let us consider the following commutative diagram
    $$\begin{CD}
	@.@.\kr f@=\kr f\\
	@.@.@Vi'VV@VViV\\
	0@>>>M@>\nu>>\dec N{S'}{P_i}@>\pi>>P_i\\
	@.@|@Vf'VV@VVfV\\
	0@>>>M@>>g>N@>>h>S'\\
      \end{CD}$$
where the lower right square is pullback, so that the morphisms
$f'$ and $\pi$ are epimorphisms since $f$ and $h$ are epimorphisms.
As $P_i$ is a projective $\cc A$-module, it follows that $\nu$ is a split
monomorphism. Let $\pi'$ be the canonical projection onto $M$, that is
$\pi'\nu=1_M$. By assumtion there exists a morphism $\beta:P_i\to M$ such
that $\pi' i'=\beta i$. But $\beta i=\beta\pi i'$, and hence
$(\pi'-\beta\pi)i'=0$. Then there exists $q:N\to M$ such that $\pi'-\beta\pi=qf'$
whence $1_M=\pi'\nu=qf'\nu=qg$, that is $g$ is a split monomorphism that
finishes the proof.
\end{proof}

Let $\cc M$ be a class of objects of the Grothendieck category $\cc C$.
We shall say that a morphism $u:X\to Y$ is {\em $\cc M$-covering\/}
provided that $\coker u$ is an \negl M object.

\begin{lem}\label{exx}
Consider a category of right $\cc A$-modules $\mc A$ with
$\cc A=\{P_i\}_{i\in I}$ a ring.
Then a morphism $u:X\to Y$ is $\cc M$-covering if and only if for any
$P_i\in\cc A$, $M\in\cc M$, and morphism $f:P_i\to Y$ the sequence
$$0\to\homa{P_i}M\bl\phi\to\homa{\dec XY{P_i}}M\bl\psi\to\homa{\kr u}M$$
induced by $u$ and $f$ is exact.
\end{lem}

\begin{proof}
Consider the commutative diagram
    $$\begin{CD}
	@.@.@.\kr (pf)\\
	@.@.@.@VVhV\\
	0@>>>\kr u@>i'>>\dec XY{P_i}@>u'>>P_i\\
	@.@|@Vf'VV@VVfV\\
	0@>>>\kr u@>>i>X@>>u>Y@>>p>\coker u@>>>0\\
      \end{CD}$$
where the rows are exact, the couple $(\kr (pf),h)$ is
the kernel of $pf$. Because $pfu'=puf'=0$, there exists a unique
morphism $q:\dec XY{P_i}\to\kr (pf)$ such that $u'=hq$. We claim
that $q$ is an epimorphism. Indeed, let us consider the folowing
diagram
     $$\begin{CD}
	 \dec XY{P_i}@>q>>\kr (pf)@>h>>P_i\\
	 @Vf'VV@VVV@VVfV\\
	 X@>>n>\kr p@>>k>Y@>>p>\coker u\\
       \end{CD}$$
where $(\kr p,k)$ is the kernel of $p$, $kn=u$, and
the middle arrow means a unique morphism that makes the diagram
commute. As the outer and right squares are pullback, it follows
that also the left square is pullback. But then the epimorphism
$n:M\to\im u=\kr p$ implies $q$ is an epimorphism.

Suppose that $u$ is $\cc M$-covering and $\alpha:P_i\to M$ a morphism
with $\alpha u'=0$. Then $\alpha hq=0$, and so $\alpha h=0$, hence
$\alpha=0$ since $\coker u$ is \negl M. Thus the homomorphism
$\phi$ induced by $u'$ is a monomorphism. Suppose now that a
morphism $\alpha:\dec XYP_i\to M$ satisfies $\alpha i'=0$. Since
$(\kr (pf),q)$ is the cokernel of $i'$, we deduce that there
exists a morphism $t:\kr (pf)\to M$ such that $tq=\alpha$. Since
$\coker u$ is \negl M, there exists $l:P_i\to M$ such that $lh=t$.
But then $lu'=lhq=tq=\alpha$.

Conversely, suppose the sequence of Lemma is exact and $g$ is an arbitrary
morphism $P_i\to\coker u$. Since $P_i$ is projective, there exists a
morphism $f:P_i\to Y$ such that $pf=g$, so that $(\kr (pf),h)$ is the
kernel of $g$. If a morphism $\alpha:P_i\to M$ satisfies $\alpha h=0$,
then $\alpha u'=\alpha hq=0$ whence $\alpha =0$ since $\phi$ is a monomorphism.
Let now $\beta\in\homa{\kr g}M$. Because $\beta qi'=0$
and $\kr\psi=\im\phi$, it follows that there
exists $\alpha:P_i\to M$ such that $\beta q=\alpha u'$, that is
$\beta q=\alpha hq$, and therefore $\beta=\alpha h$. Thus the canonical
homomorphism $\homa{P_i}M\to\homa{\kr g}M$ is an isomorphism. From
Proposition~\ref{negligible} it follows that $\coker u$ is \negl M.
\end{proof}

In proving the preceding Lemma we have used the fact that any right
$\cc A$-module $M$ has an injective envelope $E(M)$. However
there is another way (which does not use the existence
of injective envelopes\label{pri}) of describing
$\cc M$-covering morphisms. Namely it holds the following.

\begin{lem}\label{parf}
An $\cc A$-homomorphism $u:X\to Y$ of right $\cc A$-modules $X$, $Y$
is $\cc M$-covering \ifff for any finite set of indices $J$,
$M\in\cc M$, and morphism $f:P_J=\ps_{j\in J}P_j\to Y$ with
$P_j\in\cc A$ the sequence
$$0\to\homa{P_J}M\bl\phi\to\homa{\dec XY{P_J}}M\bl\psi\to\homa{\kr u}M$$
induced by $u$ and $f$ is exact.
\end{lem}

\begin{proof}
Our assertion is proved by the slight modification of the proof
of the preceding Lemma. Namely it is necessary to
substitute $P_i$ for $P_J$ in the respective
places and then apply Proposition~\ref{prf}.
\end{proof}

\section{Representation of Grothendieck categories}\label{rgc}

Consider a Grothendieck category
$\cc C$ and fix a family of generators $\cc U=\{U_i\}_{i\in I}$ of $\cc C$.
In subsection~\ref{Ab} we have defined a functor $T:\cc C\to\mc A=\fc U$ which
assigns to $X\in\cc C$ the object $(-,X)$ (here $\cc A$ denotes a
ring of finitely generated projective generators $h_{U_i}=(-,U_i)$, $i\in I$,
of the functor category $\fc U$). The following result was also obtained
by Prest~\cite{Prr}.

\begin{thm}[Popescu and Gabriel]\label{main}
Let $\cc C$ be a Grothendieck category with a family of generators
$\cc U=\{U_i\}_{i\in I}$, $\cc A=\{h_{U_i}=(-,U_i)\}_{i\in I}$, and
$T$ be the functor defined above. Then:

$(1)$ $T$ is full and faithful.

$(2)$ $T$ induces an equivalence between $\cc C$ and the quotient category
$\mc A/\cc S$, where $\cc S$ denotes the largest localizing subcategory
in $\mc A$ for which all modules $TX=(-,X)$ are $\cc S$-closed.
\end{thm}
\begin{proof}
(1): From Lemma~\ref{11} it follows that $T$ is a faithful functor. To see
that it is full, we must show that if $X$ and $Y$ are objects in $\cc C$
and $\varPhi:(-,X)\to (-,Y)$ is a functor, then $\varPhi$ is of the
form $\varPhi(f)=\phi f$ for some $\phi:X\to Y$. Denote by
$\varLambda_i=\homa {U_i}X$ the set of all morphisms $U_i\to X$ and
put $\varLambda=\cup_{i\in I}\varLambda_i$. For each $\alpha\in\varLambda_i$
we let $i_\alpha:U_i\to U_\varLambda=\ps_{i\in\varLambda}U_i$ denote the
corresponding injection. There exists a unique morphism $\lambda:
U_\varLambda\to X$ such that $\lambda i_\alpha=\alpha$ for each
$\alpha\in\varLambda$, and $\lambda$ is an epimorphism since
$\cc U$ is a family of generators. Similarly there exists a unique
morphism $\mu:U_\varLambda\to Y$ such that $\mu i_\alpha=\varPhi(\alpha)$
for each $\alpha\in\varLambda$. Let $\kappa:K\to U_\varLambda$ be the
kernel of $\lambda$. We can show that $\mu\kappa=0$, then $\mu$
factors as $\mu=\phi\lambda$ for some $\phi:X\to Y$ and for each
$i\in I$, $\alpha:U_i\to X$ we get $\varPhi(\alpha)=\mu i_\alpha=
\phi\lambda i_\alpha=\phi\alpha$, and our assertion would be proved.

So we need to check that $\mu\kappa=0$. For each finite subset
$J$ of $\varLambda$ and each $\alpha\in J$ there are canonical morphisms
$\pi'_\alpha:U_J=\ps_{i\in J} U_i\to U_\alpha$, $i'_\alpha:U_\alpha\to U_J$
and $i_J:U_J\to U_\varLambda$. Let $\kappa_J:K_J\to U_J$ be the kernel
of the composed morphism $\lambda i_J:U_J\to X$. Since $K$ is
the direct limit of the kernels $K_J$ for all finite subsets $J$
of $\varLambda$, it suffices to show that $\mu i_J\kappa_J=0$. Now
for each $U_i\in\cc U$, $\beta:U_i\to K_J$ we have, using the
fact that $\varPhi$ is a functor, that
   \begin{multline*}
     \mu i_J\kappa_J\beta=\mu i_J(\sum_{\alpha\in J}i'_\alpha\pi'_\alpha)\kappa_J\beta=
     \sum_{\alpha\in J}\mu i_\alpha\pi'_\alpha\kappa_J\beta=
     \sum_{\alpha\in J}\varPhi(\alpha)\pi'_\alpha\kappa_J\beta\\
     =\sum_{\alpha\in J}\varPhi(\alpha\pi'_\alpha\kappa_J\beta)=
     \varPhi(\sum_{\alpha\in J}\lambda i_J i'_\alpha\pi'_\alpha\kappa_J\beta)=
      \varPhi(\lambda i_J\kappa_J\beta)=0
   \end{multline*}
since $\lambda i_J\kappa_J=0$. Since this holds for arbitrary
$\beta:U_i\to K_J$, it follows that $\mu i_J\kappa_J=0$.

(2): Let $\cc S$ be the largest localizing subcategory
in $\mc A$ for which all modules of the form $TX=(-,X)$ are $\cc S$-closed.
This subcategory exists by Lemmas~\ref{largest} and~\ref{cog} and
it is cogenerated by the class of injective modules of the form
$E(TX)\ps E(E(TX)/TX)$. Let $\cc T=\{TX\}_{X\in\cc C}$; then the respective
$\cc T$-negligible objects and $\cc T$-covering morphisms will be
referred to as negligible and covering respectively, omitting
the prefix $\cc T$. Since every module $TX$ is \cl S, there is a
functor $T':\cc C\to\mc A/\cc S$ such that $T=iT'$ with
$i:\mc A/\cc S\to\mc A$ an inclusion functor. We must show that
$T'$ is an equivalence. Since $iT'=T$ is full and faithful by (1),
also $T'$ is full and faithful, and thus it suffices to show that
every \cl S module is isomorphic to a module of the form $TX$.

To see this, for each $\cc A$-module $M$ we choose an exact sequence
   $$\pp {\ps_{i\in I}h_{U_i}\bl\alpha}{\ps_{j\in J}h_{U_j}}M$$
with $U_i\in\cc U$, $I,J$ some sets of indices. Then $\alpha$ induces in a
natural way a morphism $\beta:\ps_{i\in I}U_i\to\ps_{j\in J}U_j$ in $\cc C$.
To be precise, for each $X\in\cc C$ we have a functorial isomorphism
   \begin{align*}
    \psi_0(X):\homa{\ps_{i\in I}h_{U_i}}{TX}&\iso\prod_{i\in I}\homa{h_{U_i}}{TX}\\
    &\iso\prod_{i\in I}\qc{U_i}X\iso\qc{\ps_{i\in I}U_i}X
   \end{align*}
and analogously, one has a functorial isomorphism
   \begin{align*}
    \psi_1(X):\homa{\ps_{j\in J}h_{U_j}}{TX}&\iso\prod_{j\in J}\homa{h_{U_j}}{TX}\\
    &\iso\prod_{j\in J}\qc{U_j}X\iso\qc{\ps_{j\in J}U_j}X.
   \end{align*}
The morphism $\alpha$ induces the functorial morphism
   $$\zeta=\psi_0\circ(\alpha,-)\circ\psi_1^{-1}:\qc{\ps_{j\in J}U_j}-\lra{}\qc{\ps_{i\in I}U_i}-.$$
According to~\cite[Corollary~1.7]{BD}, there exists a unique morphism
$\beta:\ps_{i\in I}U_i\to\ps_{j\in J}U_j$ such that $(\beta,-)=\zeta$.

Now let us define $\wt M$ by the exact sequence
    \begin{equation}\label{tr}
      \pp {\ps_{i\in I}U_i\bl\beta}{\ps_{j\in J}U_j}{\wt M}.
    \end{equation}
We now apply the functor $T'$ to~\eqref{tr}. If we knew that $T'$ preserves
colimits, we would have then the following isomorphisms
    \begin{multline*}
      {(\ps_{i\in I}h_{U_i}},{\ps_{j\in J}h_{U_j})}
      \iso{(\ps_{i\in I}T'U_i},{\ps_{j\in J}T'U_j)}\\
      \iso{(T'(\ps_{i\in I}U_i)},{T'(\ps_{j\in J}U_j))}
      \iso\qc{\ps_{i\in I}U_i}{\ps_{j\in J}U_j}
    \end{multline*}
and obtain the following commutative diagram
in $\mc A/\cc S$
    $$\begin{CD}
	\ps_{i\in I}h_{U_i}@>\qq\alpha S>>\ps_{j\in J}h_{U_j}@>>>\qq MS@>>>0\\
	@|@|@.\\
	\ps_{i\in I}h_{U_i}@>>T'\beta=\qq\alpha S>\ps_{j\in J}h_{U_j}@>>>T'(\wt M)@>>>0\\
      \end{CD}$$
with exact in $\mc A/\cc S$ rows, and that $\qq MS$ is \cl S would imply
$\qq MS=T'(\wt M)$. To conclude the proof it thus remains to show:

\begin{lem}\label{pff}
The functor $T':\cc C\to\mc A/\cc S$ is exact and preserves direct sums.
\end{lem}

\begin{proof}
First we prove the exactness of $T'$. The functor $T'=\qq{(-)}ST$,
where $\qq{(-)}S$ is the respective localizing functor, is
obviously left exact, and so it suffices to prove that $T'$ preserves
epimorphisms. This means that if $u:X\to Y$ is an epimorphism in $\cc C$,
then the morphism $Tu$ of $\mc A$ is covering, i.e., in view of
Lemma~\ref{exx}, that for any object $Z$ of $\cc C$ and $U_i\in\cc U$
we have the exact sequence
\begin{equation}\label{1.1}
0\to\homa{h_{U_i}}{TZ}\to\homa{\dec {h_{U_i}}{TY}{TX}}{TZ}\to\homa{\kr Tu}{TZ}
\end{equation}
induced by any $\cc A$-homomorphism $f:h_{U_i}\to TY$. Since $T$ is full
and faithful and $h_{U_i}=TU_i$, we deduce that there exists a morphism
$g:U_i\to Y$ such that $f=Tg$. Therefore we have the commutative diagram
    $$\begin{CD}
       0@>>>\kr u@>i'>>\dec {U_i}YX@>u'>>U_i@>>>0\\
       @.@.@VVV@VVgV\\
       @.@.X@>>u>Y@>>>0
      \end{CD}$$
with exact rows. The short exact sequence
    $$\les{\kr u}{\dec{U_i}YX}{U_i}$$
induces the exact sequence
    $$\ii{\qc{U_i}Z}{\qc{\dec{U_i}YX}Z}{\qc{\kr u}Z}$$
and so, since $T$ is fully faithful, we have the exact sequence
   $$0\to\homa{h_{U_i}}{TZ}\to\homa{T(\dec{U_i}YX)}{TZ}\to\homa{T(\kr u)}{TZ}.$$
But since $T$ is left exact, $T(\kr u)$ is isomorphic to $\kr(Tu)$,
$T(\dec{U_i}YX)$ is obviously isomorphic to $\dec{h_{U_i}}{TY}{TX}$
and these isomorphisms are functorial. Thus sequences~\eqref{1.1}
are exact for all $Z\in\cc C$, $U_i\in\cc U$ and now it suffices to apply
Lemma~\ref{exx}.

It remains to prove that $T'$ preserves direct sums. Actually we prove
a little more, namely that it preserves direct unions. Suppose
$\{X_\gamma\}_{\gamma\in\varGamma}$ is a directed family of
subobjects of $X\in\cc C$ such that $X=\sum_{\gamma\in\varGamma}X_\gamma$
with $\varGamma$ some set of indices.
We need to show that the  canonical monomorphism
   $$u:\sum_{\gamma\in\varGamma}TX_\gamma\lra{} T(\sum_{\gamma\in\varGamma}X_\gamma).$$
is covering. Let $q$ be an arbitrary morphism $TU_i\to\coker u$. Since
$TU_i=h_{U_i}$ is a projective object of $\mc A$, there exists a morphism
$g:TU_i\to TX$ such that $q=vg$ where $v=\coker u$. Because $T$ is fully
faithful, there exists a morphism $f:U_i\to X$ such that $g=Tf$. Let
$\dec{X_\gamma}X{U_i}$ be the fibered product associated to the scheme
   $$\begin{CD}
     @.U_i\\
     @.@VVfV\\
     X_\gamma@>>u_\gamma>X
     \end{CD}$$
with $u_\gamma$ the canonical monomorphism. Since $T$ is left exact, it
follows that $T(\dec{X_\gamma}X{U_i})$ is the fibered product associated to
the scheme
   $$\begin{CD}
     @.TU_i\\
     @.@VVTfV\\
     TX_\gamma@>>Tu_\gamma>TX.
     \end{CD}$$
We obtain thus that $\sum\limits_{\gamma\in\varGamma}T(\dec{X_\gamma}X{U_i})$
is the fibered product associated to the scheme
   $$\begin{CD}
     @.TU_i\\
     @.@VVTfV\\
     \sum\limits_{\gamma\in\varGamma}TX_\gamma@>>Tu>TX
     \end{CD}$$
that implies $\sum\limits_{\gamma\in\varGamma}T(\dec{X_\gamma}X{U_i})\iso
\dec{(\sum\limits_{\gamma\in\varGamma}TX_\gamma)}{TX}{TU_i}$. Since
   $$\sum\nolimits_{\cc C}(\dec{X_\gamma}X{U_i})\iso
     \dec{(\sum\nolimits_{\cc C}X_\gamma)}X{U_i}=\dec XX{U_i}\iso U_i$$
and $T$ is fully faithful, one obtains
    \begin{multline*}
    \homa{\sum\limits_{\gamma\in\varGamma}T(\dec{X_\gamma}X{U_i})}{TZ}\iso
    \lo\homa{T(\dec{X_\gamma}X{U_i})}{TZ}\\
    \iso\lo\qc{(\dec{X_\gamma}X{U_i})}Z\iso
    \qc{(\sum\limits_{\gamma\in\varGamma}\dec{X_\gamma}X{U_i})}Z
    \iso\qc{U_i}Z\iso\homa{h_{U_i}}{TZ}
    \end{multline*}
for all $Z\in\cc C$. Hence we get, in view of Lemma~\ref{exx},
that $\coker u$ is negligible.
\def\qed{}\end{proof}

  Thus the functor $T'$ is an equivalence. This concludes the proof of
the Popescu-Gabriel Theorem.
\end{proof}

\begin{cor}\label{uh}
Let $\cc C$ be a Grothendieck category with a family of generators
$\cc U=\{U_i\}_{i\in I}$ and $\cc M=\{M_j\}_{j\in J}$ an
arbitrary family of objects of $\cc C$. Put
$\overline{\cc U}=\cc U\cup\cc M$ and
$\cc A=\{h_U=(-,U)\}_{U\in\overline{\cc U}}$. Then
the functor $T:\cc C\to\mc A$, $TX=(-,X)$,
defines an equivalence between $\cc C$ and the quotient category
$\mc A/\cc S$, where $\cc S$ denotes the largest
localizing subcategory in $\mc A$ for which all modules
$TX$ are $\cc S$-closed.
\end{cor}

\begin{proof}
It suffices to note that $\overline{\cc U}$ is again a
family of generators for $\cc C$ (that directly follows
from definition).
\end{proof}

\begin{cor}[Popescu and Gabriel~\cite{PG}]\label{PPG}
If $\cc C$ is a Grothendieck category, $U$ is a generator of $\cc C$,
then the functor $\qc U-$ establishes an equivalence between categories
$\cc C$ and $\Rfp/\cc S$ with $A=\End U$ the endomorphism ring of $U$,
$\cc S$ some localizing subcategory of $\Rfp$.
\end{cor}

\begin{proof}
By Theorem~\ref{main} the category $\cc C$ is equivalent to the quotient
category of $\fc U$ with $\cc U=\{U\}$ the family of generators. Since
the representable functor $h_U$ is a finitely generated projective generator
for $\fc U$, it follows that $\fc U$, in view of the Mitchell Theorem,
is equivalent to the category of right $A$-modules $\Rfp$ whence
our assertion easily follows.
\end{proof}

\begin{rem}{\rm
In the proof of Lemma~\ref{pff} we essentially used
Lemma~\ref{exx}, which uses, in turn, the fact that any
$\cc A$-module has an injective envelope. There are different
ways of (independent of Popescu-Gabriel Theorem~\cite{PG})
proving this fact (e.g.\ see~\cite[Chapter~VI]{Fr}
or~\cite[\S III.3.10]{Pop}). Nevertheless, Lemma~\ref{pff} can
be easily deduced from Lemma~\ref{parf},
which, as we have noticed on p.~\pageref{pri}, does not use
injective envelopes. To be precise, if we
substitute $P_i$ for $P_J$ in the respective places of the proof
of Lemma~\ref{pff} with $J$ an arbitrary finite set of indices,
our proof will literally
repeat given one and then the fact that any object $X$ of a
Grothendieck category has an injective envelope is proved
similar to~\cite[Corollary~6.32]{BD}.
}\end{rem}

Let $\cc P$ be a localizing subcategory of $\cc C$. We shall
identify $\cc C$, via the functor $T:\cc C\to\mc A$,
$\cc A=\{h_{U_i}\}_{i\in I}$, with the quotient category $\mc A/\cc S$.
According to Proposition~\ref{dalee} there is a localizing subcategory
$\cc T$ of $\mc A$ such that $\cc T\supseteq\cc S$ and $\cc T/\cc S=\cc P$.
Moreover, the quotient category $\cc C/\cc P$ is equivalent to the
quotient category $\mc A/\cc T$.

Let $\f F$ be the Gabriel topology on $\cc A$ corresponding to
$\cc T$. Similar to the category of modules, the family
$\f G=\{\f G^i\}_{i\in I}$
   $$\f G^i=\{X\mid X\subseteq U_i, U_i/X\in\cc P\}$$
we call a {\em Gabriel topology\/} on $\cc U$. It is easily
verified that $\f G\subseteq\f F$ and
   $$\f G=\qq{\f F}S=\{\qq{\f a}S\mid\f a\in\f F\}.$$

\begin{prop}
For a Gabriel topology $\f G$ on $\cc U$ the following assertions hold:
\begin{itemize}

\item[$T1'.$] $U_i\in\f G$ for each $i\in I$.

\item[$T2'.$] If $\f a\in\f G^i$ and $\mu\in\qc{U_j}{U_i}$, $U_j\in\cc U$,
	then $\{\f a:\mu\}=\mu^{-1}(\f a)$ belongs to~$\f G^j$.

\item[$T3'.$] If $\f a$ and $\f b$ are subobjects of $U_i$ such that
	$\f a\in\f G^i$ and $\{\f b:\mu\}\in\f G^j$ for any
	$\mu\in\qc{U_j}{U_i}$ with $\im\mu\subset\f a$, $U_j\in\cc U$,
	then $\f b\in\f G^i$.
\end{itemize}
\end{prop}

\begin{proof}
$T1'$: Straightforward.

$T2'$: Because the left square of the commutative diagram
   $$\begin{CD}
      0@>>>\mu^{-1}(\f a)@>>>U_j@>>>U_j/\mu^{-1}(\f a)@>>>0\\
      @.@VVV@VV{\mu}V@VV{\delta}V\\
      0@>>>\f a@>>>U_i@>>>U_i/\f a@>>>0\\
     \end{CD}$$
is pullback, it follows that $\delta$ is a monomorphism. Hence
$U_j/\mu^{-1}(\f a)\in\cc P$ since $\cc P$ is closed under subobjects.

$T3'$: Suppose $\f a\in\f G^i$ and $\mu:h_{U_j}\to h_{U_i}$ is an
$\cc A$-homomorphism such that $\im\mu\subseteq\f a$; then $\f a\in\f F$.
Because $\qq{(-)}T=\qq{(-)}P\qq{(-)}S$, it follows that
$$\qq{(h_{U_j}/\mu^{-1}(\f b))}T=\qq{(U_j/\qq{\mu^{-1}}S(\f b))}P=0$$
since $U_j/\qq{\mu^{-1}}S(\f b)\in\cc P$. Therefore $\mu^{-1}(\f b)\in\f F$
and by $T3$ we deduce that $\f b\in\f F$. Since $\f b$ is \cl S, we
get that $\f b\in\f G$.
\end{proof}

Further, the \env P $\qq XP$ of an arbitrary object $X$ of $\cc C$ is
constructed similar to $\mc A$. Namely, since the localizing functor
$\qq{(-)}T$ factors as $\qq{(-)}T=\qq{(-)}P\qq{(-)}S$, one has
   \begin{multline*}
    \qq XP(U_i)=\qq XT(U_i)=\lp_{\f a\in\f F^i}(\f a,X/\qq tT(X))=\\
    \lp_{\f a\in\f F^i}(\f a,X/\qq tP(X))\iso
    \lp_{\qq{\f a}S\in\f G^i}(\qq{\f a}S,X/\qq tP(X))
   \end{multline*}
and this isomorphism is functorial in both arguments. Here we made
use of the fact that $\qq tT(X)=\qq tP(X)$ for any $X\in\cc C$.

\subsection{Projective generating sets}

Let $\cc U=\{U_i\}_{i\in I}$ be some family of objects of $\cc C$.
By $\cc C_\cc U$ we shall denote a subcategory of $\cc C$ consisting
of objects generated by $\cc U$. To be precise, $C\in\cc C_\cc U$
\ifff there is a presentation
   $$\pp{\ps_IU_i}{\ps_JU_j}C$$
of $C$ by objects from $\cc U$ and say that $\cc U$ is a
{\em generating set\/} for $\cc C_\cc U$. When $\cc U=\{U\}$ we write
$\cc C_U=\cc C_\cc U$. Clearly $\cc C_\cc U=\cc C$ \ifff $\cc U$
is a family of generators of $\cc C$.

To begin, we study categories $\cc C_\cc U$ with $\cc U$
a generating set consisting of projective objects. So
suppose $P$ is a projective object of $\cc C$.
For any $M\in\cc C_P$ there is a projective presentation
   \begin{equation}\label{6.1}
     \pp{\ps_IP_i}{\ps_JP_j}M
   \end{equation}
with $I$, $J$ some sets of indices, $P_i=P_j=P$ for all $i,j$.

\begin{thm}\label{gga}
Let $\cc C$ be a Grothendieck category, $P$ some projective
object of $\cc C$. Then the subcategory
$\cc S=\{C\in\cc C\mid\qc PC=0\}$ is
localizing and the localized object $\qq PS$ is
$\cc C/\cc S$-projective and a $\cc C/\cc S$-generator.
Moreover, the localizing functor $\qq{(-)}S$ induces an
equivalence of $\cc C/\cc S$ and $\cc C_P$.
\end{thm}

\begin{proof}
To begin, we shall show that $\cc S$ is localizing. Indeed, let
   $$\les {A'}A{A''}$$
be a short exact sequence in $\cc C$. If we apply the exact functor
$\qc P-$, we shall get a short exact sequence of abelian groups
   $$\les{\qc P{A'}}{\qc PA}{\qc P{A''}}$$
whence it easily follows that $\qc PA=0$ \ifff $\qc P{A'}=0$ and $\qc P{A''}=0$
that implies $\cc S$ is a Serre subcategory. Furthermore, if
we consider the map
$$\qc P{\ps A_i}\lra\phi\qc P{\prod A_i}\iso\prod\qc P{A_i}=0$$
with $A_i\in\cc S$ and $\phi$ a monomorpism induced by the
canonical monomorphism $\ps A_i\to\prod A_i$, it will be follow
that $\cc S$ is closed under taking coproducts, and therefore $\cc S$
is a localizing subcategory.

In the rest of the proof we shall show that $\qq PS$ is a projective
generator of $\cc C/\cc S$.
First let us consider \ses
   $$\les ABC$$
in $\cc C/\cc S$. It induces an exact sequence in $\cc C$
   $$\ls ABCS$$
with $S\in\cc S$.
If we apply the exact functor $\qc P-$, one gets the following
commutative diagram of abelian groups:
    $$\begin{CD}
	0@>>>{\qc PA}@>>>{\qc PB}@>>>{\qc PC}@>>>0\\
	@.@VVV@VVV@VVV@.\\
	0@>>>{\q {\qq PS}A}@>>>{\q {\qq PS}B}@>>>{\q {\qq PS}C}@>>>0\\
      \end{CD}$$
where vertical arrows are isomorphisms. Therefore $\qq PS$ is
$\cc C/\cc S$-projective.

It remains to check that $\qq PS$ is a generator of $\cc C/\cc S$.
Let $A$ be an arbitrary object of $\cc C/\cc S$, $I=\qc PA$; then
there exists a morphism $u:\ps_{i\in I}P_i\to A$ with $P_i=P$, $u_i=i$
for all $i\in I$. We have an exact sequence in $\cc C$
    $$\ps P_i\lra u A\lra{}\coker u\lra{} 0.$$
Let $w:P\to A$ be an arbitrary morphism. By construction of $u$
we have that $\im w\subseteq\im u$, and so there exists a morphism
$p:P\to\ps P_i$ such that $w=up$ (since $P$ is projective) that
evidently implies $\qc P{\coker u}=0$, i.e.\ $\coker u\in\cc S$.

But then,
   $$u_{{}_{\cc S}}:\qq {(\ps P_i)}S\iso\ps {\qq {{(P_i)}}S}\lra{}\qq AS$$
is a $\cc C/\cc S$-epimorphism. Thus any object $A\in\cc C/\cc S$ is a
quotient object of $\ps\qq{(P_i)}S$ whence it is easily follows that
$\qq PS$ is a generator of $\cc C/\cc S$.

Now let us show that the restriction of $\qq{(-)}S$ to $\cc C_P$
defines an equivalence of $\cc C/\cc S$ and $\cc C_P$. For
$M\in\cc C$ there is an exact sequence of the
form~\eqref{15.19A}
   $$\ls{A'}{M\bl\lambda}{\qq MS}{A''}$$
with $A'$, $A''\in\cc S$. Hence $\lambda$ induces an isomorphism:
   \begin{equation}\label{6.2}
    \qc PM\iso\qc P{\qq MS}\iso\q{\qq PS}{\qq MS}.
   \end{equation}
But then we have that
   $$\qc{\ps P_i}M\iso\prod\qc{P_i}M
     \iso\qc{\ps P_i}{\qq MS}\iso\q{\ps\qq{(P_i)}S}{\qq MS}$$
with $P_i=P$ for any $i\in I$. Now, if we consider projective
presentation~\eqref{6.1} for $M\in\cc C_P$, we obtain a
commutative diagram
   $$\begin{CD}
      0@>>>\qc MN@>>>\qc{\ps P_j}N@>>>\qc{\ps P_i}N\\
      @.@.@VVV@VVV\\
      0@>>>\qc M{\qq NS}@>>>\qc{\ps P_j}{\qq NS}@>>>\qc{\ps P_i}{\qq NS}\\
     \end{CD}$$
with exact rows and vertical arrows isomorphisms. Thus,
   $$\qc MN\iso\qc M{\qq NS}\iso\q{\qq MS}{\qq NS}$$
that implies $\qq{(-)}S|_{\cc C_P}$ is a fully faithful functor.
Finally let $N\in\cc C/\cc S$. Consider a $\cc C/\cc S$-projective
presentation of $N$
   $$\pp{\ps\qq{(P_i)}S\bl\alpha}{\ps\qq{(P_j)}S}N$$
where $P_i=P_j=P$. Then there exists $\beta:\ps P_i\to\ps P_j$
such that $\alpha=\qq\beta S$, hence $N\iso\qq{(\coker\beta)}S$,
as was to be proved.
\end{proof}

\begin{cor}\label{gag}
Under the conditions of Theorem~\ref{gga} the quotient category
$\cc C/\cc S$ is equivalent to
the quotient category $\Rfp/\cc P$, where $A=\End P$ is the
endomorphism ring of $P$, $\cc P$ is some localizing subcategory of $\Rfp$.
\end{cor}

\begin{proof}
It suffices to note that isomorphisms~\eqref{6.2} induce a ring
isomorphism
   $$A=\qc PP\iso\qc P{\qq PS}\iso\q{\qq PS}{\qq PS}.$$
Now our assertion follows from Corollary~\ref{PPG} since
$\qq PS$ is a generator for $\cc C/\cc S$.
\end{proof}

Now we consider an arbitrary family of projective objects
$\cc U=\{P_i\}_{i\in I}$ and let $\cc C_\cc U$ be the category
generated by $\cc U$.

\begin{cor}\label{mzm}
Let $\cc C$ be a Grothendieck category and $\cc U=\{P_i\}_{i\in I}$
some family of projective objects of $\cc C$. Then the
subcategory $\cc S=\{C\in\cc C\mid\qc PC=0\text{ for all
$P\in\cc U$}\}$ is localizing and $\cc C/\cc S$ is equivalent
to the quotient category $\mc A/\cc P$, where $\cc A=\{h_P\}_{P\in\cc U}$,
$\cc P$ is some localizing subcategory of $\mc A$.
Moreover, the localizing functor $\qq{(-)}S$ induces an
equivalence of $\cc C/\cc S$ and $\cc C_\cc U$.
\end{cor}

\begin{proof}
Denote by $Q=\ps_{P\in\cc U}P$; then $\cc S=\{C\in\cc C\mid\qc QC=0\}$.
>From the preceding Theorem it follows that $\cc S$ is localizing and
that $\qq QS=\ps_{P\in\cc U}\qq PS$ is a $\cc C/\cc S$-projective
generator that implies $\qq{\cc U}S=\{\qq PS\mid P\in\cc U\}$ is a family
of projective generators of $\cc C/\cc S$.

In view of isomorphisms~\eqref{6.2}, there is an equivalence
of $\cc U$ and $\cc U_\cc S$. Now our assertion follows
from Theorem~\ref{main}.
\end{proof}

\begin{rem}{\rm
Under the conditions of Corollary~\ref{mzm} the localizing functor
$\qq{(-)}S$ factors through $\cc C_\cc U$
   $$\begin{diagram}
      \node{\cc C}\arrow[2]{e,t}{\qq{(-)}S}\arrow{se,b}{F}
      \node{}
      \node{\cc C/\cc S}\\
      \node[2]{\cc C_\cc U}\arrow{ne,b}{G}
      \node{}
     \end{diagram}$$
where $G$ is a restriction of $\qq{(-)}S$
to $\cc C_\cc U$ and $F$ is constructed as follows.
Let $C\in\cc C$ and $I=\cup_{P\in\cc U}I_P$ with $I_P=\qc PC$;
then there is a morphism $\phi:\ps_{\mu\in I}P_\mu\to C$, $\phi_\mu=\mu$.
Similarly, there is a morphism
$\psi:\ps_{\nu\in J}P_\nu\to\kr\phi$
with $J=\cup_{P\in\cc U}J_P$, $J_P=\qc P{\kr\phi}$, so that
we have the commutative diagram:
   $$\begin{diagram}
      \node{\ps_{\nu\in J}P_\nu}\arrow[2]{e,t}{\zeta}
      \arrow{se,b}{\psi}
      \node{}
      \node{\ps_{\mu\in I}P_\mu}\arrow{e,t}{\phi}\node{C}\\
      \node[2]{\kr\phi}\arrow{ne,b}{i}
      \node{}
     \end{diagram}$$
with $i=\kr\phi$. By definition we put $F(C)$ to be equal
to $\coker\zeta$.
}\end{rem}

Following notation of section~\ref{poss}, given
$\cc A=\{P_i\}_{i\in I}$ a family of finitely generated
projective objects of $\cc C$, the category $\cc C_\cc A$
is denoted by $\mc A$. By Proposition~\ref{da} $\mc A\iso\fc A$.

\begin{thm}\label{thm2.1}
Suppose $\cc C$ is a Grothendieck category and $\cc A=\{P_i\}_{i\in I}$ is
some family of finitely generated projective objects of $\cc C$. Then the
subcategory $\cc S=\{C\in\cc C\mid\qc{P_i}C=0\text{ for all
$P_i\in\cc A$}\}$ is localizing and $\qq{\cc A}S=\{\qq{(P_i)}S\}_{i\in I}$
is a family of finitely generated projective generators of
the quotient category $\cc C/\cc S$. Moreover, the localizing
functor $\qq{(-)}S$ induces an equivalence of $\cc C/\cc S$
and $\mc A$.
\end{thm}

\begin{proof} By Corollary~\ref{mzm} $\qq{\cc A}S$ is a family
of projective generators for $\cc C/\cc S$.
It thus remains to show that every $\qq{(P_i)}S\in\qq{\cc A}S$
is $\cc C/\cc S$-finitely generated.

To see this, we consider an object $X$ of $\cc C/\cc S$.
Let $\{X_j\}_{j\in J}$ be a directed family of $\cc C/\cc S$-subobjects
of $X$ such that $X=\sum_{\cc C/\cc S}X_j$. Because $X_j$ are also
$\cc C$-subobjects, it follows that the $\cc C$-direct union
$\sum\nolimits_{\cc C}X_j$ is a subobject of $X$ and the
quotient object $A=X/\sum\nolimits_{\cc C}X_j$ belongs to
$\cc S$. Indeed, if we apply the exact localizing functor
$\qq{(-)}S$ commuted with direct limits to the short exact sequence
   $$\les {\sum\nolimits_{\cc C}X_j}XA$$
we shall obtain then the short exact sequence
$$\les {\qq{(\sum\nolimits_{\cc C}X_j)}S
     =\sum\nolimits_{\cc C/\cc S}X_j}X{\qq AS}.$$
Whence $\qq AS=0$, that is $A\in\cc S$. For any
$P_i\in\cc A$ one has then
   $$\qc {P_i}{\sum\nolimits_{\cc C}X_j}\iso\qc {P_i}X
     \iso\q{\qq{(P_i)}S}{\sum\nolimits_{\cc C/\cc S}X_j}.$$
Thus,
   $$\q{\qq{(P_i)}S}{\sum\nolimits_{\cc C/\cc S}X_j}
     \iso\qc{P_i}{\sum\nolimits_{\cc C}X_j}
     \iso\lp\qc{P_i}{X_j}
     \iso\lp\q{\qq{(P_i)}S}{X_j}.$$
By Theorem~\ref{fingen} $\qq{(P_i)}S\in\fg\cc C/\cc S$.

The fact that the localizing functor $\qq{(-)}S$ induces an
equivalence between $\cc C/\cc S$ and $\mc A$ is proved
similar to Theorem~\ref{gga}.
\end{proof}

\begin{cor}\cite[Theorem~2.1]{GG}\label{ggg1}
Let $\cc C$ be a Grothendieck category and $P$ be some finitely
generated projective object of $\cc C$. Then the subcategory
$\cc S=\{C\in\cc C\mid\qc PC=0\}$ is localizing
and the functor $\q{\qq PS}-$ gives an equivalence of
categories $\cc C/\cc S$ and $\Rfp$ with $A=\End P$ the endomorphism
ring of the object $P$.
\end{cor}

\begin{proof}
According to the preceding Theorem the localized object $\qq PS$ is a
finitely generated projective generator, and so, in view of the Mitchell
Theorem and Corollary~\ref{gag}, the quotient category
$\cc C/\cc S$ is equivalent to the category
of modules over the ring $A=\End_{\cc C}P$.
\end{proof}

\begin{example}{\rm
Consider the category of generalized right $A$-modules
$\rc=(\lfp,\text{Ab})$.
Let $M\in\lfp$, $R=\End_AM$ and ${\cc S}_M=\{F\in\rc\mid F(M)=0\}$.
By Corollary~\ref{ggg1} there exists an equivalence of categories
$\rc/{\cc S}_M\bl h\to\Mod R$ where $h(F_{{\cc S}_M})=F_{{\cc S}_M}(M)=F(M)$
for any $F\in\rc$. A quasi-inverse functor to $h$ is constructed as follows:
$g:\Mod R\to\rc/{\cc S}_M$, $E\bl g\mapsto(\ten{(M,E)}-)_{{\cc S}_M}$.
In particular, given $F\in\rc$, there is an isomorphism
$F_{{\cc S}_M}\iso(\ten{(M,F(M))}-)_{{\cc S}_M}$.
}\end{example}

Suppose $\mc A$, $\cc A=\{P_i\}_{i\in I}$, is a category
of right $\cc A$-modules and for any $P_i\in\cc A$ we put
    $$\cc S_i=\{M\in\mc A\mid\homa{P_i}M=0\}.$$
By Corollary~\ref{ggg1} $\mc A/\cc S_i\iso\Rfp_i$ with
$A_i=\End P_i$. We consider this equivalence as identification.
Then the following inclusion holds:
   $$\bigcup_{P_i\in\cc A}\Rfp_i\subseteq\mc A.$$
Generally speaking, as it shows the next example, categories
of $A_i$-modules $\Rfp_i$ do not cover $\mc A$ intirely.

\begin{example}{\rm
Consider the category of generalized abelian groups
$\lgn{\bb Z}=(\modd\bb Z,\text{\rm Ab})$. It has been said
in section~1 that the functor $M\mapsto\tn-M$ identifies
pure injective abelian groups and injective objects
of $\lgn{\bb Z}$. By Kaplansky's Theorem~\cite{Ka},
indecomposable pure injective abelian groups $\zg\lgn{\bb Z}$
are precisely abelian groups of the form:
\begin{itemize}
\item{1.} The injective modules $\bb Q$ and, for every
      prime $p$, the Pr\"ufer groups $\bb Z_{p^\infty}$.
\item{2.} Every cyclic group $\bb Z_{p^n}$ of order a prime
      power.
\item{3.} For every prime $p$, the $p$-adic completion
      $\overline{\bb Z}_p=\lo\bb Z_{p^n}$ of the intergers.
\end{itemize}
Ziegler~\cite{Zi} has shown that a similar argument holds
for Dedekind domains.

Let $\cc S_M=\{C\in\lgn{\bb Z}\mid C(M)=0\}$ with $M\in\modd\bb Z$.
As any finitely generated abelian group is isomorphic to
$\bb Z^{\ps m}\ps\bb Z_{p_i^{k_i}}^{\ps n}$ with $p_1\ldots p_n$
prime, it suffices to show that the relation
   $$\Mod\bb Z\bigcup\bigl[\bigcup_{\text{$p$ is prime}}\Mod\bb Z_{p^n}
     \bigr]\varsubsetneq\lgn{\bb Z}$$
holds. According to~\cite[Corollary~2.4]{GG} $\tn-Q\in\zg\lgn{\bb Z}$
belongs to $\lgn{\bb Z}/\cc S_M\iso\Mod A_M$ with
$A_M=\End M$ \ifff $Q\iso\Hom_{A_M}(M,E)$,
where $E$ is an indecomposable injective $A_M$-module
and in this case $E\iso\tn MQ$. Clearly that $\bb Q$,
$\bb Z_{p^\infty}\in\lgn{\bb Z}/\cc S_{\bb Z}$
and $\bb Z_{p^n}\in\lgn{\bb Z}/\cc S_{\bb Z_{p^n}}$.
However $\overline{\bb Z}_p$ does not satisfy the indicated
condition. Indeed, it is not an injective abelian group,
and therefore
$\tn-{\overline{\bb Z}_p}\notin\lgn{\bb Z}/\cc S_{\bb Z}$.
Because for prime $p\ne q$ the object
$(\bb Z_p,-)\in\cc S_{\bb Z_{q^n}}$ and
$\tn{\bb Z_p}{\overline{\bb Z}_p}\iso\bb Z_p\ne 0$,
it follows that $t_{\cc S_{\bb Z_{q^n}}}(\tn-Q)\ne 0$
where $t_{\cc S_{\bb Z_{q^n}}}$ is an
$\cc S_{\bb Z_{q^n}}$-torsion functor.
Hence, in view of Lemma~\ref{clinj},
$\tn-{\overline{\bb Z}_p}\notin\lgn{\bb Z}/\cc S_{\bb Z_{q^n}}$.
Finally
$\tn-{\overline{\bb Z}_p}\notin\lgn{\bb Z}/\cc S_{\bb Z_{p^n}}$
since
   $(\bb Z_{p^n},\tn{\bb Z_{p^n}}{\overline{\bb Z}_p})
     \iso\bb Z_{p^n}\not\iso\overline{\bb Z}_p$.
}\end{example}

\section{Finiteness conditions for localizing subcategories}\label{walker}

Let $\cc S$ be a localizing subcategory of the Grothendieck category
$\cc C$. In this section, we shall investigate how various finiteness properties
of $\cc S$ are reflected by properties
of the family of generators $\cc U=\{U_i\}_{i\in I}$ of $\cc C$.
We study the precise conditions on $\cc U$ which make the quotient
category $\cc C/\cc S$ into e.g.\ a locally finitely generated category.

Let us consider an inclusion functor $i:\cc C/\cc S\to\cc C$. We say
that $\cc S$ is {\em of prefinite type\/} (respectively {\em of finite type\/})
provided that $i$ commutes with direct unions (respectively with
direct limits), i.e.\ for any \cl S object $C$ and any directed
family $\{C_i\}_I$ of \cl S subobjects of $C$ the relation
   $$\sum\nolimits_{\cc C/\cc S}C_i=\sum\nolimits_{\cc C}C_i$$
holds (respectively $\lp_{\cc C/\cc S}C_i=\lp_{\cc C}C_i$).
It should be remarked that if $\cc S$ is of prefinite type, then in
particular every direct sum of \cl S objects is \cl S.
Recall also that a subcategory $\cc B\subseteq\cc A$ of an
abelian category $\cc A$ is {\em exact\/} provided that it is abelian
and the inclusion functor of $\cc B$ into $\cc A$ is exact.

\begin{prop}\label{ttt}\cite[Theorem~3.41]{Fr}
A subcategory $\cc B$ of $\cc A$ is an exact subcategory \ifff the
following two conditions hold:

$(1)$ If $B_1$, $B_2\in\cc B$ then the coproduct $B_1\ps B_2$ is
an object of $\cc B$.

$(2)$ If $\beta:B_1\to B_2$ is a morphism in $\cc B$, then both the
$\cc A$-kernel and $\cc A$-cokernel of $\beta$ are objects of $\cc B$.
\end{prop}

Finally a subcategory $\cc B$ of an abelian category $\cc A$ is called
{\em coexact\/} provided that the perpendicular subcategory
   $$\cc B^{\perp}=\{A\in\cc A\mid (B,A)=0, \Ext^1(B,A)=0\text{ for
     all $B\in\cc B$}\}$$
is exact. For example, if $\cc C$ is a Grothendieck category,
$\cc S\subseteq\cc C$ is localizing, then $\cc S$ is coexact
if and only if the quotient category $\cc C/\cc S$ is exact because
$\cc C/\cc S=\cc S^{\perp}$. Note also that any quotient category
$\cc C/\cc S$ is coexact since $\cc S=(\cc C/\cc S)^{\perp}$ is an
exact subcategory.

\begin{prop}\label{osya}
If $\cc S\subseteq\cc C$ is a coexact localizing subcategory, then
for any projective object $P$ of $\cc C$ the object $\qq PS$ is
$\cc C/\cc S$-projective. If, in addition, $\cc C$ has a family of
projective generarors $\cc A=\{P_i\}_{i\in I}$ and any $\qq{(P_i)}S$
is $\cc C/\cc S$-projective, then $\cc S$ is coexact.
\end{prop}

\begin{proof}
If $\cc S$ is coexact, then any short exact sequence
   \begin{equation}\label{lkg}
   \les A{B\bl\beta}C
   \end{equation}
in $\cc C/\cc S$ is also exact in $\cc C$. Suppose $P\in\cc C$
is projective; then we have the following commutative diagram
   $$\begin{CD}
    0@>>>\qc PA@>>>\qc PB@>(P,\beta)>>\qc PC@>>>0\\
    @.@VVV@VVV@VVV\\
    0@>>>\q{\qq PS}A@>>>\q{\qq PS}B@>>(\qq PS,\beta)>\q{\qq PS}C@>>>0\\
    \end{CD}$$
with exact rows and vertical arrows being isomorphisms. Hence
$\qq PS$ is $\cc C/\cc S$-projective.

Conversely, we need to show that any $\cc C/\cc S$-epimorphism
is $\cc C$-epimorphism. Consider exact sequence~\eqref{lkg}.
By Lemma~\ref{yyy}
$C/\im_{\cc C}\beta\in\cc S$. By assumption the morphism
$(\qq PS,\beta)$ is an epimorphism, where $P\in\cc A$, hence
$(P,\beta)$ is an epimorphism. Therefore $(P,C/\im_{\cc C}\beta)=0$
for every $P\in\cc A$ that implies $C/\im_{\cc C}\beta=0$
since $\cc A$ is a family of generators.
\end{proof}

The following statement charecterizes coexact localizing
subcategories of (pre)finite type.

\begin{prop}\label{coex}
For a coexact localizing subcategory $\cc S$ of a Grothendieck
category $\cc C$ the following assertions are equivalent:

$(1)$ $\cc S$ is of finite type.

$(2)$ $\cc S$ is of prefinite type.

$(3)$ The inclusion functor $i:\cc C/\cc S\to\cc C$ commutes with
      coproducts.
\end{prop}

\begin{proof}
$(1)\Rightarrow (2)\Rightarrow (3)$: Straightforward.

$(3)\Rightarrow (1)$: Let $\lp_{\cc C/\cc S}C_i$ be a
$\cc C/\cc S$-direct limit of $C_i\in\cc C/\cc S$, $i\in I$. Denote by
$\varLambda$ the subset of $I\times I$ consisting of pairs $(i,j)$
with $i\le j$ and for any $\lambda\in\varLambda$ we put
$s(\lambda)=i$, $t(\lambda)=j$. By~\cite[Proposition~IV.8.4]{St}
   $$\lp_{\cc C/\cc S}C_i=\coker_{\cc C/\cc S}\bigl[
   \oplus_{\lambda\in\varLambda}C_{s(\lambda)}\lra{\phi}
   \oplus_{i\in I}C_i\bigr]$$
with $\phi$ induced by
$\phi_{\lambda}=u_j\gamma_{ij}-u_i:C_{s(\lambda)}\to\ps_{i\in I}C_i$,
$\lambda=(i,j)$ and $\gamma_{ij}:C_{s(\lambda)}\to C_{t(\lambda)}$
the canonical morphism. By hypothesis the inclusion functor
$\cc C/\cc S\to\cc C$ is exact and commutes with coproducts,
and so we obtain
   $$\lp_{\cc C/\cc S}C_i=\coker_{\cc C}\bigl[
   \oplus_{\lambda\in\varLambda}C_{s(\lambda)}\lra{\phi}
   \oplus_{i\in I}C_i\bigr]=\lp_{\cc C}C_i$$
as was to be proved.
\end{proof}

\begin{lem}\label{tor}
Let $\cc C$ be a Grothendieck category and $\cc S$ a localizig
subcategory of $\cc C$. Then the $\cc S$-torsion functor $t=\qq tS$
commutes with direct limits \ifff any $\cc C$-direct limit of \cl S
objects is \tr S.
\end{lem}

\begin{proof}
Suppose that $t$ commutes with direct limits. Let us consider an exact
sequence
   $$\ii {t(\lp_{\cc C}C_i)}{\lp_{\cc C}C_i\bl\lambda}{\lp_{\cc C/\cc S}C_i}$$
with $C_i\in\cc C/\cc S$, $\lambda$ an \env S of $\lp C_i$. Applying
the left exact functor $t$, we have
   $$t(t(\lp C_i))=t(\lp C_i)=\lp t(C_i)=0$$
since $t(C_i)=0$. So $\lp_{\cc C}C_i$ is \tr S.

Conversely, let $\{C_i\}_I$ be a direct system of objects from
$\cc C$. As above, there is an exact sequence
   $$\ii {t(C_i)}{C_i\bl{\lambda_{C_i}}}{\qq{(C_i)}S}.$$
Because the direct limit functor is exact, one gets an exact sequence
   \begin{equation}\label{ttt}
    \ii{\lp t(C_i)}{\lp C_i}{\lp\qq{(C_i)}S}.
   \end{equation}
Since $\cc S$ is closed under taking direct limits, it follows
that $\lp t(C_i)\in\cc S$. Now, if we apply $t$ to sequence~\eqref{ttt},
one obtains
   $$\lp t(C_i)=t(\lp t(C_i))=t(\lp C_i)$$
since $t(\lp\qq{(C_i)}S)=0$.
\end{proof}

Let $\f G$ be a Gabriel topology on $\cc U$ corresponding to $\cc S$.
By a {\em basis\/} for the topology $\f G$ we mean a subset
$\f B$ of $\f G$ such that every object in $\f G$
contains some $\f b\in\f B$.

\begin{thm}\label{pref}
Let $\cc C$ be a locally finitely generated Grothendieck category
with the family of finitely generated generators $\cc U\subseteq\fg\cc C$
and suppose that $\cc S$ is a localizing subcategory of $\cc C$.
Then the following conditions are equivalent:

\begin{itemize}
\item[$(1)$] $\cc S$ is of prefinite type.

\item[$(2)$] For any $U\in\cc U$ the natural morphism
      $\lp\qc U{C_i}\to\qc U{\sum\nolimits_{\cc C/\cc S}C_i}$
      induced by the \env S $\lambda$ of $\sum\nolimits_{\cc C}C_i$
      is an isomorphism.

\item[$(3)$] $\qq{\cc U}S=\{\qq US\}_{U\in\cc U}$
      is a family of $\cc C/\cc S$-finitely generated
      generators for $\cc C/\cc S$.

\item[$(4)$] If $C$ is $\cc C$-finitely generated, then $\qq CS$
      is $\cc C/\cc S$-finitely generated.

\item[$(5)$] The torsion functor $t$ commutes with direct limits.

\item[$(6)$] Any $\cc C$-direct limit of \cl S object is \tr S.

\item[$(7)$] $\f G$ has a basis of finitely generated objects.
\end{itemize}

Thus $\cc C/\cc S$ is a locally finitely generated Grothendieck category
and in this case, any $\cc C/\cc S$-finitely generated object
$D$ is a localization $\qq CS$ of some $C\in\fg\cc C$.
\end{thm}

\begin{proof}
Equivalence $(5)\Leftrightarrow (6)$ follows from Lemma~\ref{tor}.

$(1)\Rightarrow (2)$: By definition
$\sum\nolimits_{\cc C/\cc S}C_i=\sum\nolimits_{\cc C}C_i$. Now
our assertion follows from Theorem~\ref{fingen}.

$(2)\Rightarrow (3)$: Let
$\lambda:\sum\nolimits_{\cc C}C_i\to\sum\nolimits_{\cc C/\cc S}C_i$
be an \env S of $\sum\nolimits_{\cc C}C_i$. Then the composed map
   $$\lp\qc U{C_i}\lra{\varPhi}\qc U{\sum\nolimits_{\cc C}C_i}\lra{(U,\lambda)}
     \qc U{\sum\nolimits_{\cc C/\cc S}C_i}$$
with $\varPhi$ the canonical morphism is, by hypothesis, an isomorphism. Hence,
   $$\lp\q{\qq US}{C_i}\iso\lp\qc U{C_i}
     \iso\qc U{\sum\nolimits_{\cc C/\cc S}C_i}
     \iso\q{\qq US}{\sum\nolimits_{\cc C/\cc S}C_i}.$$
By Theorem~\ref{fingen} $\qq US\in\fg\cc C/\cc S$.

$(3)\Rightarrow (4)$: If $C\in\fg\cc C$ there is an epimorphism
$\eta:\ps_{i=1}^nU_i\to C$ for some $U_1\ldots U_n\in\cc U$. Since
$\ps_{i=1}^n\qq{(U_i)}S\in\fg\cc C/\cc S$, it follows that
$\qq CS\in\cc C/\cc S$.

$(4)\Rightarrow (7)$: Suppose $\f a\in\f G$, that is $\qq{\f a}S=\qq US$
for some $U\in\cc U$. Write $\f a=\sum_{\cc C}\f a_i$ as a directed
sum of $\cc C$-finitely generated subobjects $\f a_i$. Then
$\qq US=\sum_{\cc C/\cc S}\qq{(\f a_i)}S$. By assumtion there is
$i_0$ such that $\qq US=\qq{(\f a_{i_0})}S$, and hence $\f a_{i_0}\in\f G$.

$(7)\Rightarrow (1)$: First let us show that $\qq US\in\fg\cc C/\cc S$
where $U\in\cc U$. Indeed, suppose $\qq US=\sum_{\cc C/\cc S}\f a_i$; then
$\f a=\lambda^{-1}_U(\sum_{\cc C}\f a_i)=\sum_{\cc C}\lambda_U^{-1}(\f a_i)\in\f G$
with $\lambda_U$ an \env S of $U$. By assumtion there is a finitely
generated subobject $\f b$ of $\f a$ such that $\f b\in\f G$. Then
there is $i_0$ such that $\f b\subseteq\lambda_U^{-1}(\f a_{i_0})$.
One has
   $$\qq US=\qq{\f b}S\subseteq\qq{(\lambda_U^{-1}(\f a_{i_0}))}S
     \subseteq\f a_{i_0}\subseteq\qq US.$$
So $\f a_{i_0}=\qq US$ that implies $\qq US\in\fg\cc C/\cc S$.

Further, the isomorphism
   \begin{multline*}
   \qc U{\sum\nolimits_{\cc C/\cc S}C_i}\iso
   \q{\qq US}{\sum\nolimits_{\cc C/\cc S}C_i}\\
   \iso\lp\q{\qq US}{C_i}\iso\lp\qc U{C_i}
   \iso\qc U{\sum\nolimits_{\cc C}C_i}
   \end{multline*}
is functorial in $U$, and so $\sum_{\cc C/\cc S}C_i=\sum_{\cc C}C_i$.

$(1)\Rightarrow (6)$: The direct limit $\lp_I C_i$ may be described as a
quotient object of a coproduct $\ps_IC_i$. To be precise, let $R$ be the
subset of $I\times I$ consisting of pairs $(i,j)$ with $i\le j$ and
for each $S\subseteq R$ we put
   $$C_S=\sum_{(i,j)\in S}\im(u_i-u_j\gamma_{ij})\subseteq\ps C_i,$$
where $u_i:C_i\to\ps C_i$ is the canonical monomorphism for $i\in I$,
$\gamma_{ij}:C_i\to C_j$ is the canonical morphism for $i\le j$. Then
$\lp C_i=\ps C_i/C_R=\ps C_i/\sum C_S$, where $S$ runs over all finite
subsets of $R$. By assumtion both $\ps C_i$ and $\sum C_i$ are \cl S,
and therefore $\lp C_i$, as a quotient object of \cl S objects, is \tr S.

$(5)\Rightarrow (1)$: Suppose that $X$ is \cl S. Write
$X=\sum_{\cc C/\cc S}X_i$ as a direct union of \cl S subobjects.
Then $X/\sum_{\cc C}X_i=t(X/\sum_{\cc C}X_i)=t(\lp_{\cc C}X/X_i)=\lp t(X/X_i)=0$
since $t(X/X_i)=0$. Thus $\sum_{\cc C/\cc S}X_i=\sum_{\cc C}X_i$.

In turn, let $D\in\fg\cc C/\cc S$. Write $D=\sum_{\cc C}D_i$ as a
directed sum of $D_i\in\fg\cc C$. Then $D=\qq DS=
\sum_{\cc C/\cc S}\qq{(D_i)}S$ whence $D=\qq{(D_{i_0})}S$ for some $i_0$.
\end{proof}

\begin{prop}\label{fin}
Let $\cc C$ be a locally finitely presented Grothendieck
category with the family of finitely presented generators
$\cc U\subseteq\fp\cc C$. Then the following conditions are equivalent:

\begin{itemize}
\item[$(1)$] $\cc S$ is of finite type.

\item[$(2)$] For any $U\in\cc U$ the natural morphism
      $\lp\qc U{C_i}\to\qc U{\lp_{\cc C/\cc S}C_i}$
      induced by the \env S $\lambda$ of $\lp_{\cc C}C_i$
      is an isomorphism.

\item[$(3)$] $\qq{\cc U}S=\{\qq US\}_{U\in\cc U}$
      is a family of $\cc C/\cc S$-finitely presented
      generators for $\cc C/\cc S$.

\item[$(4)$] If $C$ is $\cc C$-finitely presented, then $\qq CS$
      is $\cc C/\cc S$-finitely presented.
\end{itemize}

Thus $\cc C/\cc S$ is a locally finitely presented Grothendieck category
and in this case, any $\cc C/\cc S$-finitely presented object
$D$ is a localization $\qq CS$ of some $C\in\fp\cc C$.
\end{prop}

\begin{proof}
It suffices to note that for each $C\in\fp\cc C$ the representable
functor $(C,-)$ commutes with direct limits (see Theorm~\ref{finpr})
and there is a presentation
   $$\pp{\ps_{i=1}^nU_i}{\ps_{j=1}^mU_j}C$$
of $C$ by objects from $\cc U$; and then our proof literally
repeats Theorem~\ref{pref}.

In turn, if $D\in\fp\cc C/\cc S$ there is an epimorphism $\eta:\qq BS\to D$
with $B\in\fp\cc C$; then $\kr\eta\in\fg\cc C/\cc S$. According
to~\cite[Lemma~2.13]{He} we can choose $C\subseteq B$ such that
$C\in\fg\cc C$ and $\qq CS=\kr\eta$. Hence $D=\qq{(B/C)}S$.
\end{proof}

Now let us consider a localizing subcategory $\cc S$ of the module category
$\mc A$ with $\cc A=\{P_i\}_{i\in I}$ a ring. Then the family
$\qq {\cc A}S=\{\qq PS\}_{P\in\cc A}$ generates $\mc A/\cc S$ and we call it
the {\em ring of quotients\/} of $\cc A$ with respect to $\cc S$.

There is a naturally defined functor
$j:\mc A/\cc S\to\Mod\qq{\cc A}S=(\qq{\cc A}S^{\text{op}},\text{Ab})$:
   $$\Hom_{\qq{\cc A}S}(\qq PS,j(M)):=\homa PM$$
for every $P\in\cc A$ and $M\in\mc A/\cc S$. The next assertion is an
analog of the Walker and Walker Theorem (see~\cite[Proposition~XI.3.4]{St}).

\begin{prop}\label{greb}
The functor $j:\mc A/\cc S\to\Mod\qq{\cc A}S$ is an equivalence \ifff
the localizing subcategory $\cc S$ is of finite type and coexact.
\end{prop}

\begin{proof}
Suppose $j$ is an equivalence. Then every $\qq PS$ finitely generated and
projective in $\mc A/\cc S$. By Proposition~\ref{osya} $\cc S$ is
coexact and by Proposition~\ref{fin} it is of finite type.

Conversely, let $\cc S$ be of finite type and coexact.
By Proposition~\ref{osya} every $\qq PS$ is projective in $\mc A/\cc S$
and by Proposition~\ref{fin} it is finitely generated in $\mc A/\cc S$.
Now our assertion follows from Proposition~\ref{Freyd}.
\end{proof}

\subsection{Left exact functors}

Let $\cc C$ be a locally finitely presented Grothendieck category.
By Theorem~\ref{main} the functor $T:\cc C\to\mc A$ with
$\cc A=\{h_X\}_{X\in\fp\cc C}$ induces an equivalence of $\cc C$ and
$\mc A/\cc S$ where $\cc S$ is some localizing subcategory
of $\mc A$. Let $\f F=\{\f F^X\}_{X\in\fp\cc C}$ be the
respective Gabriel topology on $\cc A$.

Denote by $\cc L$ the subcategory of $\mc A$ consisting of $L\in\mc A$
such that for any $x\in L(X)$, $X\in\fp\cc C$, there exists an epimorphism
$f:Y\to X$ such that $L(f)(x)=0$. It is directly verified that $\cc L$
is a localizing subcategory.

\begin{lem}\label{br1}
If $\p{X'\bl p}{X\bl f}{X''}$ is an exact sequence in $\fp\cc C$,
then $\coker(Tf)\in\cc L$.
\end{lem}

\begin{proof}
Let $Y\in\fp\cc C$ and $y\in\coker(Tf)(Y)$; then there is a morphism
$g\in\qc Y{X''}$ such that $u_y=r\circ Tg$, where $r:h_{X''}\to\coker(Tf)$
is the canonical epimorphism and $u_{y,Y}(1_Y)=y$. Consider in $\cc C$
the commutative diagram
   $$\begin{CD}
       \dec X{X''}Y@>f'>>Y\\
       @Vg'VV@VVgV\\
       X@>>f>X''\\
     \end{CD}$$
in which $f'$ is an epimorphism. Because $Y\in\fp\cc C$, there is
a finitely generated subobject $Z'$ of $\dec X{X''}Y$ such that
$f'(Z')=Y$. There exists an epimorphism $h:Z\to Z'$ with $Z\in\fp\cc C$.
So $f'h$ is an epimorphism. It is easily seen that $\coker(Tf)(f'h)(y)=0$,
hence $\coker(Tf)\in\cc L$.
\end{proof}

Recall that the functor $M\in\mc A$, $\cc A=\{h_X\}_{X\in\fp\cc C}$, is
{\em left exact\/} if for any exact sequence
   $$\pp{X'}X{X''}$$
in $\fp\cc C$ the sequence of abelian groups
   $$\ii{M(X'')}{M(X)}{M(X')}$$
is exact.

\begin{prop}\label{br2}
Every left exact functor $M\in\mc A$ is \cl L.
\end{prop}

\begin{proof}
Let $M:\fp\cc C\to\text{Ab}$ be a contravariant left exact functor, $X\in\fp\cc C$
and $x\in M(X)$ such that $\kr u_x\in\f F^X$. Here $u_x$ denotes a unique
morphism such that $u_{x,X}(1_X)=x$. Then there exists an epimorphism $f:Y\to X$
in $\fp\cc C$ such that $M(f)(x)=0$. But $M(f)$ is a monomorphism of abelian
groups and thus $x=0$. So $M$ is \tr L.

Now let $\f a\in\f F^X$ and $g:\f a\to M$ a morphism in $\mc A$. There exists
an epimorphism $f:Y\to X$ such that $\im(Tf)\subseteq\f a$.
By Lemma~\ref{br1} $\im(Tf)\in\f F^X$. One has a morphism
$g\circ Tf=t:h_Y\to M$. Since $\kr f\in\fg\cc C$, there exists an epimorphism
$Z\to\kr f$ and thus one gets an exact sequence
   $$\pp {Z\bl p}{Y\bl f}X.$$
We have the following commutative diagram of abelian groups:
   $$\begin{CD}
      0@>>>h_Y(X)@>h_Y(f)>>h_Y(Y)@>h_Y(p)>>h_Y(Z)\\
      @.@Vt_XVV@VVt_YV@VVt_ZV\\
      0@>>>M(X)@>>M(f)>M(Y)@>>M(p)>M(Z)\\
     \end{CD}$$
So,
   $$(M(p)t_Y)(1_Y)=(t_Zh_Y(p))(1_Y)=(g_Z(Tf)_Zh_Y(p))(1_Y)=g_Z(fp)=0.$$
Then there exists an element $x\in M(X)$ such that $M(f)(x)=t_Y(1_Y)=
(gTf)_Y(1_Y)=g_Y(f)$. Thus $u_x|_{\im(Tf)}=g|_{\im(Tf)}$. Since
${\im(Tf)}\in\f F^X$ and $M$ is \tr L, it follows that $u_x|_{\f a}=g$.
Thus $u_x:h_X\to M$ is a morphism prolonging $g$.
The uniqueness of $u_x$ follows from the fact that $M$ is \tr L. So $M$ is \cl L.
\end{proof}

Define $\Lex((\fp\cc C)^{\text{\rm op}},\text{\rm Ab})$ to be the category of
contravariant left exact functors from $\fp\cc C$ to Ab.
So, we are now in a position to prove the following result.

\begin{thm}[Breitsprecher~\cite{Br}]\label{br3}
Let $\cc C$ be a locally finitely presented Grothendieck category. The
representation functor
$T=\qc-?:\cc C\to((\fp\cc C)^{\text{\rm op}},\text{\rm Ab})$ defined by
$X\mapsto (-,X)$ is a natural equivalence between $\cc C$ and
$\Lex((\fp\cc C)^{\text{\rm op}},\text{\rm Ab})$.
\end{thm}

\begin{proof}
Our assertion would be proved, if we showed that $\cc S=\cc L$. Because
for every $X\in\cc C$ the functor $T(X)$ is left exact, in view of
Proposition~\ref{br2}, it follows that $T(X)$ is \cl L, and so $\cc L\subseteq\cc S$.
Conversely, let $M\in\cc S$. We can choose a projective presentation of $M$
   $$\pp{\ps h_{Y_j}\bl g}{\ps h_{X_i}}M.$$
By Proposition~\ref{fin} $\cc S$ is of finite type, and so any coproduct of
\cl S objects is an \cl S object. Therefore there exists a morphism $f\in\Mor\cc C$
such that $Tf=g$. Furthermore, $f$ is an epimorphism in $\cc C$ since $\qq MS=0$.
Thus, without loss of generality, we can assume that for every $M\in\cc S$
there is an exact sequence
   $$\pp{TY\bl{Tf}}{TX}M,$$
where $f$ is a $\cc C$-epimorphism.

According to~\cite[Lemma~5.7]{Kr3} $f$ is a direct limit
$f_{\alpha}:Y_{\alpha}\to X_{\alpha}$ of epimorphisms $f_{\alpha}$ in $\fp\cc C$,
and so $M\iso\lp\coker (Tf_{\alpha})$. By Lemma~\ref{br1} every
$\coker (Tf_{\alpha})\in\cc L$, hence $M\in\cc L$.
\end{proof}

The next result describes localizing subcategories of finite type of
a locally coherent Grothendieck category.

\begin{thm}\label{fincoh}
Let $\cc C$ be a locally coherent Grothendieck category
with the family of coherent generators $\cc U\subseteq\coh\cc C$.
Then the following conditions are equivalent:

\begin{itemize}
\item[$(1)$] $\cc S$ is of finite type.

\item[$(2)$] $\cc S$ is of prefinite type.

\item[$(3)$] The torsion functor $t$ commutes with direct limits.

\item[$(4)$] Any $\cc C$-direct limit of \cl S objects is \tr S.

\item[$(5)$] For any $U\in\cc U$ the natural morphism
      $\lp\qc U{C_i}\to\qc U{\lp_{\cc C/\cc S}C_i}$
      induced by the \env S $\lambda$ of $\lp_{\cc C}C_i$
      is an isomorphism.

\item[$(6)$] $\qq{\cc U}S=\{\qq US\}_{U\in\cc U}$
      is a family of $\cc C/\cc S$-coherent
      generators for $\cc C/\cc S$.

\item[$(7)$] If $C$ is $\cc C$-coherent, then $\qq CS$ is
      $\cc C/\cc S$-coherent.

\item[$(8)$] $\f G$ has a basis of coherent objects.
\end{itemize}

Thus $\cc C/\cc S$ is a locally coherent Grothendieck category
and in this case, any $\cc C/\cc S$-coherent object
$D$ is a localization $\qq CS$ of some $C\in\coh\cc C$.
\end{thm}

\begin{proof}
$(1)\Rightarrow (5)$: It follows from Proposition~\ref{fin}.

$(5)\Rightarrow (6)$: By Proposition~\ref{fin}
$\qq{\cc U}S\subseteq\fp\cc C/\cc S$. Let us show that
for any $U\in\cc U$ the object $\qq US$ is coherent.
Let $C\subseteq\qq US$ be $\cc C/\cc S$-finitely generated
subobject of $\qq US$. From~\cite[Lemma~2.13]{He} it
follows that there is a $\cc C$-finitely generated
subobject $A\subseteq U$ of $U$ such that $\qq AS=C$.
Since $U$ is coherent, the object $A\in\coh\cc C$.
By Proposition~\ref{fin} $\qq AS\in\fp\cc C/\cc S$.

$(6)\Rightarrow (7)$: Easy.

$(7)\Rightarrow (1)$: $\cc S$ is of finite type by
Proposition~\ref{fin}.

$(2)\Leftrightarrow (3)\Leftrightarrow (4)$:
It follows from Theorem~\ref{pref}.

$(2)\Leftrightarrow (8)$: Since any finitely generated
subobject of a coherent object is coherent, our assertion
follows from Theorem~\ref{pref}.

$(1)\Rightarrow (2)$: Straightforward.

$(3)\Rightarrow (1)$: It follows from~\cite[Lemma~2.4]{Kr1}.
\end{proof}

As in~\cite{He} for an arbitrary subcategory $\cc X$ of $\cc C$,
denote by $\vec\cc X$ the subcategory of $\cc C$ consisting of direct
limits of objects in $\cc X$. Herzog~\cite{He} and Krause~\cite{Kr1}
observed that localizing subcategories of finite type in a locally
coherent category are determined by Serre subcategories
of the abelian subcategory $\coh\cc C$. Namely it holds the following.

\begin{thm}[Herzog and Krause]\label{hk}
Let $\cc C$ be a locally coherent Grothendieck category.
There is a bijective correspondence between Serre subcategories
$\cc P$ of $\coh\cc C$ and localizing subcategories $\cc S$
of $\cc C$ of finite type. This correspondence is given by
the functions
   \begin{gather*}
     \cc P\longmapsto\vec\cc P\\
     \cc S\longmapsto\coh\cc S=\cc S\cap\coh\cc C
   \end{gather*}
which are mutual inverses.
\end{thm}

Later on, given a Serre subcategory $\cc P$ of $\coh\cc C$, the
$\vec\cc P$-torsion functor will be denoted by $\qq tP$.

Recall that an object $C\in\cc C$ of a Grothendieck category $\cc C$ is
{\em noetherian\/} if any subobject of $C$ is finitely generated. $\cc C$
is called {\em locally noetherian\/} if it has a family of noetherian
generators. In that case, the relations
   $$\fg\cc C=\fp\cc C=\coh\cc C$$
hold. If in addition $\cc C$ is a locally finitely generated, it is locally
coherent and any localizing subcategory $\cc S$ of $\cc C$ is of finite type.
Any $X\in\cc S$ is a direct union $\sum_{i\in I}X_i$ of objects
$X_i\in\fg\cc C\cap\cc S=\coh\cc C\cap\cc S=\coh\cc S$. Furthermore, any
quotient category $\cc C/\cc S$ is locally noetherian with a family of
noetherian generators $\{\qq CS\}_{C\in\coh\cc C}$.

\subsection{The Ziegler topology}

The study of pure-injective (= algebraically compact) modules over
different classes of rings plays an important role in the theory of
rings and modules. It goes back to Cohn and has been further developed
by various mathematicians. Because pure-injective modules can be
defined, using some condition of solvability (in this module) of linear
equations systems, many problems of (algebraic!) structure of pure-injective
modules admit reformulations using concepts of the model theory. It is
such an approach that led Ziegler~\cite{Zi} to construction of the
topological space (``The Ziegler spectrum'') whose points are
indecomposable pure-injective modules. Recently Herzog~\cite{He} and
Krause~\cite{Kr1} have proposed an algebraic definition of the Ziegler
spectrum.

Let $\cc C$ be a Grothendieck category with the family of generators
$\cc U$. We denote by $\zg\cc C$ the set of isomorphism classes of
indecomposable injective objects in $\cc C$ and call $\zg\cc C$ the
{\em Ziegler spectrum\/} of $\cc C$. The fact that $\zg\cc C$ forms
a set follows from that any indecomposable injective object in $\cc C$
occurs as the injective envelope of some $\cc U$-finitely generated
object $X\in\fg_{\cc U}\cc C$ and $\fg_{\cc U}\cc C$ is skeletally
small. It will be convenient to identify each isomorphism class with
a representative belonging to it. If $\cc S$ is a localizing
subcategory in $\cc C$, the assignment $X\mapsto E(X)$
induces injective maps $\zg\cc S\to\zg\cc C$ and
$\zg\cc C/\cc S\to\zg\cc C$. We consider both maps as identifications.
They satisfy $\zg\cc S\cup\zg\cc C/\cc S=\zg\cc C$ and
$\zg\cc S\cap\zg\cc C/\cc S=\emptyset$~\cite[Corollaire~III.3.2]{Ga}.

Now let $\cc C$ be a locally coherent category, i.e.\ $\cc U\subseteq\coh\cc C$.
To an arbitrary subcategory $\cc X\subset\coh\cc C$, we associate
the subset of $\zg\cc C$
   $$\cc O(\cc X)=\{E\in\zg\cc C\mid\qc CE\ne 0\text{ for some }C\in\coh\cc C\}.$$
If $\cc X=\{C\}$ is singleton, we abbreviate $\cc O(\cc X)$ to
$\cc O(C)$; thus $\cc O(\cc X)=\cup_{C\in\cc X}\cc O(C)$.

We restrict the discussion to subcategories of the form $\cc O(\cc S)$
where $\cc S\subseteq\coh\cc C$ is a Serre subcategory. In that case,
   $$\cc O(\cc S)=\{E\in\zg\cc C\mid\qq tS(E)\ne 0\}.$$

\begin{thm}[Herzog~\cite{He} and Krause~\cite{Kr1}]\label{krut}
For a locally coherent Gro\-then\-dieck category $\cc C$ the following
assertions hold:

$(1)$ The collection of subsets of $\zg\cc C$
      $$\{\cc O(\cc S)\mid\cc S\text{ is a Serre subcategory}\}$$
satisfies the axioms for the open sets of a topological space on $\zg\cc C$.
This topological space we call the Ziegler spectrum of $\cc C$ too.

$(2)$ There is a bijective inclusion preserving correspondence between
the Serre subcategories $\cc S$ of $\coh\cc C$ and the open subsets $\cc O$
of $\zg\cc C$. This correspondence is given by the functions
   \begin{gather*}
     \cc S\mapsto\cc O(\cc S)\\
     \cc O\mapsto\cc S_{\cc O}=\{C\in\coh\cc C\mid\cc O(C)\subseteq\cc O\}
   \end{gather*}
which are mutual inverses.
\end{thm}

Recall that a topological space $\cc X$ is {\em quasi-compact\/} provided
that for every family $\{\cc O_i\}_{i\in I}$ of open subsets
$\cc X=\cup_{i\in I}\cc O_i$ implies $\cc X=\cup_{i\in J}\cc O_i$ for some
finite subset $J$ of $I$. A subset of $\cc X$ is quasi-compact if it is
quasi-compact with respect to the induced topology.

By~\cite[Corollary~3.9]{He} and~\cite[Corollary~4.6]{Kr1} an open subset
$\cc O$ of $\zg\cc C$ is quasi-compact \ifff it is one of the basic open
subsets $\cc O(C)$ with $C\in\coh\cc C$.

Serre subcategories $\cc S$ of $\coh\cc C$ arise in the following natural
way. An object $M\in\cc C$ is {\em coh-injective\/} if $\Ext^1_{\cc C}(C,M)=0$
for each $C\in\coh\cc C$. Then the subcategory
   $$\cc S_M=\{C\in\coh\cc C\mid\qc CM=0\}$$
generated by $M$ is Serre. Moreover, every Serre subcategory of $\coh\cc C$
arises in this fashion~\cite[Corollary~3.11]{He}.

\begin{exs}{\rm
Here are some examples of the Ziegler-closed subsets:

(1) Let $\cc C$ be a locally coherent Grothendieck category; then the functor
category $((\coh\cc C)^{\text{\rm op}},\text{\rm Ab})$ is locally coherent.
By Theorem~\ref{br3} the category $\cc C$ is equivalent to
$\Lex((\coh\cc C)^{\text{\rm op}},\text{\rm Ab})=((\coh\cc C)^{\text{\rm op}},\text{\rm Ab})/\vec\cc L$,
where $\cc L$ is a Serre subcategory of $\coh((\coh\cc C)^{\text{\rm op}},\text{\rm Ab})$
consisting of objects isomorphic to $\coker(\mu,-)$ for some epimorphism
$\mu:A\to B$ in $\coh\cc C$. Thus $\zg\cc C$ is closed in
$\zg((\coh\cc C)^{\text{\rm op}},\text{\rm Ab})$.

(2) If $\cc C=\rc$, the Ziegler spectrum of $\rc$
    $$\zg\rc=\{\ten Q-\mid Q_A\text{ is an indecomposable pure-injective module}\}.$$
We notice that $\zg\rc$ is quasi-compact since $\zg\rc=\cc O(\ten A-)$.
Prest, Rothmaler and Ziegler have shown~\cite[Corollary~4.4]{PRZ}
(see also~\cite[Theorem~2.5]{GG}) that the ring $A$ is right coherent if
and only if the set
   $$I_{\text{inj}}=\{\ten Q-\in\zg\rc\mid Q\text{ is an injective module}\}$$
is closed in $\zg\rc$.

(3) If $A$ is left coherent, from~\cite[Theorem~2.4]{GG2} it follows that
   $$I_{\text{flat}}=\{\ten Q-\in\zg\rc\mid Q\text{ is a flat module}\}$$
is closed in $\zg\rc$.

Recall that the ring $A$ is {\em weakly Quasi-Frobenius\/} if the functor
$\hm-A$ gives a duality of categories $\lfp$ and $\rfp$. Weakly Quasi-Frobenius
rings are described as (two-sided) coherent absolutely pure
rings~\cite[Theorem~2.11]{GG2}. Over weakly Quasi-Frobenius rings
$I_{\text{inj}}=I_{\text{flat}}$. And backwards,
if $A$ is a (two-sided) coherent ring and $I_{\text{inj}}=I_{\text{flat}}$,
then $A$ is a weakly Quasi-Frobenius ring (see~\cite[Corollary~2.12]{GG2}).

(4) Let $\rho:A\to B$ be a ring homomorphism; then $\rho$ induces the exact
functor $\rho^*:\coh\rc\to\coh\cc C_B$. If $\rho$ is an epimorphism of rings,
then the map $M_B\mapsto M_A$ induces a homeomorphism
$\zg\cc C_B\to\zg\rc\backslash\cc O(\cc S)$ with $\cc S=\kr\rho^*$. Thus
$\zg\cc C_B$ is closed in $\zg\rc$~\cite[Corollary~9]{Pr}.
}\end{exs}

\section{Categories of generalized modules}

The most important notions and properties of the category
of generalized $A$-modules $\rc$ can be easily
extended to an arbitrary functor category $\mc A=\fc A$.
The model-theoretic background can be found in~\cite{Bur,Pr2}.
In this section, for the most part, we adhere to the reference~\cite{He}.

\subsection{Tensor products}

Let $\mc A$ be the category of right $\cc A$-modules with the family of
finitely generated projective generators $\cc A=\{P_i\}_{i\in I}$. According
to Proposition~\ref{Freyd} $\mc A\iso\fc A$. We refer to the functor category
$\ffc A$ as a {\em category of left $\cc A$-modules\/} and denote it by
$\mmc A$. It is a locally finitely presented Grothendieck category with
the family of finitely generated projective generators
$\cc A^{\rm op}=\{h^P=\homa P-\}_{P\in\cc A}$.

\begin{prop}\cite[Theorem~III.6.3]{Pop}
Let $M:\cc A^{\rm op}\to\text{\rm Ab}$ ($N:\cc A\to\text{\rm Ab}$) be a
right (left) $\cc A$-module. Then, a unique functor
$\tena M-:\mmc A\to\text{\rm Ab}$ ($\tena-N:\mc A\to\text{\rm Ab}$)
exists such that:

$(1)$ There are functorial isomorphisms $\tena M{h^P}\iso M(P)$ and
      $\tena{h_P}N\iso N(P)$ for $P\in\cc A$.

$(2)$ $\tena M-$ and $\tena-N$ have right adjoints.
\end{prop}

Note also that the {\em tensor product functor\/} $\tena M-$ is right
exact and commutes with direct limits.

So let $\mc A=\fc A$ ($\mcl A=\ffc A$) be the category of
right (left) $\cc A$-modules with $\cc A=\{P_i\}_{i\in I}$,
$\fpl A$ ($\fpr A$) be the category of finitely presented
left (right) $\cc A$-modules. By definition every
$M\in\fpr A$ has a projective presentation
   \begin{equation}\label{www}
     \pp{\ps_{j=1}^nh_{P_j}}{\ps_{k=1}^mh_{P_k}}M
   \end{equation}
in $\mc A$. Similarly, every $M\in\fpl A$ has a presentation
   $$\pp{\ps_{j=1}^nh^{P_j}}{\ps_{k=1}^mh^{P_k}}M$$
in $\mmc A$.

Denote by $\rcc=(\fpl A,\text{Ab})$ ($\lcc=(\fpr A,\text{Ab})$)
the category of additive covariant functors defined on
$\fpl A$ ($\fpr A$) and call $\rcc$ ($\lcc$) the {\em category
of generalized right (left) $\cc A$-modules.} They naturally
extend the respective categories of generalized $A$-modules
for which $\cc A=\{A\}$ is a ring. Because $\fpl A$ is
closed under cokernels, $\rcc$ is a locally coherent Grothendieck
category. Finitely generated projective objects of
$\rcc$ are precisely of the form $\{(M,-)\}_{M\in\fpl A}$
and this family generates $\rcc$.
The tensor product functor $\tena ?-:\mc A\to\rcc$ defined
by the rule $\qq MA\mapsto\tena M-$ is fully faithful and
right exact.

Each finitely presented (= coherent) generalized module
$C\in\rcc$ has a projective presentation
   $$\pp{(K,-)\bl{(f,-)}}{(L,-)}C$$
with $K$, $L\in\fpl A$. In particular, if $M\in\fpr A$, then
$\rcohh{\tena M-}$. Indeed, consider projective presentation~\eqref{www}
of $M$. As tensoring is right exact, this gives an exact sequence
in $\rcc$
   $$\pp{\ps_{j=1}^n\tena{h_{P_j}}-}{\ps_{k=1}^m\tena{h_{P_k}}-}M$$
which is a presentation of $\tena M-$ in $\rcc$ by finitely
generated projective objects since
$\ps\tena{h_{P_j}}-\iso (\ps P_j,-)$.

\begin{lem}
An object $E\in\rcc$ is $\coh\rcc$-injective \ifff it is isomorphic to
one of the functors $\tena M-$ where $M$ is a right $\cc A$-module.
\end{lem}

\begin{proof}
Herzog~\cite[Proposition~2.2]{He} has shown that $E$ is $\coh\rcc$-injective
\ifff it is right exact. Therefore the functor $\tena M-$ is coh-injective.

Let $E$ be coh-injective. Define $M_{\cc A}$ by putting
$\homa PM=\rqq{(h^P,-)}E=E(h^P)$. Now our assertion is proved similar
to~\cite[Proposition~IV.10.1]{St}.
\end{proof}

Thus the category $\mc A$ of right $\cc A$-modules can be considered as
the subcategory of $\coh\rcc$-injective objects of the category $\rcc$.

In order to describe the points of the Ziegler spectrum $\zg\rcc$ of
$\rcc$, recall that \ses
   $$\les{L\bl\mu}{M\bl\delta}N$$
of an arbitrary locally finitely presented Grothendieck category $\cc C$
is {\em pure-exact\/} provided that the sequence
   $$\les{\qc XL}{\qc XM}{\qc XN}$$
is exact for all $X\in\fp\cc C$. In this case, $\mu$ is called a
{\em pure-monomorphism.} An object $Q\in\cc C$ is said to be {\em pure-injective}
if every pure-exact sequence with first term $Q$ splits.

\begin{prop}\label{piano}
For \ses $\epsilon:\es LMN$ in $\cc C$, the following are equivalent:

$(1)$ $\epsilon$ is pure-exact in $\cc C$.

$(2)$ $\epsilon$ is a direct limit of split exact sequences
      $\es{L_i}{M_i}{N_i}$ in $\cc C$.
\end{prop}

\begin{proof}
Write $N=\lp_I N_i$ as a direct limit of finitely presented objects $N_i$.
For every $i\in I$ consider the following commutative diagram
   $$\begin{CD}
     \epsilon_i:0@>>>L@>>>M_i@>\delta_i>>N_i@>>>0\\
     @.@|@V{\psi_i}VV@VV{\phi_i}V\\
     0@>>>L@>>>M@>>\delta>N@>>>0\\
     \end{CD}$$
in which the right square is pullback and $\phi_i$ the canonical morphism.
Since for $i\le j$ the relations
   $$\delta\psi_i=\phi_i\delta_i=\phi_j(\phi_{ij}\delta_i)$$
hold, there exists a unique $\psi_{ij}:M_i\to M_j$ such that
$\psi_i=\psi_j\psi_{ij}$. Clearly that the system $\{M_i,\psi_{ij}\}_I$
is direct and $\lp\epsilon_i=\epsilon$.
By hypothesis there exists a morphism $f:N_i\to M$ such that $\delta f=\phi_i$,
and hence there exists $g:N_i\to M_i$ such that $\delta_i g=1_{M_i}$, i.e.\
each sequence $\epsilon_i$ splits.

Conversely, if each $\epsilon_i$ splits, then the sequence
   $$\les{\homa X{L_i}}{\homa X{M_i}}{\homa X{N_i}}$$
is exact. The fact that $\homa X-$, $X\in\fp\cc C$, commutes with direct
limits implies~(1).
\end{proof}

\begin{cor}
The sequence $\epsilon:\es LMN$ of right $\cc A$-modules is pure-exact \ifff
the $\rcc$-sequence $\tena\epsilon-:\es{\tena L-}{\tena M-}{\tena N-}$ is exact.
\end{cor}

\begin{proof}
As tensoring commutes with direct limits, the necessary condition follows
from the preceding Proposition. Conversely, if $\tena\epsilon-$ is exact,
then for any $X\in\fpr A$ one has the following exact sequence
   \begin{multline*}
     0\lra{}(\tena X-,\tena L-)\lra{}(\tena X-,\tena M-)\lra{}\\
     (\tena X-,\tena N-)\lra{}\Ext^1(\tena X-,\tena L-).
   \end{multline*}
Because $\rcohh{\tena X-}$ and the object $\tena L-$ is $\coh\rcc$-injective,
it follows that $\Ext^1(\tena X-,\tena L-)=0$.
\end{proof}

Notice that as tensoring commutes with direct limits, a monomorphism $\mu:L\to M$
is pure \ifff $\tena\mu X$ is a monomorphism for any right $\cc A$-module $X$.

The next result is proved similar to~\cite[Proposition~4.1]{He}.

\begin{prop}
An object $E\in\rcc$ is an injective object \ifff it is isomorphic to one
of the functors $\tena Q-$ where $Q$ is a pure-injective right $\cc A$-module.
\end{prop}

Thus the points of the Ziegler spectrum $\zg\rcc$ of $\rcc$ are represented
by the pure-injective indecomposable right $\cc A$-modules. Note that every
indecomposable injective right $\cc A$-module $E_{\cc A}$ is pure-injective
and hence $\tena E-$ is a point of $\zg\rcc$.

It is easy to see that $\zg\rcc=\cup_{P\in\cc A}\cc O(\tena P-)$ is a union
of basic open quasi-compact subsets $\cc O(\tena P-)$.
Prest has shown~\cite[Example~1.5]{Pr2} that $\zg\rcc$ need not
be compact, in contrast to the case $\cc A=\{A\}$ where
$A$ is a ring, the whole spase need not be basic open.

\subsection{Auslander-Gruson-Jensen Duality}

Let $A$ be an arbitrary ring. Gruson and Jensen~\cite{GJ}
and Auslander~\cite{Au2} proved that there is a duality
$D:(\coh\lc)^{\text{op}}\iso\coh\rc$ between the respective
subcategories of the coherent objects of $\lc$ and $\rc$.
It is not hard to show that the same holds for an arbitrary
ring $\cc A=\{P_i\}_{i\in I}$.

Namely let the functor $D:(\coh\lcc)^{\text{op}}\to\rcc$ be given by
   $$(DC)({}_{\cc A}N)=\lqq C{\tena-N}$$
where $\lcohh C$ and $N\in\fpl A$. If $\delta:B\to C$ is a
morphism in $\coh\lcc$, then
   $$D(\delta)_N:(DC)({}_{\cc A}N)\lra{}(DB)({}_{\cc A}N)$$
is put to be $\lqq{\delta}{\tena-N}$.

We claim that $D$ is exact. Indeed, if
   $$\les {A\bl\alpha}{B\bl\beta}C$$
is \ses in $\coh\lcc$, then because $\tena-N$ is $\coh\lcc$-injective,
the sequence
   $$\les {(C,\tena-N)\bl{(\beta,\tena-N)}}{(B,\tena-N)\bl{(\alpha,\tena-N)}}{(A,\tena-N)}$$
is exact for each $N\in\fpl A$. Therefore the sequence
   $$\les{DC\bl{D\beta}}{DB\bl{D\alpha}}{DA}$$
is exact.

By construction of $D$,
   $$D(\qq MA,-)({}_{\cc A}X)=((\qq MA,-),\tena-X)\iso\tena MX.$$
Thus $D(\qq MA,-)\iso\tena M-$. Given $\lcohh C$, consider a projective
presentation
   $$\pp{(\qq NA,-)}{(\qq MA,-)}C$$
of $C$. Applying the exact functor $D$, we obtain an exact sequence
   $$\ii{DC}{\tena M-}{\tena N-}.$$
Thus $DC$ is a coherent object of $\rcc$ and therefore the functor $D$
has its image in $\coh\rcc$.

\begin{thm}[Auslander, Gruson and Jensen] The functor
$D:(\coh\lcc)^{\text{\rm op}}\to\coh\rcc$ defined above constitutes a duality
between the categories $\coh\lcc$ and $\coh\rcc$. Furthermore, for
$\qq MA\in\fpr A$ and ${}_{\cc A}N\in\fpl A$ we have that
   $$D(\qq MA,-)\iso\tena M-\text{ \ and \ }D(\tena-N)\iso ({}_{\cc A}N,-).$$
\end{thm}

\begin{proof}
Since $D(\tena-N)({}_{\cc A}M)\iso(\tena-N,\tena-M)\iso\homa{{}_{\cc A}N}{{}_{\cc A}M}$,
it follows that $D(\tena-N)\iso({}_{\cc A}N,-)$. Similarly, we can define
the functor $D':(\coh\rcc)^{\text{op}}\to\coh\lcc$ in the other direction.
Both of the compositions $DD'$ and $D'D$ are exact functors that are
equivalences on the respective categories of finitely generated projective
objects. Therefore they are both natural equivalences.
\end{proof}

Because the category $\coh\lcc$ has enough projectives, the duality
gives the following.

\begin{prop}[Auslander] The category $\coh\rcc$ has enough injectives
and they are precisely the objects of the form $\tena M-$ where
$\qq MA\in\fpr A$.
\end{prop}

Thus every coherent object $\rcohh C$ has both a projective
presentation in $\rcc$
   \begin{equation}\label{7.1}
    \pp{({}_{\cc A}K,-)}{({}_{\cc A}L,-)}C
   \end{equation}
and an injective presentation in $\coh\rcc$
   \begin{equation}\label{7.2}
    \ii C{\tena M-}{\tena N-}
   \end{equation}
where $K$, $L\in\fpl A$ and $M$, $N\in\fpr A$.

We conclude the section by Herzog's Theorem. Let
$\cc S\subseteq\coh\lcc$; then the subcategory
   $$D\cc S=\{DC\mid C\in\cc S\}$$
is Serre in $\coh\rcc$ and the restriction to $\cc S$
of the duality $D$ gives a duality $D:\cc S^{\text{op}}\to D\cc S$.
By Theorem~\ref{krut} the map $\cc O(\cc S)\mapsto\cc O(D\cc S)$
induced on the open subsets of the Ziegler spectrum is an
iclusion-preserving bijection. Similar to~\cite{He}, it
is easily shown that the functor $D$ induces an isomorphism
of abelian groups
   $${}_{\lcc/\vec\cc S}(A,B)\iso{}_{\rcc/\overrightarrow{D\cc S}}(DB,DA)$$
where $A$ and $\lcohh B$ and the assignment given by
$\qq AS\mapsto(DA)_{D\cc S}$ is functorial. Thus we have
the following.

\begin{thm}[Herzog]\label{herz}
Let $\cc A$ be a ring. There is an inclusion-preserving
bijective correspondence between Serre subcategories of
$\coh\lcc$ and those of $\coh\rcc$ given by
   $$\cc S\longmapsto D\cc S.$$
The induced map $\cc O(\cc S)\mapsto\cc O(D\cc S)$ is an isomorphism
between the topologies, that is, the respective algebras of open sets,
of the Ziegler spectra of $\lcc$ and $\rcc$. Furthermore, the
duality $D$ induces dualities between the respective subcategories
$D:\cc S^{\text{\rm op}}\to D\cc S$ and
$D:(\coh\lcc/\vec\cc S)^{\text{\rm op}}\to\coh\rcc/\overrightarrow{D\cc S}$
as given by the following commutative diagram of abelian categories:
   $$\begin{CD}
      0@>>>\cc S@>>>\coh\lcc@>>>\coh\lcc/\vec\cc S@>>>0\\
      @.@VDVV@VDVV@VVDV\\
      0@>>>D\cc S@>>>\coh\rcc@>>>\coh\rcc/\overrightarrow{D\cc S}@>>>0\\
     \end{CD}$$
\end{thm}

\section{Grothendieck categories as quotient categories of
$(\fpl A,\text{\rm Ab})$}\label{reka}

In this section we give another representation of the Grothendieck
category $\cc C$ as a quotient category of $\rcc$. To begin, let
us prove the following.

\begin{prop}\label{pravda}
Let $\mc A$ be the category of right $\cc A$-modules with
$\cc A=\{P_i\}_{i\in I}$ a ring and let $\rcc$ be the respective category
of generalized right $\cc A$-modules. Then the category $\mc A$
is equivalent to the quotient category of $\rcc$ with respect to
the localizing subcategory $\ccc P=\{F\in\rcc\mid F(P)=0\text{ for all
$P\in\cc A$}\}$.
\end{prop}

\begin{proof}
For an arbitrary functor $F\in\rcc$ by $F(\cc A)$ denote a right
$\cc A$-module defined as follows. If $P\in\cc A$, then we put
$F(\cc A)(P)=F({}_{\cc A}P)$. It is directly checked that $F(\cc A)\in\mc A$.
>From Theorem~\ref{thm2.1} it follows that the functor
$\varPhi:\rcc\to\mc A$, $F\mapsto F(\cc A)$, defines an equivalence
of categories $\mc A$ and $\rcc/\kr\varPhi$. Clearly,
$\ccc P=\kr\varPhi$.
\end{proof}

Because there is a natural equivalence of functors
$\tena P-\iso({}_{\cc A}P,-)$, from Theorem~\ref{thm2.1} it also
follows that the quotient category
$\rcc/\ccc P$ is equivalent to the subcategory
$\Mod\ol{\cc A}=\{\tena M-\mid M\in\mc A\}$.

\begin{rem}{\rm
It easy to see that the subcategory
$$\ccc P=\{\kr(\tena\mu-)\mid\mu\text{ is a monomorphism in $\mc A$}\}.$$

It should also be remarked that $\rcc/\ccc P$-injective objects,
in view of Lemma~\ref{clinj}, are precisely objects $\tena E-$ where
$E_{\cc A}$ is an injective right $\cc A$-module.
}\end{rem}

Now suppose that $\cc C$ is a Grothendieck category with the family
of generators $\cc U$. As usual, let $\mc A$ be the category of right
$\cc A$-modules with $\cc A=\{h_U\}_{U\in\cc U}$. We are now in a
position to prove the following.

\begin{thm}\label{wi}
Every Grothendieck category $\cc C$ with the family of generators
$\cc U$ is equivalent to the quotient category of $\rcc$ with
respect to some localizing subcategory $\cc S$ of $\rcc$.
\end{thm}

\begin{proof}
By Theorem~\ref{main} there is a pair of functors
$(s,q)$, where $s:\cc C\to\mc A$ and $q:\mc A\to\cc C$,
defining $\cc C$ as a quotient category of $\mc A$. In turn,
by Proposition~\ref{pravda} there is a pair of functors
$(g,h)$, where $g:\mc A\to\rcc$ and $h:\rcc\to\mc A$,
defining $\mc A$ as a quotient category of $\rcc$. It thus
suffices to show that $gs$ is a fully faithful functor,
the functor $qh$ is exact and left adjoint
to $gs$. Indeed, the composition $gs$ of fully faithful functors
$g$ and $s$ is again a fully faithful functor and the
composition $qh$ of exact functors $q$ and $h$ is an exact functor.
The fact that $gs$ is right adjoint to $qh$ follows from the
following isomorphisms:
   $$\rqq F{gs(C)}\iso\homa {h(F)}{s(C)}\iso\qc {qh(F)}C.$$
Hence $\cc C$ is equivalent to the quotient category of $\rcc$
with respect to the localizing subcategory ${\cc S}=\kr(qh)$.
\end{proof}

\begin{cor}\cite[Theorem~2.3]{GG}
Every Grothendieck category $\cc C$ with a generator
$U$ is equivalent to the quotient category of $\rc$,
$A=\End_{\cc C}U$, with
respect to some localizing subcategory $\cc S$ of $\rc$.
\end{cor}

\begin{proof}
It follows from Corollary~\ref{PPG} and Theorem~\ref{wi}.
\end{proof}

The ring $\cc A$ of the module category $\mc A$ is said to be
{\em right coherent\/} if each object $P\in\cc A$ is coherent.
Suppose now $\cc C$ is a locally coherent, i.e.\
$\cc U\subseteq\coh\cc C$ and $\mc A$ the respective module
category with $\cc A=\{h_U\}_{U\in\cc U}$. One easily verifies
that $\cc A$ is right coherent.

\begin{thm}\label{tak}
Let $\cc C$ be a Grothendieck category with the family of generators
$\cc U$. Consider the following conditions:

$(1)$ $\cc C$ is locally coherent, i.e.\ $\cc U\subseteq\coh\cc C$.

$(2)$ The localizing subcategory $\cc S$ from the preceding
      Theorem is of finite type.

$(3)$ $\cc S$ is of prefinite type.

$(4)$ $\zg\cc C=\{\tena E-\mid E\text{ is $\cc C$-injective}\}$
      is closed in $\zg\rcc$.

Then conditions $(1)$, $(2)$ and $(3)$ are equivalent and $(1), (2), (3)$
imply $(4)$. If $\cc C$ is a locally finitely generated
category, then also $(4)$ implies $(1), (2), (3)$.

\end{thm}

\begin{proof}
Equivalence~$(2)\Leftrightarrow (3)$ follows from Theorem~\ref{fincoh}.

$(1)\Rightarrow (4)$. By assumption the ring $\cc A$ is right coherent.
Therefore the category of right $\cc A$-modules $\mc A$ is
locally coherent. By Theorem~\ref{fincoh} $\zg\cc C$ is closed
in $\zg(\mc A)$. Our assertion would be proved if we showed that
$\zg(\mc A)$ is a closed subset of $\zg\rcc$. In view of
Theorem~\ref{krut} and Theorem~\ref{hk} this is equivalent to
the localizing subcategory $\ccc P$ of $\rcc$ to be of finite type.
Thus we must show that $\ccc P=\vec\ccc S$ where
$\ccc S=\ccc P\cap\coh\rcc$ is a Serre subcategory of $\coh\rcc$.
Clearly that $\vec\ccc S\subseteq\ccc P$. Let us show the
inverse inclusion.

Let $F\in\ccc P$; then there exists an exact sequence
   $$0\lra{}F\lra\iota\tena M-\lra{\tn\mu-}\tena N-\lra{\tn\nu-}\tena L-\lra{}0,$$
where $\tena M-=E(F)$, $\tena N-=E(\text{coker\,}\iota)$, $\mu:M\to N$ is a
monomorphism. By hypothesis $\cc A$ is right coherent, and so
the exact sequence $\es {M\bl\mu}{N\bl\nu}L$ is a direct limit of exact
sequences $0\to M_i\bl {\mu_i}\to N_i\bl {\nu_i}\to L_i\to 0$ with
$M_i,N_i,L_i\in\fpr A$~\cite[Lemma~5.7]{Kr3}. If
$C_i=\kr(\tn{\mu_i}-)$, then $C_i\in\ccc S$. Consider a
commutative diagram
   \begin{equation}\label{jarrett}
   \begin{CD}
     0@>>>C_i@>{\rho_i}>>\tena {M_i}-@>{\tn{\mu_i}-}>>\tn{N_i}-\\
     @.@.@VV\tn{\alpha_i}-V@VV\tn{\beta_i}-V@.\\
     0@>>>F@>\iota>>\tn M-@>{\tn\mu-}>>\tn N-
   \end{CD}
   \end{equation}
with $\alpha_i$ and $\beta_i$ the canonical homomorphisms.
 Because $(\tn{\mu\alpha_i}-)\rho_i={\tn{(\beta_i\mu_i}-)\rho_i}=0$,
it follows that there exists $\gamma_i:C_i\to F$ such that
$\iota\gamma_i=\tn{(\alpha_i}-)\rho_i$. Similarly,
given $i\le j$, we can construct the commutative diagrams
  $$\begin{CD}
     0@>>>C_i@>{\rho_i}>>\tena{M_i}-@>{\tn{\mu_i}-}>>\tena{N_i}-\\
     @.@V{\gamma_{ij}}VV@VV\tn{\alpha_{ij}}-V@VV\tn{\beta_{ij}}-V@.\\
     0@>>>C_j@>{\rho_j}>>\tena{M_j}-@>{\tn{\mu_j}-}>>\tena{N_j}-.
    \end{CD}$$
It is directly verified that the family $\{C_i,\gamma_{ij}\}$
is direct. Now, taking a direct limit on the upper
row in diagram~\eqref{jarrett}, one obtains $F=\lp C_i$.
Thus $F\in\vec{\ccc S}$.

$(2)\Rightarrow (1)$. By Theorem~\ref{fincoh} the
$\cc S$-localization $\qq{(\tena{h_U}-)}S\iso U$ of the
$\rcc$-coherent object $\tena{h_U}-$ is $\cc C$-coherent.

$(4)\Rightarrow (2)$. Suppose $\cc C$ is a locally finitely
generated category and $\zg\cc C=\zg\rcc/\cc S$ is a closed
subset of $\zg\rcc$. By Theorem~\ref{krut} there is a Serre
subcategory $\cc P$ of $\coh\rcc$ such that
$\zg\cc C=\zg\rcc/\vec\cc P$. From~\cite[Proposition~2.7]{GG}
it follows that $\cc S=\vec\cc P$. Now our assertion follows
from Theorem~\ref{hk}.
\end{proof}

\begin{cor}\label{gage}
Let $\mc A$ be the module category with $\cc A=\{P_i\}_{i\in I}$.
Then the following assertions are equivalent:

$(1)$ The ring $\cc A$ is right coherent.

$(2)$ The localizing subcategory $\ccc P$ is of finite type.

$(3)$ The localizing subcategory $\ccc P$ is of prefinite type.

$(4)$ $\zg(\mc A)=\{\tena E-\mid E_{\cc A}
      \text{ is an injective right $\cc A$-module}\}$
      is closed in $\zg\rcc$.
\end{cor}

It is useful to have available the following criterion of coherence
of a ring.

\begin{prop}\label{kach}
For a ring $\cc A=\{P_i\}_{i\in I}$ the following are equivalent:

$(1)$ $\cc A$ is right coherent.

$(2)$ For any finitely presented left $\cc A$-module $M$ the right
      $\cc A$-module $M^*=\homa M{\cc A}$ is finitely presented.

$(3)$ For any finitely presented left $\cc A$-module $M$ the right
      $\cc A$-module $M^*=\homa M{\cc A}$ is finitely generated.

$(4)$ For any coherent object $\rcohh C$ the right $\cc A$-module $C(\cc A)$
      is finitely presented.

$(5)$ For any coherent object $\rcohh C$ the right $\cc A$-module $C(\cc A)$
      is finitely generated.
\end{prop}

\begin{proof}
By Proposition~\ref{pravda} the functor $\rcc\to\mc A$, $F\mapsto F(\cc A)$,
defines an equivalence of categories $\rcc/\ccc P$ and $\mc A$. By the
preceding Corollary $\cc A$ is right coherent \ifff $\ccc P$ is of finite
type (= of prefinite type). Because the family
$\{({}_{\cc A}M,-)\}_{M\in\fpl A}$ is a family
of generators of $\rcc$, our assertion immediately follows from
Theorem~\ref{pref} and Proposition~\ref{fin}.
\end{proof}

\begin{cor}
For a ring $\cc A=\{P_i\}_{i\in I}$ the following are equivalent:

$(1)$ $\ccc P$ is coexact.

$(2)$ For any finitely presented left $\cc A$-module $M$ the right
      $\cc A$-module $M^*=\homa M{\cc A}$ is projective.

In particular, $\ccc P$ is of finite type and coexact \ifff
for any $M\in\fpl A$ the module $M^*=\homa M{\cc A}$ is a finitely generated
projective right $\cc A$-module. In that case, the ring $\cc A$ is right
coherent.
\end{cor}

\begin{proof}
It follows from Proposition~\ref{osya}. The second part follows from
Proposition~\ref{fin} and Proposition~\ref{kach}.
\end{proof}

\begin{example}{\rm
Let $\cc C_{\bb Z}$ be the category of generalized abelian groups. Then
the localizing subcategory $\cc P^{\bb Z}$ is evidently of finite type
and coexact.
}\end{example}

Now we intend to describe localizing subcategories of prefinite type
in locally finitely generated Grothendieck categories in terms of
localizing subcategories of finite type of the category $\rcc$.
Let $\ccc S$ be a Serre subcategory $\ccc P\cap\coh\rcc$ of $\rcc$.

\begin{prop}\label{1994}
For a localizing subcategory $\cc S$ of $\mc A$, $\cc A=\{P_i\}_{i\in I}$,
the following assertions are equivalent:

$(1)$ $\cc S$ is of prefinite type.

$(2)$ There is a Serre subcategory $\cc T$ of $\coh\rcc$ such that
      $\cc T\supseteq\ccc S$ and $\vec\cc T(\cc A)=\cc S$.
\end{prop}

\begin{proof}
Denote by $\cc P$ a localizing subcategory of $\rcc$ such that
$\cc P\supseteq\ccc P$ and $\cc P(\cc A)=\cc S$ (see Proposition~\ref{dalee}).
Let $\cc T=\cc P\cap\coh\rcc$; then $\cc T$ is a Serre subcategory of
$\coh\rcc$. Clearly, $\cc T\supseteq\ccc S$. Obviously that
$\vec\cc T\subseteq\cc P$. Hence
$\vec\cc T(\cc A)\subseteq\cc S$. Let us show the inverse inclusion.

Let $\f F=\{\f F^P\}_{P\in\cc A}$ be a Gabriel topology corresponding to
$\cc S$. Our assertion would be proved, if we showed that $\f F$ has a basis
consisting of those ideals $\f a$ of $P$ such that $\f a=\f b(\cc A)$, where
$\f b$ is a coherent subobject of $\tena P-$ such that $(\tena P-)/\f b\in\cc T$.

So let $\f a\in\f F^P$. Consider the following exact sequence
   $$\ii{\kr(\tn\alpha-)}{\tena{\f a}-\bl{\tn\alpha-}}{\tena P-}.$$
Because $\alpha$ is a monomorphism, $\kr(\tn\alpha-)\in\ccc P$. Let
$\wt{\f a}=\im(\tn\alpha-)$; then $\wt{\f a}(\cc A)=\f a$ and
$(\tena P-)/\wt{\f a}\in\cc P$. Write $\wt{\f a}=\sum_{i\in I}\f a_i$ as a
direct union of finitely generated subobjects $\f a_i$ of $\wt{\f a}$. Because
each $\f a_i$ is a subobject of the $\rcc$-coherent object $\tena P-$, it
follows that $\f a_i$ is coherent. One has
   $$\qq PS=\qq{\f a}S=\qq{(\wt{\f a}(\cc A))}S=\sum\nolimits_{i\in I}\qq{(\f a_i(\cc A))}S.$$
By Theorem~\ref{pref} the object $\qq PS\in\fg(\mc A/\cc S)$, and so there is
a finite subset $J$ of $I$ such that $\qq PS=\sum_{i\in J}\qq{(\f a_i(\cc A))}S$.
Let $\f b=\sum_{i\in J}\f a_i$; then $(\tena P-)/\f b\in\cc P$ and since
$(\tena P-)/\f b$ is a coherent objects, one has $(\tena P-)/\f b\in\cc T$.
Hence $\f b(\cc A)=\f a$ as was to be proved.
\end{proof}

Now let us consider a locally finitely generated Grothendieck category $\cc C$.
By Theorem~\ref{wi} there is a localizing subcategory $\cc S$ of $\rcc$ such
that $\cc C$ is equivalent to $\rcc/\cc S$ with $\cc A=\{h_U\}_{U\in\cc U}$.
We consider this equivalence as
identification. If $\cc Q$ and $\cc P$ are localizing subcategories of $\rcc$,
by $\qq{\cc Q}P$ denote a subcategory of $\rcc/\cc P$ which consist of
$\{\qq QP\mid Q\in\cc Q\}$. Denote by $\cc L$ a Serre subcategory
$\cc S\cap\coh\rcc$ of $\coh\rcc$.

\begin{prop}
Let $\cc Q$ be a localizing subcategory of a locally finitely generated
Grothendieck category $\cc C$ with $\cc U$ a family of generators of $\cc C$.
Let $\cc A=\{h_U\}_{U\in\cc U}$ be a ring generated by $\cc U$. Then the
following assertions are equivalent:

$(1)$ $\cc Q$ is of prefinite type.

$(2)$ There is a Serre subcategory $\cc T$ of $\coh\rcc$ such that
      $\cc T\supseteq\cc L$ and $\vec\qq{\cc T}S=\cc Q$.
\end{prop}

\begin{proof}
By Theorem~\ref{main} there is a localizing subcategory $\cc P$ of
$\mc A$ such that $\cc C$ is equivalent to $\mc A/\cc P$ and by
Proposition~\ref{dalee} there is a localizing subcategory $\cc V$ of
$\mc A$ such that $\cc V\supseteq\cc P$ and $\cc V/\cc P=\cc Q$. Since
both $\cc P$ and $\cc Q$ are of prefinite type, it follows that $\cc V$
is of prefinite type. From the preceding Proposition it follows that
$\cc V=\vec\cc T(\cc A)$ for some localizing subcategory of finite type
$\vec\cc T$ of $\rcc$. Then
$\cc Q=\qq{\cc V}P=\qq{(\vec\cc T(\cc A))}P=\vec\cc T_{\cc S}$.
\end{proof}

\begin{question}{\rm
Is it true that if $\vec\cc T$ is a localizing subcategory of finite type
in $\rcc$ containing the subcategory $\vec\ccc S$, then the subcategory
$\vec\cc T(\cc A)$ is localizing and
of prefinite type in $\mc A$? If this was true, we could define the Ziegler
topology of an arbitrary locally finitely generated Grothendieck category.
}\end{question}

\section{Absolutely pure and flat modules}

In this section we sketch how classes of absolutely pure, flat
$\cc A$-modules may be studied with the help of some torsion/localization
functors in the category $\rcc$.

\begin{defs}{\rm
(1) Let $\cc C$ be a locally finitely presented Grothendieck
category. An object $C\in\cc C$
is said to be {\em absolutely pure\/} (or {\em $FP$-injective\/}) if
$\qq{\Ext^1}C(X,C)=0$ for any $X\in\fp\cc C$. An object $C\in\cc C$ is called
{\em $fp$-injective\/} if for any monomorphism $\mu:X\to Y$ in $\fp\cc C$
the morphism $\qc{\mu}C$ is an epimorphism. Evidently that every absolutely
pure object is $fp$-injective and every $fp$-injective finitely presented
object is absolutely pure. The ring $\cc A=\{P_i\}_{i\in I}$ is {\em right
absolutely pure\/} if every right $\cc A$-module $P\in\cc A$ is absolutely
pure.

(2) Let $\mc A$, $\cc A=\{P_i\}_{i\in I}$, be the category of right
$\cc A$-modules. A module $M\in\mc A$ is {\em flat\/} if the tensor functor
$\tena M-$ is exact. $M$ is called {\em $fp$-flat\/} if for any monomorphism
$\mu:{}_{\cc A}K\to{}_{\cc A}L$ in $\fpl A$ the morphism $\tena M\mu$ is a
monomorphism. Evidently that every flat module is $fp$-flat.
}\end{defs}

\begin{rem}{\rm
Some authors refer to absolutely pure objects as $fp$-injective.
Therefore, to make not terminological displacements, every time
we refer to somebody the reader should make more precise the terminology.
}\end{rem}

One easily verifies:

\begin{lem}\label{cave}
Let $\cc C$ be a locally finitely presented Grothendieck category.
The following statements are equivalent for an object $C\in\cc C$:

$(1)$ $C$ is absolutely pure.

$(2)$ Every exact sequence $\es C{C'}{C''}$ is pure.

$(3)$ There exists a pure-exact sequence $\es C{C'}{C''}$ with $C'$
      absolutely pure.
\end{lem}

We shall adhere to the following notation:
   \begin{gather*}
    \ccc P=\{F\in\rcc\mid F(\cc A)=0\}\\
    \ccc S=\{C\in\coh\rcc\mid C(\cc A)=0\}\\
    \qq{\cc S}A=\{C\in\coh\rcc\mid (C,\tena P-)=0\text{ for all $P\in\cc A$}\}.
   \end{gather*}
Subcategories ${}^{\cc A}\cc S$ and ${}_{\cc A}\cc S$ of $\coh\lcc$ are
defined similar to subcategories $\ccc S$ and $\qq{\cc S}A$ respectively.
By Theorem~\ref{hk} $\vec\ccc S$ and $\vec\qq{\cc S}A$
are localizing subcategories of finite type. By Corollary~\ref{gage}
$\ccc P=\vec\ccc S$ \ifff the ring $\cc A$ is right coherent.
If we consider presentations~\eqref{7.2}, it is easily seen that
$$\ccc S=\{\kr(\tn\mu-)\mid\mu:M\to N\text{ is a monomorphism in $\fpr A$}\}.$$
Similarly,
$$\qq{\cc S}A=\{\coker(\mu,-)\mid\mu:L\to K\text{ is a monomorphism in $\fpl A$}\}.$$

\begin{prop}\label{waits}
Let $K\in\mc A$; the the following assertions hold:

$(1)$ $K$ is absolutely pure \ifff the functor $\tena K-$ is $\ccc P$-torsionfree.

$(2)$ $K$ is $fp$-injective \ifff the functor $\tena K-$ is $\vec\ccc S$-torsionfree.

$(3)$ $K$ is $fp$-flat \ifff the functor $\tena K-$ is $\qq{\vec\cc S}A$-torsionfree.
\end{prop}

\begin{proof}
Adapt the proof for modules over a ring $\cc A=\{A\}$~\cite[Proposition~2.2]{GG2}.
\end{proof}

\begin{cor}
The set of indecomposable pure-injective $fp$-injective ($fp$-flat) modules
is closed in $\zg\rcc$.
\end{cor}

It is well-known~(see e.g.~\cite{St}) that for a ring $A$ a right $A$-module
$M$ is finitely presented (finitely generated) \ifff the natural map
$\ten M{(\prod_{i\in I}N_i)}\to\prod_{i\in I}(\ten M{N_i})$
is an isomorphism (epimorphism) for every family $\{N_i\}_{i\in I}$ of
right $A$-modules. This generalises to arbitrary module categories $\mc A$
as follows.

\begin{lem}\cite[Lemma~7.1]{Kr3}\label{plant}
Let $\cc A$ be a ring. For $M\in\mc A$ the following are
equivalent:

$(1)$ M is finitely presented (finitely generated).

$(2)$ The natural morphism
      $\tena M{(\prod_{i\in I}N_i)}\to\prod_{i\in I}(\tena M{N_i})$
      is an isomorphism (epimorphism) for every family $\{N_i\}_{i\in I}$
      in $\mc A$.

$(3)$ The natural morphism
      $\tena M{(\prod_{i\in I}P_i)}\to\prod_{i\in I}(\tena M{P_i})=\prod_{i\in I}M(P_i)$
      is an isomorphism (epimorphism) for every family $\{P_i\}_{i\in I}$
      in $\cc A$.
\end{lem}

Now the next result may be proved similar to~\cite[Proposition~2.3]{GG2}.

\begin{prop}\label{fripp}
Let $\cc A=\{P_i\}_{i\in I}$ be a ring; then:

$(1)$ For every family of right $\cc A$-modules $\{M_i\}_I$ the module
$\prod_I M_i$ is absolutely pure (respectively $fp$-injective, $fp$-flat)
\ifff every $M_i$ is absolutely pure (respectively $fp$-injective, $fp$-flat).

$(2)$ The direct limit $\lp M_i$ of $fp$-injective (respectively $fp$-flat)
right $\cc A$-modules $M_i$ is an $fp$-injective
(respectively $fp$-flat) module.
\end{prop}

Let us consider now a locally finitely presented Grothendieck category $\cc C$
with the family of generators $\cc U\subseteq\fp\cc C$. As usual, consider
the category of modules $\mc A$ with $\cc A=\{h_U\}_{U\in\cc U}$. By
Theorem~\ref{main} $\cc C$ is equivalent to the quotient category $\mc A/\cc S$.
Furthermore, by Proposition~\ref{fin} $\cc S$ is of finite type. By
Theorem~\ref{wi} there is a localizing subcategory $\cc P$ of $\rcc$ such that
$\cc C$ is equivalent to the quotient category $\rcc/\cc P$. Similar to the
category of modules, absolutely pure/$fp$-injective objects of $\cc C$
can be described in terms of torsion functors in $\rcc$. To begin, let us
prove the following.

\begin{prop}\label{brubek}
For an object $C\in\cc C$ the following assertions hold:

$(1)$ $C$ is absolutely pure \ifff it is absolutely pure as a right $\cc A$-module.

$(2)$ $C$ is $fp$-injective \ifff it is $fp$-injective as a right $\cc A$-module.
\end{prop}

\begin{proof}
(1). Let $C$ be a $\cc C$-absolutely pure object and $M\in\fpr A$. We must
show that $\Ext^1_{\cc A}(M,C)=0$. Equivalently, any short exact sequence
   $$\les{C\bl\alpha}XM$$
of right $\cc A$-modules splits. By Proposition~\ref{fin} $\qq MS\in\fp\cc C$.
By assumption, the morphism $\qq\alpha S$ splits, i.e.\
there exists $\beta:\qq XS\to C$ such that $\beta\qq\alpha S=1_C$. Then
$(\beta\lambda_X)\alpha=\beta\qq\alpha S=1_C$ where $\lambda_X$
is an \env S of $X$. So $\alpha$ splits.

Conversely, let $C$ be an absolutely pure right $\cc A$-module and let
   $$\epsilon:\les{C\bl\alpha}{E\bl\beta}X$$
be a $\cc C$-exact sequence with $E=E(C)$ and $X=E/C$. By assumption,
the short exact sequence
   $$\bar\epsilon:\les{C\bl\alpha}{E\bl{\bar\beta}}{\im\beta}$$
is pure-exact in $\mc A$. Clearly that $\qq{\bar\epsilon}S=\epsilon$.
>From Proposition~\ref{piano} it follows that $\bar\epsilon$ is a direct
limit of split exact sequences
   $$\bar\epsilon_i:\les{C_i}{E_i}{M_i}$$
in $\mc A$. Then $\epsilon$ is a direct limit of split exact sequences
$\epsilon_i=\qq{(\bar\epsilon_i)}S$. Thus $C$ is $\cc C$-absolutely pure.

(2). Suppose $C$ is an $fp$-injective object in $\cc C$ and $\mu:M\to N$
is a monomorphism in $\fpr A$. Since $\cc S$ is of finite type, the
morphism $\qq\mu S$ is a monomorphism in $\fp\cc C$. Consider the
commutative diagram
   $$\begin{CD}
      \homa NC@>(\mu,C)>>\homa MC\\
      @VVV@VVV\\
      \qc{\qq NS}C@>>(\qq\mu S,C)>\qc{\qq MS}C
     \end{CD}$$
where vertical arrows are isomorphisms. Because $(\qq\mu S,C)$
is an epimorphism, it follows that $(\mu,C)$ is an epimorphism.

Coversely, suppose $\mu:X\to Y$ is a monomorphism in $\fp\cc C$.
Then there is a monomorphism $\gamma:M\to N$ in $\fpr A$ such
that $\qq\gamma S=\mu$. Indeed, we can embed $\mu$ into the
commutative diagram in $\cc C$ with exact rows:
   $$\begin{CD}
      (-,\ps_{i=1}^nU_i)@>\psi>>(-,\ps_{j=1}^mU_j)@>>>X@>>>0\\
      @VVV@VVV@VV{\mu}V\\
      (-,\ps_{k=1}^sU_k)@>>\phi>(-,\ps_{l=1}^tU_l)@>>>Y@>>>0.
     \end{CD}$$
Because each $\ps h_{U_i}$ is \cl S and finitely generated projective
in $\mc A$, both $\coker\psi$ and $\coker\phi$ are finitely
presented right $\cc A$-modules. We put $M=\coker\psi$
and $N=\coker\phi$. There is a unique morphism
$\gamma:M\to N$. Since $\qq MS=X$ and $\qq NS=Y$, it follows
that $\qq\gamma S=\mu$. Consider the commutative diagram
   $$\begin{CD}
      \homa NC@>(\gamma,C)>>\homa MC\\
      @VVV@VVV\\
      \qc YC@>>(\mu,C)>\qc XC
     \end{CD}$$
where vertical arrows are isomorphisms. Because $(\gamma,C)$
is an epimorphism, it follows that $(\mu,C)$ is an epimorphism.
So $C$ is $fp$-injective in $\cc C$.
\end{proof}

\begin{cor}
The ring $\cc A=\{h_U\}_{U\in\cc U}$ is right absolutely pure \ifff
each $U\in\cc U$ is an absolutely pure object in $\cc C$.
\end{cor}

\begin{proof}
It suffices to note that each $h_U$ is \cl S (see Theorem~\ref{main}) and
then apply the preceding Proposition.
\end{proof}

Denote by $\cc T=\coh\cc C\cap\cc P$ and let $\qq tP$ and
$\qq tT$ be the torsion functors corresponding to the
localizing subcategories $\cc P$ and $\vec\cc T$ of $\rcc$.

\begin{prop}
Let $C\in\cc C$; the the following assertions hold:

$(1)$ $C$ is an absolutely pure object of $\cc C$ \ifff
      $\qq tP(\tena C-)=0$.

$(2)$ $C$ is an $fp$-injective object of $\cc C$ \ifff
      $\qq tT(\tena C-)=0$.
\end{prop}

\begin{proof}
(1). Let $C$ be absolutely pure. By the preceding Proposition
it is an absolutely pure right $\cc A$-module. Now let $E$ be an
injective envelope of $C$. Then $\tena C-$ is a subobject
of $\tena E-$. Because $\tena E-$ is \tr P, it follows that
$\tena C-$ is \tr P. Conversely, since $\cc P\supseteq\ccc P$,
our assertion follows from Propositions~\ref{waits} and~\ref{brubek}.

(2). Let $C$ be $fp$-injective and $T\in\cc T$. Consider exact
sequence~\eqref{7.2}
   $$\ii T{\tena M-\bl{\tn\mu-}}{\tena N-}$$
where $M$, $N\in\fpr A$. Because $0=\qq TP=\qq{T(\cc A)}S$, it
follows that the morphism $\qq\mu S$ is a monomorphism in
$\fp\cc C$. Consequently, the morphism $(\qq\mu S,C)$ is an
epimorphism, and so the morphism $(\mu,C)$ is an epimorphism
too. As $\tena C-$ is a $\coh\rcc$-injective
object, one has an exact sequence
   $$\pp{(\tena N-,\tena C-)\bl{(\tn\mu-,\tena C-)}}{(\tena M-,\tena C-)}{(T,\tena C-)}.$$
But $(\mu,C)$ is an epimorphism, hence $(T,\tena C-)=0$. So $\qq tT(\tena C-)=0$.
Because $\cc T\supseteq\ccc S$, the converse follows from the preceding
Proposition and Proposition~\ref{waits}.
\end{proof}

For a ring $A$ the Chase Theorem asserts that $A$ is left coherent \ifff
any direct product $\prod M_i$ of flat right $A$-modules $M_i$ is flat.
This generalizes to arbitrary module category as follows.

\begin{prop}[Chase]\label{page}
Let $\cc A=\{P_i\}_{i\in I}$ be a ring. Then the following are a equivalent:

$(1)$ $\cc A$ is left coherent.

$(2)$ Every product of flat right $\cc A$-modules is flat.

$(3)$ Every product $\prod_{j\in J}P_j$ of $P_j\in\cc A$ is a flat right
      $\cc A$-module for every set $J$.
\end{prop}

\begin{proof}
$(1)\Rightarrow (2)$: Let $\{M_j\}_{j\in J}$ be the family of flat right
$\cc A$-modules. By Proposition~\ref{fripp} the module $\prod_{j\in J}M_j$
is $fp$-flat. Let $\phi:K\to L$ be a monomorphism in $\mcl A$. As $\cc A$
is left coherent, it follows that $\phi=\lp\phi_i$ is a direct limit of
monomorphisms $\phi_i$ in $\fpl A$~\cite[Lemma~5.9]{Kr3}. Then the morphism
$\tn\phi{\prod M_j}=\lp(\tn{\phi_i}{\prod M_j})$ is a direct limit of
monomorphisms $\tn{\phi_i}{\prod M_j}$. So it is a monomorphism too.

$(2)\Rightarrow (3)$ is trivial.

$(3)\Rightarrow (1)$: Let $\qq{}AK$ be a finitely generated submodule of
finitely presented module $\qq{}AL$. For each index set $J$ we have a
commutative diagram
   $$\begin{CD}
      \tena K{\prod P_j}@>>>\tena L{\prod P_j}\\
      @V\phi_KVV@VV\phi_LV\\
      \prod K(P_j)@>>>\prod L(P_j)\\
     \end{CD}$$
where the horizontal arrows are monomorphisms. Since $\phi_L$ is a
monomorphism by Lemma~\ref{plant}, also $\phi_K$ is a monomorphism.
Thus $K$ is finitely presented by Lemma~\ref{plant}.
\end{proof}

In contrast to absolutely pure right $\cc A$-modules a class of flat right
$\cc A$-modules is realized in $\rcc$ as a class of those $\tena M-$ for
which $\qq tS(\tena M-)=0$ for some localizing subcategory $\cc S$ of
$\rcc$ \ifff $\cc A$ is left coherent.

\begin{thm}\label{jons}
For a ring $\cc A$ the following assertions are equivalent:

$(1)$ $\cc A$ is left coherent.

$(2)$ In $\rcc$ there is a localizing subcategory $\cc S$ such that any right
      $\cc A$-module $M$ is flat \ifff the functor $\tena M-$ is \tr S.

$(3)$ Every left $fp$-injective $\cc A$-module is absolutely pure.

$(4)$ Every right $fp$-flat $\cc A$-module is flat.

$(5)$ A direct limit of absolutely pure left $\cc A$-modules is absolutely pure.
\end{thm}

\begin{proof}
$(1)\Leftrightarrow (5)$: It follows from~\cite[Lemma~9.3]{Kr3}.

The rest is proved similar to~\cite[Theorem~2.4]{GG2}.
\end{proof}

\begin{thm}\label{bonem}
For a ring $\cc A$ the following conditions are equivalent:

$(1)$ $\cc A$ is right absolutely pure.

$(2)$ $\ccc S\subseteq\cc S_{\cc A}$.

$(3)$ ${}_{\cc A}\cc S\subseteq{}^{\cc A}\cc S$.

$(4)$ Every $fp$-flat right $\cc A$-module is $fp$-injective.

$(5)$ Every indecomposable pure-injective $fp$-flat right $\cc A$-module is
      $fp$-injective.

$(6)$ Every pure-injective $fp$-flat right $\cc A$-module is $fp$-injective.

$(7)$ Every $fp$-injective left $\cc A$-module is $fp$-flat.

$(8)$ Every indecomposable pure-injective $fp$-injective left $\cc A$-module
      is $fp$-flat.

$(9)$ Every pure-injective $fp$-injective left $\cc A$-module is $fp$-flat.
\end{thm}

\begin{proof}
Adapt the proof for modules over a ring (see~\cite[Theorem~2.5]{GG2}).
\end{proof}

\begin{example}{\rm
Let $\cc A=\{A\}$ be a ring and $\rc$ the category of generalized right
$A$-modules. Then the ring $\cc B=\{(M,-)\}_{M\in\lfp}$ is right
absolutely pure \ifff $A$ is (von~Neumann) regular.

Indeed, let $\cc B$ be right absolutely pure; then each $(K,-)$ with
$K\in\lfp$ is $\coh\rc$-injective. Therefore $(K,-)$ is isomorphic to
the object $\ten{K^*}-$ where $K^*=\Hom_A(K,A)$. Since every coherent
object $\rcoh C$ is a cokernel
   $$\pp{(K,-)\bl{(\alpha,-)}}{(L,-)}C$$
of $(\alpha,-)$ (see sequence~\eqref{7.1}), it follows that $C$ is
isomorphic to $\ten{(\coker{\alpha^*})}-$. So every $\rcoh C$ is
$\coh\rc$-injective, and hence $A$ is a regular ring~\cite[Theorem~4.4]{He}.
The converse follows from~\cite[Theorem~4.4]{He}.
}\end{example}

To conclude, we shall give a criterion of a duality for categories of finitely
presented left and right $\cc A$-modules. The ring $\cc A$ over which the
functor $\homa-{\cc A}$ constitutes a duality of indicated categories
one calls {\em weakly Quasi-Frobenius.} For the case $\cc A=\{A\}$ with $A$
a ring, we refer the reader to~\cite{GG2}.

\begin{thm}
For a ring $\cc A=\{P_i\}_{i\in I}$ the following assertions are equivalent:

$(1)$ $\cc A$ is weakly Quasi-Frobenius.

$(2)$ $\cc A$ is (left and right) absolutely pure and (left and right) coherent.

$(3)$ Classes of flat right $\cc A$-modules and absolutely pure right
      $\cc A$-modules coincide.

$(4)$ $\cc A$ is left absolutely pure and left coherent and any
      flat right $\cc A$-module is absolutely pure.

$(5)$ $\cc A$ is right absolutely pure and right coherent, and any
      absolutely pure right $\cc A$-module is flat.

And also assertions $(3')-(5')$, obtained from $(3)-(5)$ respectively
by sub\-sti\-tu\-ting the word ``right'' for ``left'' and vice versa.
\end{thm}

\begin{proof}
$(1)\Rightarrow (2)$: By assumption, given a finitely presented left
$\cc A$-module $M$, the right $\cc A$-module $M^*=\homa M{\cc A}$ is
finitely presented. From Proposition~\ref{kach} it follows that $\cc A$
is right coherent. Symmetrically, $\cc A$ is left coherent. Because the
functor $\homa-P$ where $P\in\cc A$ is exact both on $\fpl A$ and on $\fpr A$,
it follows that $P$ is an absolutely pure both left and right $\cc A$-module.
So $\cc A$ is an (two-sided) absolutely pure ring.

$(2)\Rightarrow (1)$: Since $\cc A$ is a (two-sided) coherent ring, from
Proposition~\ref{gage} it follows that $\ccc P=\vec\ccc S$, and hence
there is an equivalence of categories $\fpr A$ and $\coh\rcc/\vec\ccc S$.
Similarly, there is an equivalence of categories $\fpl A$ and
$\coh\lcc/\overrightarrow{{{}^{\cc A}\cc S}}$. In view of Theorem~\ref{bonem}
we have the following relations:
   \begin{gather*}
     \ccc S=D(\qq{}A\cc S)=D({}^{\cc A}\cc S)\\
     \qq{\cc S}A=D({}^{\cc A}\cc S)=D(\qq{}A\cc S).
   \end{gather*}
Now our assertion follows from Theorem~\ref{herz}.

$(2)\Rightarrow (3)$, $(2)\Rightarrow (4)$: Apply Theorems~\ref{jons}
and~\ref{bonem}.

$(3)\Rightarrow (5)$: It suffices to show that $\cc A$ is right coherent. To
see this, consider a direct system of absolutely pure right
$\cc A$-modules $\{M_i\}_{i\in I}$. Since each $M_i$, by assumption, is flat,
it follows that the module $\lp M_i$ is flat, and so it is absolutely pure.
Therefore $\cc A$ is right coherent by Theorem~\ref{jons}.

$(4)\Rightarrow (3)$: By Theorem~\ref{bonem} any absolutely pure right
$\cc A$-module is $fp$-flat, and hence flat by Theorem~\ref{jons}.

$(5)\Rightarrow (2)$: Since the ring $\cc A$ is right absolutely pure, the
module $\prod_J P_j$, where $P_j\in\cc A$ and $J$ is some set of indices,
is also absolutely pure, and therefore it is flat.
In view of Proposition~\ref{page} $\cc A$ is left coherent. By
Theorem~\ref{jons} any $fp$-injective right $R$-module is absolutely pure,
and hence flat. From Theorem~\ref{bonem} it follows that $\cc A$ is
left absolutely pure.
\end{proof}

\end{document}